\documentclass[11pt]{amsart}
\usepackage{amsmath,amsthm,amsfonts,amssymb}
\usepackage{ifthen}

 \provideboolean{showlabels}
 \setboolean{showlabels}{false}
\newcommand{\mylabel}[1]{\ifthenelse{\boolean{showlabels}}{{\tt{[{#1}]}}\label
{#1}}{\label{#1}}}

\theoremstyle{plain}

\newtheorem*{A}{Theorem A}
\newtheorem*{B}{Theorem B}
\newtheorem*{C}{Theorem C}
\newtheorem*{D}{Theorem D}
\newtheorem*{E}{Theorem E}

\newtheorem*{Q1}{Question 1}
\newtheorem*{Q2}{Question 2}

\newtheorem{theorem}[equation]{Theorem}
\newtheorem{lemma}[equation]{Lemma}
\newtheorem{proposition}[equation]{Proposition}
\newtheorem{corollary}[equation]{Corollary}

\newtheorem*{ia}{Inductive Assumption}

\newtheorem*{Co}{Condition X}
 \newtheorem{claim}{Claim}

\theoremstyle{definition}

\newtheorem{remark}[equation]{Remark}
\newtheorem{defn}[equation]{Definition}

\numberwithin{equation}{section}

\newcommand{\pn}{P^{\nu}}
\newcommand{\qn}{Q^{\nu}}

\newcommand{\an}{\alpha^{\nu}}
\newcommand{\bn}{\beta^{\nu}}
\newcommand{\gn}{G^{\nu}}
\newcommand{\sn}{S^{\nu}}
\newcommand{\kn}{K^{\nu}}
\newcommand{\en}{E^{\nu}}
\newcommand{\hn}{H^{\nu}}

\newcommand{\br}[2]{ \overline{#1}_{#2}}

\newcommand{\ZZ}{\mathbb{Z}}
\newcommand{\CC}{\mathbb{C}}
\newcommand{\F}{\mathcal{F}}
\newcommand{\V}{\mathcal{B}}
\newcommand{\SB}{\mathcal{S}}
\newcommand{\T}{\mathcal{T}}
\newcommand{\U}{\mathcal{U}}
\newcommand{\W}{\mathcal{W}}

\newcommand{\Z}{\operatorname{Z}}
\newcommand{\Syl}{\operatorname{Syl}}

\newcommand{\tin}{\text{\textrm{\textup{ in }}}}

\newcommand{\Aut}{\operatorname{Aut}}

\DeclareMathOperator{\Irr}{Irr}
\DeclareMathOperator{\Lin}{Lin}

\DeclareMathOperator{\Hom}{Hom}
\DeclareMathOperator{\Ker}{Ker}

\begin{document}

\title{ Extendible  characters and monomial groups of odd order  }

\author{Maria Loukaki} 

\address{School of Mathematics, Georgia Institute of Technology, 686 Cherry St NW, Atlanta, GA 30332, USA}
\email{loukaki@math.gatech.edu}

\thanks{ The author was 
 partially supported by NSF, grants DMS 96-00106 and DMS 99-70030}

\maketitle

\begin{abstract}
Let $G$ be a finite $p$-solvable group, where $p$ is an odd prime. 
We establish a connection between extendible irreducible 
characters of subgroups of $G$ that lie 
under monomial characters of $G$ and nilpotent subgroups of $G$. We also provide a way to 
get ``good'' extendible irreducible characters inside subgroups of $G$. As an 
application, we show that every normal subgroup $N$ of a finite monomial odd  $p, q$-group $G$, that has 
nilpotent length less than or equal to $3$, is monomial.
\end{abstract}

\section{Introduction}
A finite group $G$ is called monomial  ($M$-group)
if each of its complex  irreducible characters can be induced from some linear character  of some subgroup of $G$.
One of the outstanding  questions   in the theory of monomial groups   is 
\begin{Q1}
Is a normal subgroup $N$ of an $M$-group $G$ itself an $M$-group? 
\end{Q1}
It was shown by  Dornhoff  \cite{do} and  independently by Seitz \cite{se}  that the answer is yes when $N$ is
 a normal Hall subgroup of $G$, and conjectured that the answer is always yes.
For  $M$-groups of  even order, Dade \cite{da1} and van der Waall \cite{wa} showed separately
that the answer could be no.  They constructed an  example
of a monomial group of  order $7 \cdot  2^9$ which has a normal subgroup of index two that is not  monomial.
In their common example both $N$ and $G / N$ have even order.
So Question 1 remains open  when $N$ or $G/ N$ or both have odd order.  We remark, that there has been
 evidence   (see for example  \cite{isa3, isa5, lew2, na, pa} )  suggesting   that the answer to the above question
 is yes if $G$ is an odd $M$-group.

In \cite{lo} we  prove
\begin{A}
If $G$ is  a monomial group of order $p^a q^b$, where 
$p$ and $q$ are odd primes and $a, b$ are non--negative integers, 
then any normal subgroup  $N$ of  $G$ is again monomial.
\end{A}
That is, for monomial  odd $p^aq^b$-groups $G$  the answer to Question 1 is yes.

For the proof of  Theorem A     a special type of reductions was followed. These reductions are based on an
 observation of M. Isaacs, according to which the  Clifford theory  for abelian normal subgroups $L$ of a
group $G$ preserves monomiality  of characters (see exercise (6.11) in \cite{is}). Firstly we fix 
a monomial group $G$, a normal subgroup  $N$ of $G$ and an irreducible character
$\psi \in \Irr(N)$. This way we form the  triple $T =(G, N, \psi)$.     Now we 
apply Isaacs observation  to normal subgroups $L$
 of $G$ that are contained in $N$ and linear characters of $L$ that lie under $\psi$.
In particular, if $L$ is any normal subgroup of $G$ contained in $N$,  and $\lambda$ is a linear character of $L$
lying under $\psi$, then we may pass from $G$  to
the stabilizer $G(\lambda)$ of $\lambda$ in $G$, without loosing the monomiality of
 those irreducible characters of $G(\lambda)$
that lie  above $\lambda$. This way we get a new triple $T_1 = (G(\lambda), N(\lambda), \psi_{\lambda})$, that we call a
 {\em direct  linear reduction of $T$ } , 
where $\psi_{\lambda} \in \Irr(N(\lambda))$ is the $\lambda$-Clifford correspondent of $\psi$, and thus induces $\psi$. 
Clearly $G(\lambda)$ may not be a monomial group, but  every irreducible character of $G(\lambda)$ lying above 
$\lambda$ is still monomial.  Repeated applications of  
the same type of reductions  leads to a ``minimal''  triple $T' = (G', N', \psi')$, where   $ N' \unlhd G' \leq G$ and 
$\psi' \in \Irr(N')$ induces $\psi$ to  $N$, and where no more reductions can be performed to $T'$ (a
 more detailed analysis  on the triples and their reductions is given in Section 3 below). We call $T'$ a {\em linear 
limit } of $T$. If $Z(T')$ is the center of the induced character $(\psi')^{G'} $,  then  there is a 
unique linear  character $\zeta'$ of $Z(T')$  lying under $(\psi')^{G'}$. We call $\zeta'$ the {\em central character } 
 of the  triple $T'$  and $Z(T')$ the  { \em center }  of  $T'$.
Furthermore, Isaacs observation implies (see Proposition \ref{l4} below) that  
 every irreducible  character in $\Irr(G')$ that lies above $\zeta'$ is monomial. 
Also  the kernel $\Ker(\zeta')$ of $\zeta'$ is a  normal subgroup of $G'$, while $\Ker(\zeta') \leq \Ker(\psi')$.
The question partially answered  in \cite{lo} is

\begin{Q2}
Assume that  $N$ is a normal subgroup of an $M$-group  $G$. Let $\psi $  be any irreducible character of $N$.
Does there exists a linear limit $(G',  N', \psi')$ of $(G,  N, \psi)$ with 
 the quotient group  $N'/ (\Ker{\zeta'})$  nilpotent?
\end{Q2}

It is clear that a positive answer to Question 2  implies a positive answer to Question 1.
What we actually prove in \cite{lo} is that  Question 2 has a positive answer when $G$ is an odd $p, q$-group.

In this paper we publish  two of the main tools, Theorems B and D below, needed  for  the proof of Theorem A, 
that we think, are interesting  in themselves,
and provide an explanation of the approach we have used in \cite{lo}.  
As an easy consequence of these two theorems we prove that Question 2 has
 a positive answer when $G$ is an odd $p, q$-group 
and $N$ is a normal subgroup of nilpotent length $\leq 3$, (Theorems C and E).

For the general case we used in \cite{lo},
apart from Theorems B and D,  what we called there ``triangular sets''.
This is quite a complicated machinery. Fortunately E. C. Dade
came up with an easier correspondence than the one  the triangular sets provide, 
thus we are able to  prove the general case 
without their  use as we will see in a forthcoming paper.

When applying the linear reductions described above,
 we often reach a situation were $G$ satisfies the following
\begin{Co}
$P$ is a normal $p$-subgroup of $G$, for some odd prime $p$, such that
its center $Z(P)$ is maximal among the abelian $G$-invariant subgroups of $P$.
Furthermore,  $\zeta \in \Irr(Z(P))$ is a  $G$-invariant faithful
irreducible character of $Z(P)$,   and thus  $Z(P)$ is a cyclic central subgroup of $G$.
\end{Co}
In particular, suppose that $P$ is a normal subgroup of $G$ and $\alpha \in \Irr(P)$.  Let 
$(G', P', \alpha')$ be  a linear limit of $(G, P, \alpha)$, and assume  that $\zeta'$  is the center of $(G', P', \alpha')$. 
Then the groups $G'/ \Ker(\zeta')$ and $P'/ \Ker(\zeta')$ satisfy Condition X (see Proposition \ref{l1} below).
Assume further that  $N \unlhd G$ contains $P$ while $N/ P$ is a $q$-group for some prime $q \ne p$.
In  order to answer Question 2 in this special case, it would be enough to show that a $q$-Sylow 
subgroup $Q'$ of $G' \cap N'$ satisfies $[Q', P'] \leq \Ker(\zeta')$. This is actually true, and it
follows from the fact that  every irreducible character of $G'$ lying above $\zeta'$ is monomial. 
Theorem B handles this situation.
 \begin{B}
Assume that $G$ is  a finite $p$-solvable  group where $p$ is some odd prime.
 Let $P \unlhd G$ be a normal $p$-subgroup of $G$
that along with  $\zeta \in \Irr(Z(P))$ satisfies  Condition X.
Assume further that $S\unlhd G$ is a nilpotent normal $p'$-subgroup of $G$ and $\beta \in \Irr(S)$.
Let $\chi \in \Irr(G)$ be an irreducible monomial character of $G$, that lies above $\zeta \times \beta$
and satisfies $\chi(1)_{p'}=\beta(1)$. If $Q$ is any $p'$-subgroup of $G$ such that $PQ \unlhd G$, then $Q$
centralizes $P$.
\end{B}
Based on Theorem B  we actually show
\begin{C}
 Assume that $G$ is a finite monomial group. Assume further that $G$ has normal subgroups $M 
\leq N$ such that $M$ is nilpotent  with odd order and $N/ M$ is nilpotent.
Then  the answer to both Questions 1 and 2 above is yes.
\end{C}
Note that in Theorem C the group $G$ need not be a $p,q$-group, not even odd.

Now assume that $N$ has nilpotent length $3$. In particular, assume that 
$Q \unlhd M \unlhd N$ are all normal subgroups of $G$ with $Q$ being a $q$-group,
 $M/ Q$ a $p$-group  and $N/M$ being a $q$-group, for two odd primes $p \ne q$.
By induction we may assume that, after performing the necessary linear reductions, the group $M$ is nilpotent.
So the obstacle  this time is of the form $Q_1 \unlhd P  \times Q_1  \unlhd P   \ltimes Q  \unlhd G$,  where  now 
both $G, P$ and $G, Q_1$ satisfy Condition X and in addition, 
every irreducible character of $G$ lying above two specific $G$-invariant  characters 
$ \alpha \times \beta \in \Irr(P \times Q_1)$ is monomial.
 Again, in order to answer Question 2 we need to show that $Q$ centralizes $P$. 
But this time we can't so easily guarantee the existence of a monomial character of $G$ with the correct degree. Observe that 
according to Theorem B we need a monomial  character of $G$ whose degree has the   $q$-part 
equal to $\beta(1)$. If the character $\beta$ extends to $G$ then this problem is solved using some basic  $\pi$-theory.
Unfortunately, there is no reason   for $\beta$ to extend, but we can replace  him with  another  ``good'' one as  
the following theorem shows.
\begin{D}
Let $P$ be a $p$-subgroup, for some  odd prime $p$, of a finite group $G$.
Let $Q_1,  Q$ be   $q$-subgroups of  $G$, for some odd prime
$q\ne p$, with $Q_1 \leq  Q$.
Assume that $P$  normalizes $Q_1$,
while $Q$ normalizes the product  $P \cdot Q_1$.
Assume further that  $\beta$ is an irreducible character of $Q_1$.
 Then there exists an  irreducible character $\bn$ of $Q_1$
such that
\begin{align*}
P(\beta ) &= P(\bn), \\
Q(\beta) \leq Q(\bn) \, &\text{ and } N_Q(P(\beta)) \leq Q(\bn),  \\
\bn &\text{  extends to }  Q(\bn).
\end{align*}
\end{D}

Theorem  B  along with Theorem D,    enables us to prove:
\begin{E}
The answer to both Questions 1 and 2 is yes if $G$  is an
 odd monomial  $p, q$-group and $N$ has nilpotent length $3$.
\end{E}

 Section 2 below contains the proof of Theorem B   that is the key step  for the proof of Theorem C. The proof of Theorem C
can be found in Section 4,   while in Section 3 we go through the basic definitions and properties of linear limits 
and we state related  theorems, needed  for the proof of Theorems C and E. 
In sections 5  and 6 we  prove Theorems D and E, respectively. 
All the groups of this paper are assumed to be finite. In addition,  all the modules have finite dimension.
The notation and terminology  follows \cite{is}, with a  few exceptions. That is,  we write  $N_M(K)$ or $N (M \tin K)$ 
for the normalizer of $M$ in $K$, whenever $M, K$ are subgroups of a finite group $G$.  Also if $\phi \in \Irr(M)$ we 
denote by  $K(\phi)$ the stabilizer of $\phi$ in $K$. 
 In addition, we use the terminology of   \cite{da} when symplectic
 modules are concerned.
  So,  if $\F$ is any finite field
of characteristic $p$, and $G$ is  any  finite group, we say that  a finite--dimensional  $\F G$-module $\V$
is a  {\em symplectic } $\F G$-module if
 $\V $ carries a symplectic bilinear form $<\cdot, \cdot>$    that is invariant   by $G$.
 For any   $\F G$-submodule  $\SB$ of $\V$,  we denote by
  $\SB ^{\perp} :=  \{t \in \V |<\SB , t> = \{ 0\} \}$
the { \em perpendicular  }  $\F G$-submodule to  $\SB$.
 The $\F G$ submodule $\SB$ of $\V$ is called
 {\em isotropic }  if $\SB \leq \SB^{\perp}$, and it is called {\em self--perpendicular }
if $\SB = \SB^{\perp}$.
We say that
 $\V $ is    { \em  anisotropic }
if it contains no non--trivial isotropic
$\F G$-submodules. Furthermore, $\V$ is called  { \em hyperbolic }
if it contains some self--perpendicular $\F G $-submodule $\SB$.

{\bf Acknowledgment }
Most of the work of this paper is part of my thesis, done under the 
guidance  of my adviser E. C. Dade. 
I thank him for the enormous amount of hours 
he has spent on  this thesis, 
all the inspiring  discussions and his endless support.
I would also like to thank the Mathematics Department of the University of
 Illinois for its support.

\section{ Proof of Theorem B}
We begin with an  equivalent form of
Theorem 3.2  in \cite{da}.
\begin{theorem}\mylabel{m.t1}
Suppose that $\F$ is a finite field  of odd characteristic $p$,
 that $G$ is a finite $p$-solvable group,
that $H$ is a subgroup of $p$-power index in $G$,
 that $\V$ is an anisotropic  symplectic  $\F G$-module
 and that $\SB$ is an $\F G$-submodule of $\V$. Then the
$G$-invariant symplectic form on $\V$ restricts to a $G$-invariant
symplectic form on $\SB$. If $\SB$, with this form, restricts to a hyperbolic
symplectic $\F H$-module $\SB|_H$, then $\SB = 0$.
\end{theorem}
\begin{proof}
  Since $\V$ is symplectic and $\F G$-anisotropic, so is its
$\F G$-submodule $\SB$. Theorem 3.2 of \cite{da}, applied to $\SB$,
tells us that $\SB$ is $\F G$-hyperbolic
if $\SB|_H$ is $\F H$-hyperbolic. In that case $\SB$ is both
$\F G$-anisotropic and $\F G$-hyperbolic. So it must be $0$.
\end{proof}

%\begin{theorem}[Dade]\mylabel{td}
%Let $G$ be a $p$-solvable group, where $p $  is an odd prime, and  let $\chi \in \Irr{G}$ be an irreducible monomial
%character of $G$ whose degree is a power of $p$. If $N$ is a subnormal subgroup of $G$ and $\psi \in  \Irr(N)$
%is an irreducible constituent of $\chi_N$, then $\psi$ is monomial.
%In the same paper [C[C \cite{da}, Dade was able to prove  the following   result.
%\end{theorem}

We can now  prove  Theorem B

\begin{proof}
In view of Condition X, $Z(P)$  is a cyclic central subgroup of $G$,
and it is maximal among the abelian
subgroups of $P$ that are  normal in $G$. So
every characteristic abelian subgroup of $P$ is contained in $Z(P)$ and thus
 is cyclic.  Hence P. Hall's  theorem  (see Theorem 4.9 in \cite{go})
implies that either  $P$ is an abelian group or it
is the central product
\begin{subequations}\mylabel{m.e11}
\begin{equation}
P = T \odot Z(P),
\end{equation}
 where $T = \Omega_1(P)$ is an extra special $p$-group of exponent $p$,  and
 \begin{equation}
T \cap Z(P) = Z(T).
\end{equation}
\end{subequations}

In the case that $P = Z(P)$ is an abelian group, Theorem B holds trivially,
as $P = Z(P) \leq Z(G)$ is centralized by $G$.
Thus we may assume that $P > Z(P)$ and \eqref{m.e11} holds.

Since $\chi$ lies above $\zeta \times \beta \in \Irr(Z(P) \times S)$, Clifford's theorem implies
the existence of a unique irreducible character $\Psi $ of $G(\beta ) = G(\zeta \times \beta)$,
that  also lies above $\zeta \times \beta$  and induces $\chi $ in $G$. Furthermore, the  hypothesis that
$\chi(1)_{p'} = \beta(1)$, implies that $|G:G(\beta)|$ is a power of $p$, since
$|G:G(\beta)|$ divides $\chi(1)/ \beta(1)$.  So we get
\begin{equation}\mylabel{eex}
\Psi(1)_{p'} =  (\Psi^G(1))_{p'} = \chi(1)_{p'} = \beta(1).
\end{equation}
By hypothesis,  $\chi$ is a monomial character. Furthermore,  $Z(P) \times S $ is a nilpotent normal subgroup of $G$.
Hence we can apply Theorem 3.1 in \cite{pa}. We conclude that $\Psi \in \Irr(G(\zeta \times \beta  )) $ is also monomial.
Therefore there exists a subgroup $H$ of $G(\beta)$,
  and a linear character $\lambda \in \Lin(H)$ that
induces $\Psi = \lambda^{G(\beta)}$. Clearly $\lambda$ also induces $\chi = \lambda^G$.

The product $HS$ forms  a subgroup of $G$.
Furthermore,
\setcounter{claim}{0}
\begin{claim}\mylabel{m.st1}
$|G: HS| $ is a  power of  $p$, and $(\lambda^{HS}) |_S = \beta$.
\end{claim}
\begin{proof}
Clearly Clifford's theorem implies
 $$
\Psi |_S = m \cdot \beta,
$$
for some integer $m \geq 0$.  Hence $\deg (\Psi)= m \deg (\beta)$.
By \eqref{eex} we have  $\Psi(1)_{p'}= \beta(1)$.  Thus $m$ is a power of $p$.
As $ H \leq HS  \leq G(\beta) $, the induced character $\lambda^{HS}$
lies in $\Irr(HS)$ and induces $(\lambda^{HS})^{G(\beta)}  = \lambda^{G(\beta)} = \Psi$.
 So
$$
\deg(\lambda^{HS}) \cdot |G(\beta):HS| = \deg(\Psi)= m \deg(\beta).
$$
Clifford's theorem also implies that $\lambda^{HS}|_S = r \beta$,
 for some integer $r$.
As $ \deg(\lambda^{HS}) = |HS :H|= |S : H \cap S|$ we get that
both  $\deg(\lambda^{HS}|_S)$ and  $r$ are  $p'$-numbers.
But
$$
r \deg(\beta)  \cdot |G(\beta) :HS| = \deg(\lambda^{HS}) \cdot |G(\beta):HS| = m \deg(\beta),
$$
with $m$ a $p$-number.
Hence $r=1$, while  $|G(\beta):HS|$ is a power of  $p$. This, along with the fact that
$G(\beta)$ has $p$-power index in $G$,
 completes the proof of the claim.
\end{proof}

The  fact that $\lambda \in \Lin(H)$ induces irreducibly to  $G$ implies that
the center, $Z(G)$,
of $G$ is a subgroup  of $H$. This, along with the fact that $Z(P ) \leq Z(G)$, implies
\begin{equation}\mylabel{m.e6}
Z(P) \leq Z(G) \leq H.
\end{equation}

Let $Q$ be any $p'$-subgroup of $G$  that satisfies  $N:= P Q \unlhd G$.  In order to show
 that $Q$ centralizes $P$, we can, without loss, assume that $S$ is a subgroup of $Q$,
 or else we may  work with the $p'$-group $Q S$  that also satisfies $P (QS) = (PQ) S \unlhd G$.
Let $E := [P, Q]$. Then $E$ is a characteristic subgroup of
 $N$ and thus a normal subgroup of $G$.
 Furthermore $S$ centralizes $E$, since $E$ is a subgroup of $P$.
Even more, we have
\begin{claim}\mylabel{n.st2}
$E = [P, Q]$ is an abelian group.
\end{claim}
\begin{proof}
 Suppose not. Then $E$ is a non-abelian normal subgroup of $G$
contained in  $P = T \cdot Z(P)$.
 As $Z(P) \leq Z(G)$ (by \eqref{m.e6}), we have
$E = [P, Q] = [T, Q] \leq T$, where $T  = \Omega_1(P)$.
 Furthermore, $Z(E)$ is an abelian  normal subgroup of $G$,
contained in $T \leq P$. Hence $Z(E)$ is contained in
$T \cap Z(P) = Z(T)$. As $E$ is non-abelian and $Z(T)$ has order $p$,
we conclude that  $Z(E)= Z(T) \leq Z(P) \leq Z(G)$.
Therefore $E= [T, Q]$ is an extra special subgroup of $T$ of exponent $p$,
and its  center is central in $G$.
Hence the group $E$ satisfies condition (4.3a)
 in \cite{da}. In addition, $Q$ is a $p'$-subgroup of $G$ such that
$QE$ is normal in $G$ (as  $P =[P,Q]C_P(Q )$ and thus
$G = EN_G(Q)$). Since $P$ is a $p$-group,
 the commutator subgroup $[E, Q]=[[P,Q], Q]$ coincides with $E=[P,Q]$.
Hence  (4.3b) in \cite{da} holds with $Q$ here, in the place of $K$ there.

As  the  index of $HS$ in $G$ is a power of $p$, and $PQ=N$ is a normal
 subgroup of $G$, we conclude that $HS$ contains a $p'$-Hall subgroup of
$PQ$. Hence $HS$ contains a $P$-conjugate of $Q$. Therefore, we may
replace $H$ and $\lambda$  by some $P$-conjugates, and
assume that $HS$ contains $Q$. (Observe that because $P \leq G(\beta)$,
  the P-conjugate of $H$ is still a subgroup of $G(\beta)$
while the corresponding $P$-conjugate of $\lambda$ induces $\Psi$ in $G(\beta)$.)

The subgroup  $H \cap (E \times S)$ of $E \times S$ is equal to
$(H \cap E ) \times  (H \cap S)$, since $|E|$
and $|S|$ are relatively prime.  This implies that
$$
HS \cap (E \times S)=
(H \cap (E \times S)) \cdot S= (H \cap E) \times S.
$$
Hence
$H \cap E = HS  \cap E$.
 Thus  $H \cap E$ is a normal subgroup of $HS$.
Furthermore, the restriction $\lambda|_{H\cap E}$ of $\lambda$ to $H\cap E$
is a linear character of $H \cap E$ that is clearly $H$-invariant.
It is also $S$-invariant, as $S$ centralizes $E  \geq H \cap E$.
Hence $\lambda|_{H \cap E}$ is $HS$-invariant.  We conclude that the
restriction  of the irreducible character  $ \lambda^{HS} $
of $HS$ to $H \cap E$ is a multiple of the linear character
$\lambda|_{H \cap E}$.
Of course the irreducible character $\lambda^{HS}$  of $HS$ induces
 irreducibly to $\chi \in \Irr(G)$, and lies above a non-trivial character
of $Z(E)$ (as $Z(E) \leq Z(P)$ and $\zeta \in \Irr(Z(P))$ is  faithful).
Hence we can apply Lemma (4.4) and its Corollary   (4.8) of \cite{da},
using $HS$ here in the place of $H$ there,  and $\lambda^{HS}$ here
 in the place of $\phi$ there. We conclude that
$HS \cap E = H \cap E$ is a maximal abelian subgroup of $E$.

Let $\bar{P}:=P / Z(P)$. Then $\bar{P}$ is a symplectic $\ZZ_p G$-module, since it affords the 
$G$-invariant symplectic bilinear form $c$ defined as 
$c(\bar{x}, \bar{y})= \zeta([x, y])$, for all $x, y \in P$  (where $\bar{x}$ and $\bar{y}$ 
are the images of $x, y$ in $\bar{P}$). 
According to the hypotheses of the theorem,
 $Z(P)$ is the maximal abelian  $G$-invariant  subgroup of $P$.
 Hence $\bar{P}$ is an anisotropic
$\ZZ_p G$-module. If $\bar{E}$ is the image of $E$ in $\bar{P}$, i.e.,
$\bar{E} \cong E /Z(E)$, then $\bar{E}$ is a symplectic $\ZZ_p G$-submodule
of $\bar{P}$, as $E$ is normal in $G$. Furthermore, $\bar{E}$
 is $\ZZ_p HS$-hyperbolic as $HS \cap E$ is a maximal abelian $HS$-invariant
subgroup of $E$. Since  the index $[G: HS]$ is a power of $p$, Theorem
\ref{m.t1} forces $\bar{E}$ to be trivial.
Hence $E = Z(E)$ is abelian, and the claim follows.
\end{proof}

Now  $E = [P,Q]$ is an abelian subgroup of $P$  normal in $G$.
According to the hypotheses of the theorem,  $Z(P)$   is a   maximal
such subgroup of $P$.
Therefore  $1 \leq [P,Q] \leq Z(P) \leq Z(G)$.
So  $Q$ centralizes $[P,Q]$, which implies that
$[P,Q,Q]=1$ and thus $[P,Q]=[P,Q,Q]=1$.

This completes the proof of the theorem.
\end{proof}

An immediate consequence is
\setcounter{eqn}{0}
\begin{corollary}\mylabel{coo1}
Assume that a finite $p$-solvable group  $G$ satisfies Condition X  for its normal $p$-subgroup $P$
and the character $\zeta \in \Irr(Z(P))$.
Assume further that $\chi$ is a monomial character of $G$ that lies above $\zeta$ and its
 degree is a power of $p$.
If $Q$ is any $p'$-subgroup of $G$ so that $P Q \unlhd G$ then $Q$ centralizes $P$.
\end{corollary}

\section{Linear limits}

In this section we  give the basic definitions and properties of ordered triples and their linear limits.
A more detailed approach on the subject that actually contains most of the results
 that follow, can be found in  \cite{dl}

We denote  by $\mathfrak T$ the family of all ordered triples $(G, N, \psi)$, where $G$ is
a finite group, $N$ is a normal subgroup of $G$ and $\psi$ is an irreducible character of $N$. 
Until the end of the section we fix  an element    $T= (G, N, \psi)$ of $\mathfrak T$. 
We define the {\em center}  $Z(T)$ of $T$ to be the center of the induced character $\psi^G$, and 
the {\em central character } $\zeta^{(T)}$ to be the unique  linear character of $Z(T)$ lying under $\psi^G$. 
Then $Z(T)$ is a normal subgroup of $G$ contained in $N$,
 while $\zeta^{(T)}$ is a $G$-invariant linear 
character of $Z(T)$. So the restrictions of both $\psi$ and $\psi^G$ to $Z(T)$ are  multiples
 of $\zeta^{(T)}$.
Even more, $Z(T)$ is fully characterized by  Proposition 2.3 in \cite{dl}, that we partly restate here.
\begin{proposition}\mylabel{char.z}
The center $Z(T)$  is the largest normal subgroup $L$ of $G$ contained in $N$  such that $\psi \in \Irr(N)$ lies 
over some $G$-invariant linear character $\lambda$ of $L$. Any
other such $L$ is a subgroup of $Z(T)$, and the corresponding $\lambda$ is teh restriction of 
$\zeta^{(T)}$ 
to $L$, while the restriction of $\psi$ to $L$  equals $\psi(1) \lambda$. 
\end{proposition}
We also define the {\em kernel } $\Ker(T)$ to be the kernel of $\psi^G$. So $\Ker(T) $ is 
$\text{Core}_G(\Ker(\psi))$, and thus
  is contained in $\Ker(\psi)$. In addition, $\Ker(T) = \Ker(\zeta^{(T)})$ (see the second section in \cite{dl} for a 
detailed analysis of character triples).  

A {\em subtriple } of $T$ is any triple $T' = (G', N', \psi')$  contained in $\mathfrak T$, with $G'$ being a subgroup of $G$,
 and $\psi' $ being any irreducible character of $N'= G' \cap N$ lying under $\psi$. We write  $T' \leq T$
to denote that $T'$ is a subtriple of $T$. 
 If  there exists some normal subgroup $L$ of $G $ 
contained in $N$ and a linear character $\lambda \in \Lin(L)$ lying under $\psi$, then the stabilizer  $G'= G(\lambda)$ 
of $\lambda$ in $G$ is a subgroup of $G$, while $N' = N \cap G'$  is the stabilizer $N' = N(\lambda)$ of $\lambda$ in $N$. 
Furthermore,  the $\lambda$-Clifford correspondent  $\psi'$ of $\psi$ in $N'$, 
 is the unique irreducible character of $N'$ lying  above $\lambda$  and inducing  $\psi$. 
So  the triple $T(\lambda)= (G', N', \psi')$  is a subtriple of $T$.
By  a {\em direct linear reduction }  of $T$
we mean any subtriple $T' \leq T$ 
 of the form $T' = T(\lambda)$,
for some $L$ and $\lambda$ satisfying the above conditions.
A direct linear reduction $T'$ is called {\em proper} if $T' \ne T$. 
If  the only possible direct linear reduction of $T$ is  $T$ itself, then  
$T$ is called { \em linearly irreducible}.  
 A subtriple $T'$ of $T$ is called 
a {\em linear reduction }  of  $T$ if there is some finite chain $T_0 , T_1, \dots, T_n$ of subtriples  $T_i  \leq T$, starting with
$T_0 = T$ and ending with $T_n = T'$, such that each $T_i$ is a direct linear reduction of $T_{i-1}$ , for all 
$i=1, 2, \dots, n$.  
A {\em linear limit} of $T$ is a linearly  irreducible linear reduction of $T$, i.e., it  is a linear reduction of $T$ that 
has no proper direct linear reductions.

If $T'= T(\lambda)$ is a direct linear reduction of $T$ then without loss (see Proposition 2.18 in \cite{dl}) we may assume that 
$\lambda$ is a linear character of some normal subgroup $L$ of $G$ satisfying  $Z(T) \leq L \leq N$. In addition,  $\lambda$ 
not only lies under $\psi$ but also lies above  the central character $\zeta^{(T)}$ of $T$. 
Furthermore, for any linear reduction $T'$ of $T$ we have $Z(T) \leq Z(T') $ while $\zeta^{(T)} $ is the restriction of 
$\zeta^{(T')}$ to $Z(T)$, see equation 2.24 in \cite{dl}. In addition $\Ker(T) = \Ker(\zeta^{(T)}) \leq  \Ker(\zeta^{(T')}) = \Ker(T')$.

The first remarks follow easily from the above definitions.
\begin{remark}\mylabel{l00}
If  $T'= (G', N', \psi')$ is a linear reduction of $T= (G, N, \psi)$,   then
any linear limit $T''$  of $T'$ is also a linear limit of $T$. In addition, 
$Z(T'' ) \geq Z(T') \geq Z(T)$ and $\Ker(T'' ) \geq \Ker(T') \geq \Ker(T)$.
\end{remark}

\begin{remark}\mylabel{l0}
Assume that $H$ is a subgroup of $G$ with $N \leq H \leq G$. Let $S$ be the ordered triple $S= (H, N, \psi)$.
If $T'= (G', N', \psi')$ is a linear limit of $T$ then $S'= (G'\cap H, N', \psi')$ is a linear reduction of $S$ with $Z(T') \leq Z(S')$.
So we can get a linear limit $S''= (H'', N'', \psi'')$ of $S'$,  and thus of $S$,
 with  $N'' \leq N'$ and $Z(T') \leq Z(S') \leq Z(S'')$. 
Furthermore, the central character $\zeta^{(T')}$ of $T'$ is the restriction of  $\zeta^{(S')}$, and thus of $\zeta^{(S'')}$. 
\end{remark}

\begin{remark}\mylabel{l01}
Let $T'= (G', N', \psi')$  be a  linear limit of $T=(G, N, \psi)$. If
 $B$ is a subgroup of $G$ that centralizes $N$,  (that is  $B \leq C_G(N)$) then $B \leq G'$.
\end{remark}

Now assume that $N \leq H  \leq G $, while the irreducible character $\theta  \in \Irr(H)$ of $H$ lies above 
 $\psi \in \Irr(N)$. If $T' = (G', N', \psi')$ is a linear reduction of $T$, 
  then we form the group 
$H' = H \cap G'$ that we call the  {\em $T'$-reduction  } of $H$.
Since $T'$ is a linear reduction of $T$, there is a  chain of linear subtriples $T_0, T_1, \dots, T_n$  of $T$, starting with $T_0 = T$ and
 ending with $T_n = T'$ such that   $T_i$ is a direct linear reduction
 of $T_{i-1}$, for all $i =1, \dots, n$.  So $T_i = T_{i-1}(\lambda_i)$, where   $\lambda_i$  is a linear character of some 
normal subgroup $L_i$ of $G_{i-1}$ with $Z(T_{i-1}) \leq L_i \leq N_{i-1}$. Furthermore, $\psi_i \in \Irr(N_i)$ is 
the $\lambda_i $-Clifford correspondent of $\psi_{i-1} \in \Irr(N_{i-1})$.
Because $\theta \in \Irr(H) $ lies above $\psi$,    there is also a finite chain 
of irreducible characters $\theta_i$ of  $H_i = G _i \cap H$,  
 starting with $\theta_0 = \theta$ and ending with $\theta_n = \theta'$, such that  
$\theta_i$ is the unique $\lambda_i$ -Clifford correspondent of $\theta_{i-1}$, 
for all $i=0, 1, \dots, n$.  So $\theta_i$ lies above $\psi_i$, 
for all $i=0, 1, \dots, n$. 
We call  $\theta'$ the {\em $T'$-reduction } of $\theta$.
It is clear from the definition of $\theta'$ that it lies above $\psi'$ and induces $\theta$ 
in $H$.  
Actually $\theta'$ is the unique irreducible character of $H'$ 
with these properties by 
\begin{remark}\mylabel{r10}
Assume that $N \leq H\leq G$, and let $T'= (G', N', \psi')$ be any  linear 
reduction   of $T$. 
If  $\phi$ is an irreducible character of $H' = H \cap G'$ that lies above  $\psi'$
then $\phi$ is the  $T'$-reduction of  the irreducible character $\phi^H$ of $H$.
\end{remark}

\begin{proof}
According to \cite{dl} induction is a bijection between the irreducible characters of $H' $ 
lying above 
$\psi'$ and those of $H$ lying above $\psi$.  Hence  $\phi^H$  is an irreducible character of $H$ 
and its $T'$-reduction  is $ \phi$.
\end{proof}

If,  in addition, $H$ is a normal subgroup of $G$,  then  we clearly have 
\begin{remark}\mylabel{l8}
Assume that $N\leq H $ are  normal subgroups of $G$, 
 while $\theta \in \Irr(H | \psi)$.
So we can form the triple $S = (G, H, \theta)$. Then $Z(S) \geq Z(T)$ and $\Ker(S) \geq \Ker(T)$. 
If    $T' = (G', N', \psi')$ is     a linear limit of $T$,  $H' = H \cap G'$ and 
$\theta' \in \Irr(H')$   is the $T'$-reduction  of $\theta$,  then 
the triple $S' = (G', H', \theta')$ is a linear reduction of $S = (G, H, \theta)$. 
\end{remark}

The following is Proposition 2.21 in \cite{dl}.
\begin{proposition}\mylabel{l1}
Let $T' = (G', N', \psi')$ be a linear limit of $T$. Then $Z(T')/\Ker(T')$ is the center
of $N'/\Ker(T')$ and  is a cyclic group, which  affords a faithful $G'/\Ker(T')$-invariant linear 
character that inflates to $\zeta^{ (T') }
\in \Irr(Z(T')) $. Furthermore, $Z(T')/\Ker(T')$ is maximal among the abelian normal subgroups of  $G'/\Ker(T')$ contained in
 $N'/\Ker(T')$.
\end{proposition}

If $N$ is a nilpotent group then we can easily see
\begin{remark}\mylabel{l2}
Assume that $N$ is a nilpotent group. Let $p_1, \dots, p_n$ be the distinct primes dividing $|N|$.
Then   $N= N_{1} \times \cdots \times N_{n}$, where $N_{i}$ is the $p_i$-Sylow
subgroup of $N$, for each  $i=1, \dots, n$. Let   $\psi = \psi_1 \times \cdots \times \psi_n$ with
$\psi_i \in \Irr(N_i)$ for $i=1, \dots, n$, 
be the corresponding factorization of $\psi$. We write $T_i $ for the triples
  $T_i = (G,  N_i, \psi_i)$, for all $i=1, \dots, n$. If
 $T_i'= (G_i',N_i', \psi_i')$ is a linear limit of $T_i$ for all such $i$,
 then the group $G' = G_1'  \cap \cdots \cap  G_n'$, its normal subgroup 
$N' = N_1' \times  \cdots \times N_n'$ and the irreducible character $\psi' = \psi_1' \times 
\cdots \times \psi_n'$ of $N'$ form a linear limit $T' = (G', N', \psi')$ of $T$.
Furthermore, if $Z_i'$ is the center $Z(T_i')$ of $T_i'$, and $\zeta_i' \in \Irr(Z_i')$ the
 corresponding central character of $T_i'$, for all $i=1, \dots, n$,  then $Z' = Z_1' \times \cdots \times Z_n'$
is the center $Z(T')$  of $T'$,  while $\zeta' = \zeta_1' \times \cdots \times \zeta_n'$  is the central character of $T'$.
In addition,  the kernel of $T'$ satisfies, $\Ker(T') =  \Ker(T_1') \times \cdots \times  \Ker(T_n')$.
\end{remark}

Assume now that $N / Z(T)$ is an abelian group. Then we can introduce  an alternating  bilinear form  (see section 5 in \cite{dl})
$c $ from $ N/ Z(T) \times N / Z(T) $ to  
the multiplicative group $\CC^{\times}$ of complex numbers, defined as 
\begin{equation}\mylabel{form}
c (\bar{x}, \bar{y}) = \zeta^{(T)}([x, y]), 
\end{equation}
for all elements  $\bar{x}, \bar{y}$  of $N/ \Z(T)$, where $x$ and $y$ are pre-images of $\bar{x}$ and $\bar{y}$ respectively, 
in $N$. The action of $G$ on $N$ via conjugation  makes $c$  a  $G/N$-invariant bilinear form. 
As in  the introduction we define the {\em  perpendicular subgroup} 
 $B^{\perp}$ for any subgroup $B \leq N/ Z(T)$, to be 
\begin{equation}\mylabel{perp}
B^{\perp} = \{  \bar{x} \in N/ Z(T)  \bigm | c (B, \bar{x}) = 1 \}. 
\end{equation} 
 We also call any subgroup $B$ of $N/ Z(T)$ { \em isotropic } if 
$c(B, B)= 1$, and we say that $N/Z(T)$ is { \em $G/N$-anisotropic}  if $1$ is the only $G/N$-invariant isotropic subgroup of $N/Z(T)$.
So the form $c$ is non-singular if and only if  the perpendicular subgroup of $N / Z(T)$ equals $1$, and  in this case  
$N/Z(T)$ becomes a symplectic $G/N$-group.
Because the perpendicular subgroup of $N/Z(T)$ is a $G/N$-invariant isotropic subgroup of $N/Z(T)$, 
 if $N/Z(T)$ is anisotropic as a $G/N$-group then it is also symplectic. 
Also,  Proposition 5.2 and Proposition 5.8 in \cite{dl}  imply 
\begin{proposition}\mylabel{aniso}
If $N/ Z(T)$ is a nilpotent group and $T$ is linearly irreducible,  then $N/ Z(T)$ is an abelian anisotropic $G/N$-group.
\end{proposition}
This,  along with  Proposition 5.9 in \cite{dl},  implies
\begin{proposition}\mylabel{l1.5}
Assume that $T'= (G', N', \psi')$ is  a linear limit of $T$, where $N'/ Z(T')$ is a nilpotent group.
Then $\psi'$ vanishes on $N' - Z(T')$, and is a multiple of $\zeta^{(T')}$ on $Z(T')$. Hence $G'(\psi')= G'$.
\end{proposition}

Assume now that the abelian factor group $N/ Z(T)$ is a symplectic $G/ N$-group, i.e., assume that 
the form $c$ defined above is non-singular. 
(As we noted above,  if $T$ is linearly irreducible then $N/ Z(T)$ is a symplectic $G/N$-group.)
Then the following proposition suggests another  way to look at  linear  limits of $T$.
\begin{proposition}\mylabel{l6}
Assume that  $N/ Z(T)$ is an abelian  and symplectic $G/ N$-group.  Assume further that 
 $T'= (G', N', \psi')$  is a linear limit  of $T$.  Then  $G = G' N $  and $G' \cap N = N'$.
Thus inclusion  $G' \hookrightarrow G$  induces an isomorphism of 
$G'/ N'$ onto $G/ N$. This  isomorphism 
turns $N'/ Z(T')$  into a symplectic $G/ N$-group.  In addition,  $N'/ Z(T')$ is naturally 
isomorphic,  as a symplectic $G/ N$-group,  to the factor group $L ^{\perp} / L$, where 
$L$ is maximal among the  $G/ N$-invariant isotropic subgroups  of $N / Z(T)$,  and $L^{\perp}$ 
is the perpendicular subgroup  to $L$ with respect to the bilinear form $c$. Furthermore, $L = Z(T') / Z(T)$ 
while $L^{\perp} = N' / Z(T)$.
\end{proposition}

\begin{proof}
This is Proposition 5.23 in \cite{dl}, combined with the earlier results of  Propositions 5.18  and 5.22 in \cite{dl}.
\end{proof}

\begin{proposition}\mylabel{l7}
Assume that $T'= (G', N', \psi')$ is  a   linear reduction  of $T= (G, N, \psi)$.
If $Z' = Z(T')$ is the center of $T'$ and $\zeta'$ the central character of $T'$, then
  $G' = G(\zeta')$ and $\psi'$ is  the unique irreducible character of $N'= N(\zeta')$ lying above $\zeta'$ 
and  inducing $\psi$.  If in addition $T'$ is a linear limit of $T$ and 
$N'/Z(T')$ is nilpotent then  $\psi'$ is the  only  character in $\Irr(N' | \zeta')$. 
\end{proposition}

\begin{proof}
Let $T(\lambda)= T_1= (G_1, N_1, \psi_1)$ be a direct linear reduction of $T$, where $\lambda$ is a linear character 
of a normal subgroup $L$ of $G$ contained in $N$. Then $L$ is contained in  the center  $Z(T_1)$ of $T_1$,
 see Proposition 2.3 in \cite{dl}. Hence
 $ L \leq Z(T_1)  \leq N$.  Furthermore, the same proposition 
in  \cite{dl} implies that $\lambda$  is  a restriction of
 the $G_1$-invariant linear character $\zeta^{(T_1)}$ of $Z(T_1)$.
Since $L \unlhd G$ we have  
$$
G(\zeta^{(T_1)})  \leq G(\lambda) = G_1 = G_1 (\zeta^{(T_1)}) \leq G(\zeta^{(T_1)}).
$$
Therefore, $G_1 = G(\zeta^{(T_1)})$ which in turn implies that $N_1 = N(\zeta^{(T_1)})$.  Furthermore, since 
 $\psi_1$ is the $\lambda$-Clifford correspondent of $\psi$, while 
$\zeta^{(T_1)}$ is a $G_1$-invariant irreducible character  of  $Z(T_1) \leq N_1$ lying above
 $\lambda$ and under $\psi_1$,
  we conclude that $\psi_1$ is the unique irreducible character of $N_1$ that lies 
above $\zeta^{(T_1)}$ and induces $\psi$.  Hence the first part of the proposition holds for a direct linear reduction.
A linear limit $T'$ of $T$ is a series of direct linear reductions with starting 
triple $T= T_0$  and ending triple  $T'$. Furthermore, $Z(T) = Z(T_0) \leq
Z(T_1) \leq \cdots \leq Z(T_n)= Z(T')$.
Hence the first part of the proposition follows.

Assume now that $T'$ is a linear limit of $T$ while 
 $N'/Z(T') $ is a nilpotent group. Then Proposition \ref{l1.5} implies that $\psi'$ is fully ramified  with respect to $N'/ Z(T')$. 
Hence, see Lemma 2.6 in \cite{isol},  the unique irreducible character  of $N'$ lying above $\zeta^{(T')}$ is  $\psi'$. 
This completes the proof of the proposition.
\end{proof}

We conclude  this section by  proving that linear limits preserve monomial characters
(see Proposition \ref{l4} below).
Its proof is based   on the following lemma, that is  actually
the exercise (6.11) in \cite{is}.
\begin{lemma}\mylabel{l3}
Let $B$  be a  normal subgroup of a finite group $G$ and $\gamma $   be a linear character of $B$.
Assume further that  $\chi \in \Irr(G|\gamma)$ is an irreducible  character of $G$ lying above
$\gamma$.
If $\chi_{\gamma} \in \Irr(G(\gamma))$ is the $\gamma$-Clifford correspondent
of $\chi$ in the stabilizer $G(\gamma)$ of $\gamma$ in $G$, then
$\chi$ is monomial if and only if
$\chi_{\gamma}$ is  monomial.
\end{lemma}

\begin{proof}
It is clear that if $\chi_{\gamma}$ is monomial then $\chi$ is monomial, as
$\chi_{\gamma}$ induces $\chi$ in $G$.

So we assume that $\chi$ is a monomial character, and we will show that
$\chi_{\gamma}$  is also monomial.
Let $K= \Ker(\chi)$.  Of course  $K \unlhd G$.
It is clear that $\chi_{\gamma}$ is monomial if and only if the irreducible
character  $\chi_{\gamma}/ K$ of the factor group $G(\gamma)/K$ that
inflates to $\chi_{\gamma}$, is monomial.
Hence it  suffices to prove the lemma
in the case of a faithful irreducible character $\chi$, as we can pass to
the quotient  groups $G/K$ and $(BK)/ K$.
So in the rest of the proof we assume that $K=1$.

Clifford's Theorem implies that the restriction  $\chi|_B$ of $\chi$ to
$B$  is a sum of $G$-conjugates of
$\gamma$. Thus  $1=\Ker(\chi|_B)= \bigcap_{s \in
G/G(\gamma)}(\Ker(\gamma^s))$.
But the derived group $[B, B]$  of $B$ is contained in the kernel of
$\gamma^s$ for every $s\in G$, as $\gamma $ is linear.
Thus $[B, B] \leq \Ker(\chi|_B)=1$. So $B$ is abelian.

We can now  follow the hint of problem 6.11 in \cite{is}.
As $\chi$ is monomial,  there exists $H \leq G$ and $\lambda \in \Lin(H)$
with $\chi= \lambda^G$. Thus the irreducible character $\lambda^{HB}$ of
$HB$ lies above a $G$-conjugate  $\gamma^s$ of $\gamma$, where $s \in
G $. As the $G$-conjugate  $\lambda^{s^{-1}} \in
\Lin(H^{s^{-1}})$  of $\lambda$ also induces $\chi$,
we can replace $H$ by $H^{s^{-1}}$ and $\lambda$ by $\lambda^{s^{-1}}$.
This way $\lambda^{HB} $ is replaced
by $(\lambda^{s^{-1}})^{H^{s^{-1}} B}=(\lambda^{HB})^{s^{-1}}$, which lies above $\gamma$.

According to
  Mackey's Theorem
\begin{equation}\mylabel{lim.e1}
\lambda^{HB}|_B = (\lambda|_{H\cap B}) ^B.
\end{equation}
As $B$ is abelian, the right hand side of
 \eqref{lim.e1} equals the sum of $|B:H \cap B|$
distinct character extensions of $\lambda|_{H \cap B}$ to $B$,
each one appearing with multiplicity one.
Thus every irreducible  constituent of
$\lambda^{HB}|_B $ appears with multiplicity one.
This, along with Clifford's  Theorem, (as
 $\lambda^{HB}$ lies above $\gamma$),
implies that
$$
\lambda^{HB}|_B= e\cdot\sum_{s \in  S}
 \gamma^s=\sum_{s \in S} \gamma^s ,
$$
where $S$ is a family of representatives for the cosets $H(\gamma)Bs$ of
$H(\gamma) B = (HB)(\gamma)$ in $HB$, and
$e$ is a positive integer.
Furthermore, Clifford's Theorem implies the existence of
 an irreducible character
$\theta \in \Irr((HB)(\gamma))$ lying above $\gamma$
and inducing $\lambda^{HB}$. The fact that $e=1$
implies that $\theta|_{B}= \gamma$, i.e.,
$\theta \in \Irr((HB)(\gamma))$ is an extension
of $\gamma \in \Irr(B)$ to $(HB)(\gamma)$.
Thus $\theta \in \Lin((HB)(\gamma)| \gamma)$
  induces $\lambda^{HB}$. Hence $\theta^G= \chi$, as
$\lambda$ induces $\chi$. Therefore, $ \theta^{G(\gamma)}$
is an irreducible character of $G(\gamma)$
lying above $\gamma$ and inducing $\chi$.
As the $\gamma$-Clifford correspondent $\chi_{\gamma}$  of $\chi$
is unique, we conclude that  $\theta^{G(\gamma)} = \chi_{\gamma}$.
Hence $\chi_{\gamma}$ is  induced from the linear character  $\theta$,
 and thus is monomial.

This completes the proof of the lemma in the case of an abelian $B$.
So the lemma follows.
\end{proof}

The above lemma implies
\begin{proposition}\mylabel{l4}
Assume that $T'= (G', N', \psi')$  is a linear limit of $T$. Assume further that $\chi'  \in \Irr(G'| \psi')$ 
 is the $T'$-reduction  of $\chi \in \Irr(G | \psi)$.
Then $\chi$ is monomial  if and only if $\chi'$ is monomial.
In particular, if $G$ is a monomial group, then every irreducible character of $G'$ that lies 
above $\psi'$ is monomial.
\end{proposition}

\begin{proof}
Since $T'$ is a linear limit of $T$, there exists some chain 
 $T_0 = T \geq  T_1 \geq  \cdots  \geq T_n = T'$ of linear subtriples of $T$, such that $T_i$ is a direct linear reduction
 of $T_{i-1}$, for all $i =1, \dots, n$.  So $T_i = T_{i-1}(\lambda_i)$, where   $\lambda_i$  is a linear character of some 
normal subgroup $L_i$ of $G_{i-1}$ with $Z(T_{i-1}) \leq L_i \leq N_{i-1}$.
  If $\chi$ is an irreducible character of $G$ lying above $\psi$, and thus above $\zeta^{T}$, then there  is also a finite chain 
of irreducible characters $\chi_i \in \Irr(G_i)$ starting with $\chi_0 = \chi$ and ending with $\chi_n = \chi'$, such that  
$\chi_i$ is the unique $\lambda_i$ -Clifford correspondent of $\chi_{i-1}$.
Hence Lemma \ref{l3} implies that $\chi_i$ is monomial if and only if $\chi_{i-1}$ is monomial.
We conclude that $\chi_0 = \chi$ is monomial if and only if $\chi'= \chi_n$ is monomial.

The rest of the proposition follows from Remark \ref{r10}.
\end{proof}

\section{Proof of Theorem C }

We begin with a straightforward lemma
\begin{lemma}\mylabel{C.l1}
Assume that $N$ is a finite solvable  group and let $M$ be a nilpotent normal subgroup of $N$
whose quotient group $N/M$ is also nilpotent. Assume further that every
$p$-Sylow  subgroup of $N$ centralizes the $q$-Sylow subgroup of $M$, for  all primes
$p \ne q $ where $p \bigm | |N|$ and $q \bigm |   |M|$.
Then $N$ is also nilpotent.
\end{lemma}

\begin{proof}
Let $P, Q$ be a $p$- and a $q$-Sylow subgroup of $N$, for two distinct
primes $p$ and $q$.  Let  $x \in P$ and $y \in Q$ be two elements
of $P$ and $Q$ respectively. It is enough to show that $[x, y]=1$.
Because $N/M$ is nilpotent while $(PM)/M$ and $(QM)/M$ are a $p$-and a $q$-Sylow subgroup, 
 respectively, of $N/M$, we have  $x\cdot  y = y\cdot  x \cdot m$, for some $m \in M$.
Hence $y^{-1} \cdot x \cdot y = x \cdot m_p \cdot m_{p'}$, where $m = m_p \times m_{p'}$
is the decomposition of $m$ to its $p$-part $m_p$ and its  $p'$-part $m_{p'}$.
So $x^{y} \cdot m_{p'}^{-1}= x \cdot m_{p}$. The right hand side of the last equation
is an element of $P$, and thus has order a power of $p$. On the other hand
$x^{y}$ is an element of some $p$-Sylow subgroup of $N$.
So $x^y$  commutes with  $m_{p'}$ by hypothesis.
 Therefore  the order of $x^{y} \cdot m_{p'}$ can be
 a power of $p$ only if $m_{p'}= 1$.
Similarly we have  $m_{q'} = 1$, where $m = m_q \times m_{q'}$ is the decomposition of
$m$ to its $q$-and $q'$-parts.
We conclude that $m= 1$, and the lemma follows.
\end{proof}

The following is the main tool for the proof of Theorem C 
\begin{lemma} \mylabel{C.l2}
Assume that $G$ is a monomial finite  group.  Assume further that $G$ has normal subgroups $M 
\leq N$ such that $M$ is nilpotent  with odd order and $N/ M$ is nilpotent.
 If $\phi$ is any irreducible character of $M$, and $T$ is the triple
$T= (G, M, \phi)$, then   there exists a linear limit $T' = (G', M', \phi')$  of $T$ so that
the factor  group $(G' \cap N)/ \Ker(T') $ is nilpotent.
\end{lemma}

\begin{proof}

Let $p_1, \dots, p_k$ be the distinct primes dividing $|N|$. 
Let $\{H_{i}\}_{i=1}^k$ be a Sylow system of
$N$.
 So $H_{i}$ is a $p_i'$-Hall subgroup of $N$. Furthermore,
$\cap_{j \ne i}H_{j} = Q_{i}$ is a $p_i$-Sylow subgroup of $N$.
We also write  $M_{i}$  for the  $p_i$-Sylow subgroup of $M$
(some could be trivial, and by hypothesis those that are not  have odd order).
So $M_{i}$ is a normal subgroup of $G$, for all $i=1, \dots, k$,  and
 $M = M_{1} \times M_{2} \times \cdots \times M_{k}$.
Note that  $M \unlhd H_{i}M_{i}$ while the quotient group $(H_{i} M_{i})/ M$
is the  $p_i'$-Hall subgroup of $N/M$. Hence $(H_{i} M_{i})/ M$
is a characteristic subgroup of $N/M$ and thus a normal subgroup of $G/M$. So
$H_{i}M_{i}$ is a normal subgroup of $G$.

Clearly the irreducible  character $\phi $ of $M$ can be written as
$\phi = \phi_{1} \times \cdots \times \phi_{k}$, where $\phi_{i} \in \Irr(M_{i})$
 for all $i=1,\dots, k$.
For every arbritrary but fixed $i=1, \dots, k$, we form the triple $T_i = (G, M_i,  \phi_{i})$.
Let $T_i' = (G_i' ,  M_{i}' , \phi_{i}')$ be a linear limit of  $T_i$.
We write $K_i' = \Ker(T_i')$ for  the kernel of $T_i'$, and   $Z_i' = Z(T_i')$ for the center of $T_i'$.
If $\zeta_i'  \in \Lin(Z_i')$ is the central character of $T_i'$, then
according to Proposition \ref{l7},  we have $G_i' = G_i(\zeta_i')$.

For all $j \ne i $ the $p_j$-Sylow subgroup $M_j$ of $M$  is a subgroup of $G_i'$,
since it centralizes $M_i$ (see Remark \ref{l01}). Hence the group
$L_{i}' = M \cap G_i'$ is a normal nilpotent subgroup of $G'_i$ whose $p_i$-Sylow
subgroup is $M_i'$ and whose  $p_j$-Sylow subgroup, for any  $j\ne i$, equals  $M_j$.
So
$$
L_i' = M_1 \times \dots \times M_{i-1} \times M_i' \times M_{i+1} \times \dots \times M_k.
$$

Let $N_i'$ be  $N \cap G_i'= N(\zeta_i')$ and  $H_{i}'$ be  a  $p_i'$-Hall subgroup of $N_i'$.
 Then $(H_i' \cdot M_{i}')/ L_i'$ is the  $p_i'$-Hall subgroup of the nilpotent factor group $N_i'/L_i'$,
and thus $H_i' \cdot M_i'= H_i'  \cdot L_i' $
is a normal subgroup of $G_i'$.

According to Proposition  \ref{l1},  the group $Z_i'/K_i'$ is the center
of the $p_i$-group $M_i'/K_i'$.
If $\zeta_i'/K_i'$ is the unique character  of the quotient group
 $Z_i'/K_i'$ that inflates to $\zeta_i'$, then  the same proposition implies that
for all $i=1, \dots, k$ the groups $M_i'/K_i',  G_i'/ K_i'$ and
the character $\zeta_i'/K_i'$ satisfy Condition X.
In addition, the group $H_i' K_i'/K_i'$  is a $p_i'$-Hall subgroup
 of $N_i'/K_i'$ and its  product with $M_i' /K_i'$  is a normal subgroup of  $G_i'/K_i'$.
According to Proposition  \ref{l4}    every irreducible character of $G_i'$ that lies above $\phi_i'$ 
is monomial.  But $\phi_i'$ is the only character of $M_i'$ lying above $\zeta_i'$, by Proposition \ref{l7}.
So every irreducible character of $G_i'$ lying above $\zeta_i'$ is monomial.

Hence  every irreducible character of $G_i'/K_i'$ that lies above $\zeta_i'/K_i'$ is monomial.

The character $\zeta_i'/K_i'$ is a $G_i'/K_i'$-invariant $p_i$-special
 character of $Z_i'/K_i'$, (one could see \cite{isa4}
 for the basic definitions of  $\pi$-special characters).
Hence   there exists an irreducible  $p_i$-special
 character of $G_i'/K_i'$ that lies above $\zeta_i'/K_i'$. Therefore, that character
 is monomial and the  $p_i'$-part of its degree is $1$.
 We can now apply Corollary \ref{coo1} to the groups $G_i'/K_i', M_i'/K_i'$
 and $(H_i' K_i')/ K_i'$, for all $i=1, \dots, k$. We conclude that
 $(H_i' K_i') /K_i'$ centralizes $M_i'/K_i'$ for all such $i$.
Hence the commutator subgroup  $[H_i' K_i', M_i']$ lies inside $K_i'$. Therefore,
\begin{equation}\mylabel{eq.c1}
[H_i', M_i'] \leq K_i',  \text{ for all } i=1, \dots, k,
\end{equation}
where $H_i'$  is any   $p_i'$-Hall  subgroup    of $N_i'$.

Let $G'= G_1'  \cap G_2' \cap \cdots \cap G_k'$. Then $G'=
G(\zeta_1', \dots, \zeta_{k}')$. We also define $N' =  N \cap G'=
 N(\zeta_1', \dots, \zeta_{k}')$ and $M' = M \cap G'=
  M(\zeta_1',  \dots, \zeta_k')$. Of course  $M' = M_1'  \times \cdots  \times M_k'$.
 Also $M '  \unlhd N'  \unlhd G'$. 
Furthermore, the group $Z' = Z_1' \times  \cdots  \times Z_k'$ 
 is a normal subgroup of $G'$ contained in $M'$.  The character $\phi'= \phi_1' \times \cdots  \times \phi_k'$ is an irreducible character 
of $M'$ that lies above the $G'$-invariant linear character 
$\zeta' = \zeta_1' \times \cdots \times \zeta_k'$ of $Z'$.
In view of Remark \ref{l2} the quintuple $T'= (G', M' , \phi')$ is a linear limit of $T = (G, M , \phi)$, 
while $Z' = Z(T') $,  $\zeta' $ is the central character $\zeta^{(T')}$ of $T'$, 
  and $K' = K_1' \times \cdots \times K_k'$ is the kernel $\Ker(T')= \Ker((\phi')^{G'})$ of $T'$.

In order to complete the proof of the lemma,  it suffices to show that $N' / K'$  is a nilpotent group.
According to Lemma \ref{C.l1}, it is enough to prove that every $p$-Sylow subgroup of $N'/K'$  centralizes 
every $q$-Sylow subgroup of $M'/K'$, whenever $q$ is a prime divisor of 
of $|M'/K'|$ and and  $p \ne q$ is a prime divisor of $|N'/K'|$. 
We fix a prime $p = p_r$   that divides $|N'/ K'|$, for some $r = 1, \dots, k $.
 Let $S'$ be a $p$-Sylow subgroup 
of $N'$. Clearly $N ' \leq N_i'$ for all $i=1, \dots, k$. Hence  for all $i \ne r$, there exists a $p_i'$-Hall subgroup $H_i ' $
of $N_i'$  so that $S' \leq H_i'$. Therefore \eqref{eq.c1} implies  
$$
[S', M_i'] \leq K_i', \text{ for all } i\ne r, i=1,\dots, k.
$$
So $[(S'K')/ K', (M_i'K')/ K') ] = 1$  for all $i \ne r$. Hence $N '/ K'$ is nilpotent. This proves Lemma \ref{C.l2}.

\end{proof}

Now we can prove Theorem C that we restate here.
\begin{C}
Assume that $G$ is a finite monomial group. 
Assume further that  $G$ has normal subgroups $M 
\leq N$ such that $M$ is nilpotent  with odd order and $N/ M$ is nilpotent.
 Let $\psi$ be  an irreducible character of $N$. Then there exists a linear limit 
$T' = (G', N', \psi')$ of the triple $T= (G, N, \psi)$  so that the factor group $N'/ \Ker(S')$  is nilpotent. 
Hence $N$ is a monomial group.
\end{C}

\begin{proof}
 We fix the character $\psi \in \Irr(N)$ and an irreducible character $\phi$ 
of $M$ lying under $\psi$. 
This way we can form two triples $S= (G, M , \phi)$ and $T= (G, N, \psi)$.
According to Lemma \ref{C.l2}, there exists a linear limit 
$S'' = (G'', M'', \phi'')$ of $S = (G, M, \phi)$ so that the factor group $N'' / \Ker(S'')$ is nilpotent, where 
 $N''= G'' \cap N$. If  $\psi'' \in \Irr(N'')$ is   the  $S''$-reduction of  $\psi$, then the 
 triple $T''= (G'', N'', \psi'')$ is  a linear reduction of $T= (G, N, \psi)$, see 
Remark \ref{l8}. Because $M'' \unlhd N''$ while  $\psi'' $ lies above $\phi''$ 
the same remark implies that 
$$
Z(S'') \leq Z(T'') \text{ and } \Ker(S'') \leq \Ker(T'').
$$
Of course, the direct linear reduction  $T''$ of $T$ does not need to be a linear limit of the latter,
 but certainly any  linear limit of $T''$ is also a linear limit of $T$. 
Let $T'= (G', N', \psi')$ be a linear limit of $T'' $ and $T$. 
Then Remark \ref{l00} implies that $\Ker(T'')  \leq \Ker(T')$, while $N' \leq N''$. Hence 
$\Ker(S'') \leq \Ker(T') \leq N' \leq N''$.
 This, along with the fact that  $N'' / \Ker(S'')$ is nilpotent,  implies that $N' / \Ker(T')$  is also a nilpotent group.
Therefore the first part of Theorem C is proved.

For the rest of the theorem, observe that $\psi'$ is a monomial character of $N'$, because $N'/ \Ker(T') $ is nilpotent
and $\Ker(T') \leq \Ker{\psi'}$. Because $\psi'$ induces $\psi$ in $N$, Theorem C follows. 
\end{proof}

\section{Proof of Theorem D}

The proof of Theorem D is heavily based on
\begin{theorem}\mylabel{cc:l0.5}
Let $Q$ be a $q$-group acting on a $p$-group $P$, with $p \ne q$ odd primes. We
identify both $P$ and $Q$ with their images in the semidirect product $QP = Q
\ltimes P$.
 Let $\T$ be a finite-dimensional right
 $\ZZ _q Q P$-module such that the action of $P$ on
$\T$ is faithful.
Then there exists an element $\tau \in \T$
  such that  its stabilizer $(QP)(\tau)$ in $Q\ltimes P$ equals $Q$.
\end{theorem}

\begin{proof}[Proof of Theorem \ref{cc:l0.5}]
We will prove a series of claims under the
\begin{ia}\mylabel{cc:1}
  $Q ,P, \T$  are  chosen among all the triplets satisfying the hypothesis,
but not the conclusion, of Theorem  \ref{cc:l0.5}, so as to minimize first the
 order $|QP|$  of the semidirect product $Q\ltimes P$, and then the
 $\ZZ_q$-dimension
$\dim_{\ZZ_q}\T$ of $\T$.
\end{ia}
\noindent
\setcounter{claim}{0}
These claims will lead to a contradiction, thus proving the theorem.
First note that  $\T \ne 0$. 
\begin{claim}\mylabel{cc:c1}
$\T$ is an indecomposable $\ZZ_qQP$-module.
\end{claim}

\begin{proof}
Suppose not. Let $\T=\T_1\, \dot{+} \, \T_2$ be a direct decomposition of $\T$,
where $\T_1 ,\T_2$ are nontrivial $\ZZ_q QP$-submodules of $\T.$
For $i=1,2$ let $K_i$ be the kernel of the action of $P$ on $\T_i.$
Hence $\T_i$ is a $\ZZ_q Q\ltimes (P/K_i)$-module such that $P/K_i$ acts
faithfully on it. As $\dim _{\ZZ_q} \T_i$ is strictly
smaller than  $\dim _{\ZZ_q}\T$, the minimality in Inductive Assumption
provides  an  element  $\tau_i \in \T_i$
such that $(Q \ltimes (P/ {K_i}))(\tau _i)=Q$.
(Here we have identifying $Q$ with its image in the semidirect product
$Q \ltimes (P/ K_i)$.)
If we take as $\tau$ the sum, $\tau =\tau_1 +\tau_2$,
then $\tau$ is an element of $\T$ fixed by $Q$, as $Q$ fixes
each one of the $\tau_i$ for $i=1,2.$
Furthermore for the stabilizer of $\tau$ in $P$ we have
$$
P(\tau)=\cap_{i=1}^2 P(\tau_i) = \cap_{i=1}^2 K_i.
$$
Since $P$ acts faithfully on $\T$ the last intersection is trivial. Therefore
$(QP)(\tau) =Q$, which contradicts the Inductive Assumption. Hence
$\T$ is an indecomposable $\ZZ_qQP$-module.
\end{proof}

\begin{claim}\mylabel{cc:c2}
The restriction $\T_P$ of $\T$ to  $P$ is a multiple of an irreducible
$Q$-invariant
$\ZZ_qP$-module.
\end{claim}

\begin{proof}
Because $q$ does not divide $|P|$  we can write
$\T_P$ as a direct sum of its $\ZZ_qP$-homogeneous
components, i.e.,
$$
\T_P  = \U_1 \dotplus  \U_2 \dotplus \dots \dotplus  \U_r.
$$
So $r \ge 1$ and there exist  distinct simple $\ZZ_qP$-modules
$\V_1, \V_2, \dots, \V_r$ and positive integers $m_1, m_2, \dots , m_r$  such
that $\U_i \cong  m_i \V_i$  as $\ZZ_qP$-modules, for all $i=1, \dots, r$.
Observe that right multiplication by any element in $QP$ permutes among
themselves the $\U_i$. So each $QP$-orbit  $\Omega$  of the $\U_i$
leads to a $\ZZ_qQP$-direct summand  $\sum _{\U_i  \in \Omega} \U_i$ of $\T_P$.
According to Claim \ref{cc:c1},  the group $QP$ acts transitively
on the $\U_i$. Hence $m_1= m_2 = \dots = m_r=m$.
Furthermore, if $\V_1= \V$ and $(QP) (\V)$ is the stabilizer of the isomorphism
class of  $\V$ in $QP$, then $V_i = \V^{\sigma_i}$, were ${1=\sigma_1, \dots,
\sigma_r}$ are  representatives for the cosets in $Q\cdot P$ of $(QP)(\V)$.
Thus,  $\U:=  \U_1 \cong m {\V} = m {\V}^{\sigma_1} ,\U_2 \cong m
{\V}^{\sigma_2}, \dots , \U_r \cong m \V^{\sigma_r}$.
We may pick ${\sigma_1 , \dots ,\sigma_r}$ to be representatives
of the cosets in $Q$ of  the stabilizer, $Q (\V)$, of the  isomorphism class of
$\V$  in  $Q$.
Note that $Q(\V) = Q(\U)$ as ${\U} \cong m\V$, where $Q(\U)$ is the stabilizer
in $Q$  of ${\U}$ under multiplication in $\T$.
If $\T_P$ is not homogeneous, then $r > 1$ and $Q(\U) = Q(\V) < Q.$
For $i=1, \dots , r$ let $K_i$  be the kernel of the action of $P$
on ${\U}_i$.
 Then  for every $i=1, \dots , r$ the stabilizer $Q({\U}_i)$ of ${\U}_i$ in $Q$
equals the $\sigma_i$-conjugate, $Q(\U)^{\sigma_i}$ ,  of $Q(\U) = Q(\V)$.
For the corresponding kernels we similarly have
 $K_i = K_1 ^{\sigma_i}$.

As  ${\U}$ is a faithful $\ZZ_q P/K_1$-module
and $Q(\U) <  Q$,  the minimality of $|QP|$ in the  Inductive Assumption
implies that there exists an element  $\mu \in  {\U}$ such that
$$
(Q(\U) \ltimes ( P/K_1))(\mu) = Q(\U).
$$
For every $i=1, \dots ,r$ we can define an element
$\mu_i = \mu  \sigma_i$ of
${\U}_i$.
 Then $Q(\U_i)= Q(\U_i)^{\sigma_i}$ fixes $\mu_i$ as $Q(\U)$ fixes $\mu$.
Furthermore if $x$ is any element of $P$ fixing $\mu_i$ then
$x^{\sigma_i^{-1}}$ is an element of $P$ fixing $\mu$. Therefore
$x^{\sigma_i^{-1}} \in K_1$, which implies that
 $x\in K_i$. Thus
$$
(Q(\U_i)\ltimes (P/K_i))(\mu_i) = Q(\U_i)
$$
for every $i=1, \dots , r$.

Let $\tau$ be the sum of the  $\mu_i$ for $i=1, \dots , r.$
Then $\tau $ is an element  of $\T$ fixed by $Q$, since multiplication by
any element in $Q$ permutes the ${\U}_i$  and the  $\mu_i$
among themselves.
The stabilizer $P(\tau)$ of $\tau$ in $P$ equals the intersection of
the stabilizers of $\mu_i$ in $P$ for $i=1, \dots ,r$.
Since $(Q(\U_i)\ltimes (P/K_i))(\mu_i) = Q(\U_i) $ for every
such $i$, the latter equals the intersection of $K_i$ for $i=1, \dots ,r.$ The
faithful action of $P$ on $\T$ implies that $$
P(\tau) = \cap_{i=1}^r K_i = {1}.
$$
Hence $\T$ has an element $\tau$ with $(QP)(\tau) = Q$,
contradicting the  Inductive Assumption.
This contradiction proves Claim  \ref{cc:c2}.
\end{proof}

\begin{claim}\mylabel{cc:c3}
There are  no $Q$-invariant subgroup, $H < P$, and $\ZZ_qQH$-submodule, ${\SB}$,
of
$\T_{QH}$ such that $\T$ is the $\ZZ_qQP$-module
 ${\SB}^{QP}$ induced from ${\SB}$, i.e.,
$$
\T = \sum_{1\leq i\leq n}^{\cdot}{\SB}\sigma_i,
$$
where the  $\sigma_i$ are representatives for the cosets $H \sigma_i$  of $H$ in $P$.
\end{claim}

\begin{proof}
Suppose Claim \ref{cc:c3} is false. We choose $H$ to have maximal order among
all  $Q$-invariant subgroups of $P$ that contradict Claim \ref{cc:c3}.
Hence $\T_{QH}$ has a 
$\ZZ_qQH$-submodule, ${\SB}$, such that
${\SB}^{QP} = \T$.
If $H$ is not normal in $P$ then its normalizer, $N_P(H)$, 
in $P$ satisfies $H \lhd N_P(H) < P$. Since $H$ is $Q$-invariant, $N_P(H)$ 
is also $Q$-invariant. Hence ${\SB}^{Q N_P(H)}$ is a $\ZZ_qQN_P(H)$-submodule of
 $\T_{QN_P(H)}$. Furthermore ${\SB}^{QN_P(H)}$ induces $\T$.
Thus $N_P(H)$ is among the $Q$-invariant subgroups of $P$ that contradict
Claim \ref{cc:c3},  while  $|N_P(H)| > |H|$. So the maximality of $|H|$ implies that
$H$ is normal in $P$. 

Let $1=\sigma_1 , \dots , \sigma_k$ be coset representatives of
$H$ in $P$, and let $\bar{\sigma}_m$ denote the
 image of $\sigma_m$ in $P/H$ for $m=1,\dots,k$.
 Then $\bar{1} = \bar{\sigma_1}, \bar{\sigma_2},\dots,\bar{\sigma_k}$
are the distinct elements of $P/H$.
As $Q$ acts on $P/H$, it has to divide  the
$\bar{\sigma}_m$, for $ m=1, \dots , k$, into  orbits, $\overline{R_1}, \overline{R_2}, 
\dots,\overline{R_l}$, for some
$l\in \{1, \dots ,k\}$. We may choose $\overline{R_1}$ to be
equal to  $\{{\bar{\sigma}_1}\} = \{{1}\} $. For every $i = 1, \dots ,l$, we
pick some element $\bar{\sigma}_{i,1} \in \overline{R_i}$. Then
 $\overline{R_i} =\{ \bar{\sigma}_{i,1}^{q_j}\}_{j=1}^{j=k_i} $  where
$k_i = |\overline{R_i}|$ and
$q_j $ runs over a set $Q_j$ of coset  representatives of the stabilizer,
$C_Q(\bar{\sigma}_{i,1})$, in $Q$.
For every $i=1, \dots , l$ the stabilizer $C_Q(\bar{\sigma}_{i,1})$ acts
by conjugation on $H$ and on $\sigma_{i,1}H$, where $\sigma_{i,1} \in P$ has
image $\bar{\sigma}_{i,1} \in P/H$.
Furthermore, $H$ acts transitively by right multiplication on 
$\sigma_{i,1}H$ and $(x h)^c = x^c h^c$  for all
 $x\in \sigma_{i,1}H , h\in H , c \in C_Q(\bar{\sigma}_{i,1}) $.
Hence Glauberman's Lemma (13.8 in \cite{is}) provides an element 
$t_{i,1} \in \sigma_{i,1}H$ that is fixed by $C_Q (\bar{\sigma}_{i,1})$.
So $C_Q(t_{i,1}) \geq  C_Q(\bar{\sigma}_{i,1})$.
Furthermore, the opposite inclusion, $C_Q(t_{i,1}) \leq 
C_Q(\bar{\sigma}_{i,1})$, 
 also holds as $\bar{\sigma}_{i,1} = t_{i,1}H$. Hence,
$$
C_Q(t_{i,1}) = C_Q(\bar{\sigma}_{i,1}).
$$
In this way we can pick a $t_{i,1} \in \sigma_{i,1}H$, for every ${i= 1, \dots , l}$, 
 such that $C_Q (t_{i,1}) = C_Q (\bar{\sigma}_{i,1})$.
We can even assume that  $t_{1,1} = 1$. 
Let $t_{i,j}$ denote the $q_j$-conjugate, $ t_{i,1}^{q_j}$, 
of $t_{i,1}$ for every $j=1, \dots , k_i$.
Hence the set of all $t_{i,j}$,  for  $ i=1, \dots l $ and for 
 $ j=1, \dots , k_i$,  is a complete 
set of coset representatives of $H$ in $P$.
Furthermore the  $Q$-orbit  $\overline{R_i}$  corresponds
  to a $Q$-orbit  $R_i = \{t_{i,1}, \dots , t_{i_{k,i}} \}$, 
for every $i=1, \dots , l$.

Let $K_{\SB}$ be the kernel of the action of $H$ on ${\SB}$.
As $|H/K_{\SB}| < |P|$, the minimality of $|QP|$  in  the Inductive Assumption 
implies that there exists  $\mu  \in {\SB}$ such that its stabilizer, 
$(Q\ltimes H/K_{\SB})(\mu)$, in $Q\ltimes H/K_{\SB}$ equals $Q$, or equivalently 
$(QH  )(\mu) =QK_{\SB}$. 
We note here that $K_{\SB} < H$. Indeed,  if $H$ acts trivially on ${\SB}$, then
$\T$ is  induced from a trivial module and thus  contains both trivial and
non--trivial irreducible
$\ZZ_qP$-submodules, contradicting Claim \ref{cc:c1}. We also
have that 
 $\mu \ne 0$ since $Q= (Q\ltimes H/K_{\SB})(\mu) < QH/K_{\SB}$.
We denote by $\mu t_{i,j}$ the $t_{i,j}$-translation  
of $\mu$, for every $i= 1, \ldots, l $ and for every 
$j=1, \dots , k_i$. Then $\mu t_{i,j}$ is an element of
 $\SB t_{i,j}$ such that
$$
(Q^{t_{i,j}}H)(\mu t_{i,j}) = Q^{t_{i,j}}K_{\SB}^{t_{i,j}}.
$$
Since  ${\SB}^{QP} = \T$  we get that
\begin{subequations}\mylabel{cc:eq1}
\begin{equation}\mylabel{cc:eq1a}
 \T = {\SB}^{QP} =  \sum_{1\leq i\leq l}^{\cdot} \sum_{1\leq j \leq k_i}^{\cdot}
\SB t_{i,j}={\SB} \dotplus \sum_{2\leq i\leq l}^{\cdot} \sum_{1\leq j \leq k_i}^{\cdot}
\SB t_{i,j}.
\end{equation}

Let $\tau$ be the element  of $\T$ defined by
\begin{equation}\mylabel{cc:eq2}
\tau = -\mu + \sum_{i=2}^{l} \sum_{j=1}^{k_i} \mu{t_{i,j}}=
 -\mu + \sum_{i=2}^l  \sum_{t_{i,j} \in R_i} \mu{t_{i,j} }.
\end{equation}
\end{subequations}
We claim that $\tau $ satisfies the condition in Theorem  \ref{cc:l0.5}, i.e., that
$(QP)(\tau) = Q$. This will  contradict the Inductive Assumption,
and thus prove Claim \ref{cc:c3}.
Indeed, $R_i = \{t_{i,1}, \dots ,t_{i,k_i} \}$ is a $Q$-orbit for every
$i=2, \dots , l$.
Also $\mu$ and  $-\mu$ are $Q$-invariant as 
$(QH)(-\mu) = (QH)(\mu) = QK_{\SB}$.
Hence $\sum_{t_{i,j} \in R_i} \mu t_{i,j}$  is   $Q$-invariant.
Thus $\tau$ is  a $Q$-invariant element of $\T$, which implies that
$(QP)(\tau)= Q \cdot P(\tau)$.

If $x\in  H(\tau)$ then, since $H\lhd P$, we get
that $(\mu {t_{i,j}})x =\mu x^{(t_{i,j})^{-1}}  t_{i,j}$ 
is an element of ${\SB}t_{i,j}$,  for all $i=2, \dots, l$ and $j= 1,\dots,k_i$,
while $(-\mu)x$ is an 
element  of ${\SB}$.
Since $\tau x = \tau$, it follows from \eqref{cc:eq1}  that   $(-\mu)x =
-\mu$ and 
$(\mu {t_{i,j}})x =\mu {t_{i,j}}$ for every $i= 2, \dots , l$
and for every $j=1, \dots , k_i$.
Hence $x$ is an element of:
$$
H(\mu) \cap   \bigcap_{i=2}^l \bigcap_{j=1}^{k_i}(P(\mu {t_{i,j}}) \cap H) = 
\bigcap_{i=1}^l\bigcap_{j=1}^{k_i} H(\mu {t_{i,j}}) = \bigcap_{i=1}^l
\bigcap_{j=1}^{k_i} K_{\SB}^{t_{i, j}}.
$$
As $H$ acts faithfully on $\T$, we get that $\bigcap_{i=1}^l
\bigcap_{j=1}^{k_i} K_{\SB}^{t_{i, j}} = {1}$.
Hence $H(\tau)  = {1}$.

Now let  $x \in P \smallsetminus H$. We claim that $\tau x \ne \tau$.
Indeed any  $x \in P $  permutes the ${\SB}t_{i, j}$ among themselves. 
If $x$ fixes $\tau$, then  it also  permutes 
among themselves the  summands  $-\mu$ and $\mu t_{i,j}$, for $i \ne 1$, of
$\tau$. Since ${\SB}x \ne {\SB}$ we have
 $(-\mu) x = \mu t_{i,j}$ for some $i=2,\dots,l$ and 
some $j = 1,\dots,k_i$. But as
 $x \in P \smallsetminus H$ we have that  $x = h t$
 for some coset representative 
$t = t_{i_0,j_0}$ of $H$ in $P$ with $i_0 = 2, \dots ,l$
and some element  $h \in H$.
 Hence $\mu t_{i,j}= (-\mu)x = (-\mu)h t \in {\SB} t$, which implies that 
$t_{i,j} = t$ and  $(-\mu) h = \mu$.
This last equation  leads to a contradiction  as $h$ has odd order 
($|P|$ is odd)
and $\mu \ne -\mu$ ( as ${\SB} \leq \T$ has odd order, while $\mu \ne 0$).
Therefore  $\tau^x \ne \tau $ whenever $x \in P\smallsetminus H$.
Hence $P(\tau) = H(\tau) =1$ and $(QP)(\tau) =Q$, 
contradicting the Inductive Assumption.
This contradiction proves Claim \ref{cc:c3}.
\end{proof}

\begin{claim}\mylabel{cc:c4}
The restriction $\T_A$ of $\T$ to any normal subgroup $A \leq P$ of $QP$ 
is a multiple $e \V$ of a single faithful $QP$-invariant 
$\ZZ_p A$-module $\V$.
  Hence 
every normal abelian subgroup $A$ of $QP$ contained in $P$ is cyclic.
\end{claim}

\begin{proof}
Let $A$ be a normal subgroup of $QP$ contained in $P$, and let $\T_A$ be the
restriction of
 $\T$ to $A$. According to  Claim \ref{cc:c2},
 and Clifford's Theorem, 
$\T_A$ can be written as a direct sum of its $\ZZ_qA$-homogeneous
components, i.e.,
$$
\T_A = {\W}_1 \dotplus {\W}_2 \dotplus \dots \dotplus {\W}_s.
$$
Furthermore, $P$ acts transitively on  the ${\W}_i$ for all $i=1,\dots,s$,
while $Q$ permutes the ${\W}_i$ among themselves (as $\T$ is a $\ZZ_qQP$-module).
Hence Glauberman's lemma implies that $Q$ fixes some
$\ZZ_qA$-homogeneous component, ${\W}$,  of $\T_A$. Note that Clifford's theorem
implies that the homogeneous component ${\W}$ of $\T_A$
is  a $\ZZ_qQP(\W)$-submodule of $\T$, where $P(\W)$ is the stabilizer
of $\W$ under multiplication of elements of $\T$ by
elements of $P$. Furthermore, the $\ZZ_qQP(\W)$-submodule
$\W$ induces the $\ZZ_qQP$-module $\T = \W^{QP}$.
  In view of  Claim \ref{cc:c3}, we must have
${\W} =  \T$. Hence $\T_A \cong e {\V}$  where  ${\V}$ is an irreducible
$QP$-invariant $\ZZ_qA$-submodule of $\T$. As $P$ acts faithfully  on $\T$, the
$\ZZ_qA$-module ${\V}$ is also faithful. If $A$ is abelian, the existence
of a faithful irreducible
 $\ZZ_qA$-module implies that $A$ is cyclic.
Therefore,  the claim is proved.
\end{proof}

The $q$-group $Q$ acts on the non--trivial  $q$-group
$\T$, fixing the trivial element  $0$ of $\T$. Hence the group
$Q$ fixes at least  $q$ elements of $\T$.
So $Q$ fixes some $\tau$ with
\begin{equation}\mylabel{cc:eq3}
 \tau \in \T \text{ and }  \tau \ne 0.
\end{equation}
Hence, to complete the proof of Theorem  \ref{cc:l0.5}, by contradicting
the Inductive Assumption,
 it  is enough to show that $P(\tau) = 1$

By Claim \ref{cc:c4} every characteristic abelian subgroup of $P$ is
cyclic. Since $p$ is odd, Theorem 4.9  of \cite{go} implies
that either $P$ is cyclic or $P$
 is the central product $E \odot C$,  of  the extra--special $p$-group $E= \Omega_1(P)$
of exponent $p$,   and the  cyclic group $C= Z(P)$.

According to Claim \ref{cc:c4}, the $\ZZ_qZ(P)$-module $\T_{Z(P)}$ is a multiple of
a faithful irreducible $QP$-invariant $\ZZ_qZ(P)$-module ${\V}$, i.e.,
$\T_{Z(P)} = m {\V}$. Hence $Z(P)$ acts fix point freely  on $\T$, as it acts fix point freely on ${\V}$ (or else ${\V}$ 
wouldn't be simple and faithful).

If $P$ is cyclic, then $P = Z(P)$.  Thus $P$ acts fix  point freely on $\T$.
Hence  no element of $  Z(P) - \{ 1\}$ could fix $\tau$. Hence 
$P(\tau) = 1$. So $(QP)(\tau) = Q$, 
contradicting the Inductive Assumption.
Therefore,  $P$ can't be cyclic. 

Hence, 
\begin{equation}\mylabel{cc:e2}
P = E \odot C = \Omega _1(P) \odot Z(P), 
\end{equation}
where $E = \Omega _1(P)$ is an extra special $p$-group of exponent $p$  
 and $C = Z(P)$ is cyclic. 
Therefore the quotient group $\overline{P} = P/ Z(P)$ is an elementary 
abelian $p$-group. Furthermore $\overline{P}$ affords a bilinear form
 $c: \overline{P} \times \overline{P} \to Z(E)$ defined,  for every 
$\bar{x}, \bar{y} \in \overline{P}$, as $c(\bar{x} , \bar{y}) = [x,y]$,
 where
$x, y$ are  elements of $P$ whose  images  in $\overline{P}$ are 
$\bar{x}$ and $  \bar{y}$ respectively. 
With respect to that form $\overline{P}$ is a symplectic $\ZZ _p (Q)$-module.

\begin{claim}\mylabel{cc:c5}
The symplectic $\ZZ _p(Q)$-module $\overline{P}$ is anisotropic.
 \end{claim}

\begin{proof} 
Assume not. Then there is an isotropic non--zero $\ZZ _p(Q)$-submodule
$\bar{A}$ of $\overline{P}$. Hence $c(\bar{a} , \bar{b}) = 0$ for every
$\bar{a}, \bar{b} \in
 \bar{A}$,  because  $\bar{A} \subseteq \bar{A}^{\perp}$.
Therefore,  the definition of the symplectic
form $c$  implies  that the inverse image $A$ of $\bar{A}$ in $P$ is
an abelian subgroup of $P$ containing $Z(P)$.
Since $\bar{A}$ is a $\ZZ _p(Q)$-submodule of $\overline{P}$, the abelian
group $A$ is a normal
 subgroup of $QP$ contained in $P$.
Hence by Claim \ref{cc:c4} , $A$ is cyclic and properly
 contains  $Z(P)$.
Therefore there exists an element $a  \in A\smallsetminus Z(P)$
such that $a^p$ is a generator of $Z(P)$. On the other hand,  
equation  \eqref{cc:e2} implies that  $a = \omega \cdot c $
where $\omega \in \Omega_1(P)$ and $c \in C=Z(P)$.
Hence $a^p = \omega^p \cdot c^p = c^p$. Since $a^p$ is a generator 
of the cyclic non--trivial 
 $p$-group $Z(P)$ and $c \in Z(P)$, this last equation leads  to a
contradiction.
This proves the claim.
\end{proof}

Now we can complete the proof of  Theorem  \ref{cc:l0.5}.
If $(QP)(\tau ) \ne Q$ then there exists a $Q$-invariant subgroup 
$D = P(\tau)\ne 1$
 of $P$ such that $(QP)(\tau) = QD$. 
Hence the center $Z(D)$ of $D$ is a non--trivial $Q$-invariant abelian subgroup of $P$.
Therefore its
 image $\overline{Z(D)} = Z(D)Z(P)/Z(P)$ in
 $\overline{P}$ is an isotropic $\ZZ _p(Q)$-submodule of $\overline{P}$. 
Since $\overline{P}$ is anisotropic, $\overline{Z(D)} = \bar{1}$, i.e.,
$Z(D)$ is contained in $Z(P)$.

As we saw,  $Z(P)$ acts fix point freely  on $\T$.
This implies that no element of $Z(P) - \{1\}$ could fix $\tau$.
Hence $Z(D) = 1$, contradicting the fact that $Z(D) \ne 1$.
 So $(QP)(\tau) = Q$, contradicting the Inductive Assumption.
This final  contradiction completes the proof of Theorem  \ref{cc:l0.5}. 
\end{proof}

In terms of characters,  Theorem  \ref{cc:l0.5} implies
\begin{corollary}\mylabel{cc:c0.5}
Let $Q$ be a $q$-group acting on a $p$-group $P$ with $p\ne q$ odd primes.
Suppose that the semi--direct product $Q\ltimes P$ acts on a $q$-group $S$
such that the action of $P$ on $S$ is faithful. Then there exists a linear
character $\lambda$ of $S$
whose kernel $\Ker(\lambda)$ contains the Frattini subgroup
$\Phi(S)$ and whose  stabilizer $(QP)(\lambda)$ in $Q\ltimes P$
is $Q$.
\end{corollary}

\begin{proof}
Let $\T$ be the quotient group $\T:= S/ \Phi(S)$.
 Then $\T$ is a $\ZZ_qQP$-module.
 We write  $\T^*$ for  its  dual $\ZZ_qQP$-module, i.e.,
$\T^*= \Hom_{\ZZ_q}(\T,\ZZ_q)$. Then $P$ acts faithfully on both
$\T$ and $\T^*$. Furthermore, according to  Theorem   \ref{cc:l0.5}
there is an element $\tau  \in \T^*$ whose stabilizer in $QP$ equals $Q$.
Since the linear characters of $\T$ can be considered as the elements of $\T^*$
 composed with some faithful  linear character of $\ZZ_q$, we conclude that
there is a linear character $\lambda^* \in \Lin(\T)$ whose stabilizer in
$QP$ is $Q$.
Let $\lambda$ be the linear character of $S$ to which
$\lambda^*$ inflates. Then $\Phi(S) \leq \Ker(\lambda)$.
Furthermore, $(QP)(\lambda) = (QP)(\lambda^*) = Q$,
and the corollary follows.
\end{proof}

The following is a straightforward lemma.
\begin{lemma}\mylabel{lemmaA}
Let $P$ be a $p$-subgroup of a finite group $G$ and let  $Q_1 \leq  Q$
be $q$-subgroups
of $G$, for some distinct odd primes $p$ and $q$.
If $P$ normalizes $Q_1$, and $Q$ normalizes their product $Q_1P$, then
$ QP$ is also a subgroup of $G$ with  $Q \in \Syl_q(Q P) ,
P \in \Syl_p(QP )$, while
$Q_1P  \unlhd QP$ and  $Q_1 \unlhd QP$.
 Furthermore,  $Q$ is the product
$Q = [Q_1 , P] N_{Q}(P)$, where $[Q_1,P] \unlhd QP$ and
 $[Q_1, P] \cap N_{Q}(P) = C_{[Q_1, P]}( P ) \leq \Phi([Q_1, P]).$
\end{lemma}
\begin{proof}
Since $P$ normalizes $Q_1$, the latter is a characteristic subgroup of $Q_1P$.
Therefore, the fact that  $Q$ normalizes $Q_1P$ implies that $Q$ normalizes
$Q_1$. So $Q_1 \unlhd Q$.

The product,  $Q P = Q (Q_1 P)$, is a subgroup of $G$, since
 $Q$ normalizes the semidirect product $Q_1 \rtimes P$.
That same  product $Q_1 P$ is a normal subgroup of $QP =Q(Q_1 P)$.
We obviously have that  $Q \in \Syl_q(Q P)$ and  $ P\in \Syl_p(Q P)$.

By Frattini's argument for the Sylow $p$-subgroup  $P$ of  $Q_1P\unlhd QP$
we get
\begin{subequations}
\begin{equation}\mylabel{A1}
Q P = Q_1 P N_{QP} (P ).
\end{equation}
The normalizer, $N_{QP}(P )$, of $P$  in $Q P$  contains $P$. So
it is equal to $P N_{Q}(P)$.
 Hence  \eqref{A1} can be written as  $Q P = Q_1 N_{Q}(P) P$.
Since $Q_1 N_{Q}(P) \leq Q$  and $Q \cap P = 1$,  we
conclude  that
\begin{equation}\mylabel{A2}
Q = Q_1 N_{Q}(P).
\end{equation}
\end{subequations}
Because $(|Q_1|, |P|) = 1$, and $P$ acts on $Q_1$, we can write $Q_1$ as the
product  $Q_1 =[Q_1 , P] N_{Q_1}( P)$.
The commutator subgroup $[Q_1, P]$ is a characteristic subgroup of $Q_1 P$
and thus is also a
normal subgroup of $Q$,  as  $Q$ normalizes $Q_1 P$.
Therefore, \eqref{A2} implies
$$Q = [ Q_1, P] N_{Q}(P).$$
That $ [Q_1, P] \cap N_{Q}(P) = C_{[Q_1, P]}( P ) $ is obvious as
$(|Q_1|, |P|) = 1$.
 Also
the quotient group $K:= [Q_1,P]/\Phi([Q_1,P])$ is abelian and thus
$K=[K,P] \times C_K(P )$. As $[Q_1,P,P] = [Q_1,P]$
(by Theorem 3.6   in \cite{go}),
 we get that $K = [K,P]$ and
$C_K(P )= 1$. This  implies that
 $C_{[Q_1, P]}(P) \leq \Phi([Q_1,P])$.
\end{proof}

As an  easy consequence of Corollary  \ref{cc:c0.5} and
 Lemma \ref{lemmaA} we have:
\begin{proposition}\mylabel{cc:l1}
Let $Q$ be a $q$-group acting on a $p$-group $P$ with $p\ne q$ odd primes.
Suppose that the semi--direct product $Q \ltimes P$ acts on a $q$-group $S$
such that the action of $P$ on $S$ is faithful. Then there exists a linear
character $\lambda$ of $S$
such that  $C_{S}(P)  \leq \Ker(\lambda)$ and
 $(QP)(\lambda) = Q$.
\end{proposition}

\begin{proof}
As $P$ acts on $S$ we can write $S$ as the product $S=[S,P] \cdot
C_{S}(P)$. It is clear that the product $ Q C_{S}(P)$ forms a group.
Furthermore, $Q C_{S}(P)$ normalizes $P$ and
the semidirect product $(Q C_{S}(P))\ltimes P$
acts on $[S,P]$, while the action of $P$ on $[S,P]$ is faithful.
Then according to Corollary \ref{cc:c0.5} there exists a linear character
$\lambda_1$ of $[S,P]$ such that
$(Q C_{S}(P)P)(\lambda_1) =  Q C_{S}(P)$, while
 $\Phi([S,P]) \leq \Ker(\lambda_1)$.

  As we have seen in Lemma \ref{lemmaA}
$$
[S,P] \cap C_{S}(P) = C_{[S, P]}( P )  \leq \Phi([S,P]).
$$
Since $\lambda_1$ is a  linear character of $[S,P]$
that is trivial on $\Phi([S,P])$ and $C_{S}(P)$-invariant,
 the above inclusion implies that
 $\lambda_1$ has a unique extension to a linear character
$\lambda$ of $S$ trivial on $C_{S}(P)$.
Furthermore,
$(QP)(\lambda) = (QP)(\lambda_1) = Q$,
 and  the proposition follows.
\end{proof}

We can now prove a special case of Theorem D  where $P(\beta)$ is trivial. In this special case  the new character
 we get is linear.  In particular we have
\begin{lemma}\mylabel{le1}
Let $P$ be a $p$-subgroup  of a finite group $G$, where $p$ is an odd prime.
Let  $Q_1, Q$ be $q$-subgroups of $G$ for some odd prime $q \ne p$, with $Q_1 \leq Q$.
Assume that $P$  normalizes $Q_1$ while $Q$ normalizes the product  $ Q_1 P$.
Assume further that $\beta$ is an irreducible character of $Q_1$ such that
 $P(\beta)= 1$. Then there exists  a linear character
$\lambda$ of $Q_1$  such that $1= P(\beta ) = P(\lambda)$,  while  ${Q}(\beta)
\leq {Q}(\lambda) =Q$ and  $\lambda$ extends to
$Q$.
\end{lemma}
\begin{proof}
Because $Q$ normalizes the product $Q_1 \rtimes  P$, it normalizes its
characteristic subgroup $Q_1$. Let
  $C = N_{Q}(P)$ be the normalizer of $P$ in $Q$. Frattini's argument
implies that
$$
 Q = N_{Q}(P)  Q_1 = C   Q_1.
$$
In addition, $C $ normalizes  $P$ and  the semidirect product $C P $ acts on $Q_1$.
The fact that $P(\beta)= 1$ implies that $P$ acts faithfully on $Q_1$.
Therefore we can apply Proposition \ref{cc:l1} to the groups $C , P$ and $Q_1$
here in the place of $Q, P$ and $S$ there respectively. We conclude that there
exists a linear  character $\lambda \in \Lin(Q_1)$ such that $C_{Q_1}(P) \leq
\Ker(\lambda)$ and $(CP)(\lambda) = C$.

This last equation implies that $P(\lambda) = P(\beta) = 1$. Since $Q= C Q_1$ and $C$  fixes $\lambda$ we conclude
that $Q$ also fixes $\lambda$.  Furthermore we have
$$
C \cap Q_1 = N_{Q}(P) \cap Q_1 = C_{Q_1}(P) \leq \Ker(\lambda).
$$
Therefore $\lambda$ can be extended to $Q $.

As $Q(\beta) \leq Q =Q(\lambda)$, the proof  of Lemma \ref{le1} is complete.
\end{proof}

Finally we  can now prove  Theorem D, that we restate here.
\begin{D}
Let $P$ be a $p$-subgroup, for some  odd prime $p$, of a finite group $G$.
Let $Q_1,  Q$ be   $q$-subgroups of  $G$, for some odd prime
$q\ne p$, with $Q_1 \leq  Q$.
Assume that $P$  normalizes $Q_1$,
while $Q$ normalizes the product  $Q_1 \rtimes  P$.
Assume further that  $\beta$ is an irreducible character of $Q_1$.
 Then there exists an  irreducible character $\bn$ of $Q_1$
such that
\begin{align*}
P(\beta ) &= P(\bn), \\
Q(\beta) \leq Q(\bn) \, &\text{ and } N_Q(P(\beta)) \leq Q(\bn),  \\
\bn &\text{  extends to }  Q(\bn).
\end{align*}
\end{D}

\begin{proof}[Proof of Theorem D]
Let $P(\beta)$ be the stabilizer of $\beta$ in $P$ and $P_1$ be the
normalizer of
$P(\beta)$ in $P$.  Let $\overline{P_1}$ denote the quotient group $P_1 /
P(\beta)$.
We write  $C_1$ for  the centralizer, $C_1 = C(P(\beta) \tin Q_1)$,  of $P(\beta)$
in $Q_1$. Then it is clear that $\overline{P_1}$ acts
 on $C_1$.

The Glauberman correspondence (Theorem 13.1 in \cite{is}), applied to
the groups $P(\beta)$ and $Q_1$, provides an irreducible character
$\beta^*$ of $C_1$ corresponding to the irreducible character $\beta$
of $Q_1$.
Because  $P_1$ normalizes both  $P(\beta)$ and $Q_1$
 we get that  $P_1(\beta^*) = P_1(\beta) = P(\beta)$.
If $P_1(\beta^*)  < P(\beta^*)$ then  
$P_1(\beta^*) < N (P_1(\beta^* ) \tin P(\beta^*)) = N(P(\beta) \tin P(\beta^*)) =
 P_1(\beta^*) $, which is impossible. Therefore 
$$
P(\beta^*) = P_1(\beta^*) = P(\beta)= P_1(\beta).
$$ 
We remark here that because 
$P(\beta)$ centralizes $C_1 =C(P(\beta) \tin Q_1)$, we have 
$P(\beta) \leq C(C_1\tin P_1)\leq P_1(\beta^*) = P(\beta)$.
Hence $C(C_1\tin P_1) = P(\beta)$ and $\overline{P_1}$ acts faithfully 
on $C_1$. 

Let $C :=N(P(\beta) \tin Q)$ be the normalizer of $P(\beta)$ in $Q$.
Then $C_1$ is a normal subgroup of $C$ as $Q_1 \unlhd Q$.
 Furthermore, $C$  normalizes    $N(P(\beta) \tin PQ_1)$ 
 because  $Q$ normalizes  the product $PQ_1$. 
 As $P_1 C_1 = N(P(\beta) \tin PQ_1)$ 
 we conclude that $C$ normalizes the product 
$P_1 C_1$. Hence Frattini's argument implies that 
$$
C = N(P_1 \tin C) C_1.
$$

Now we can apply Lemma \ref{le1} to the groups $C, C_1$ and $\bar{P_1}$ and the  character $\beta^* \in \Irr(C_1)$,
in the place of $Q, Q_1, P$ and $\beta$ respectively. We conclude that 
 there exists a linear character
$\lambda \in \Lin(C_1)$  such that
\begin{subequations}\mylabel{le3}
\begin{align}
\bar{P_1}(\beta^*)= 1 &= \bar{P_1}(\lambda), \mylabel{le3a}\\
{C}(\beta^*) \leq {C}(\lambda)&= N_C(C_1) = C   \mylabel{le3b}\\
\lambda &\text{ extends to } C(\lambda)= C. \mylabel{le3c}
\end{align}
\end{subequations}

Equation \eqref{le3a} above  implies that $P_1(\lambda) = P(\beta)$.
 Thus $P(\beta) = P_1(\lambda) \leq P(\lambda)$.
We actually have that $P(\lambda) =P(\beta)$.
Indeed,  if $P(\beta) < P(\lambda)$, 
 then $P(\beta)$ would be a proper subgroup of  
$N(P(\beta) \tin P(\lambda))$.
Thus $P(\beta) < N(P(\beta) \tin P(\lambda)) = N(P(\beta) \tin P)
(\lambda) = P_1(\lambda) = P(\beta)$. So 
$$
P_1(\lambda) = P(\beta) = P(\lambda).
$$

 Let $\bn \in \Irr(Q_1)$ be the Glauberman  $P(\beta)$-correspondent  to $\lambda$.
Because  $C P_1$ normalizes both $P(\beta)$ and $Q_1$ we get 
$(CP_1)(\bn) = (CP_1)(\lambda) = CP(\beta)$.
Hence $P(\bn) \geq P_1(\bn) = P_1(\lambda) = P(\beta)$.
If $P(\bn) >  P(\beta)$ then  $
P(\beta)  < N(P(\beta) \tin P(\bn)) =
P_1(\bn)= P(\beta)$.
Thus    
$$
P(\bn) = P(\lambda) = P(\beta)
$$ and 
$( C(PQ_1))(\bn) = C P(\beta)Q_1$.
Since $C$ fixes $\bn$ and normalizes $P(\beta)$ we have
$C \leq N(P(\beta) \tin Q(\bn)) \leq N(P(\beta) \tin Q) = C$.
Hence
$$
 N(P(\beta) \tin Q) = C = N(P(\beta) \tin Q(\bn))\leq Q(\bn).
$$

In order to show that $\bn$ extends to $Q(\bn)$, we first observe that 
$Q_1 P(\beta) = (Q_1 P ) (\bn)= (Q_1 P_1)(\bn)$ as $P(\bn) = P(\beta) = P_1(\bn)$.
So the group $  Q_1P(\beta)  = (Q_1P)(\bn)$ is a normal subgroup
 of $Q(\bn) P(\beta) $ as
$Q$ normalizes the product $Q_1 P$.
Because $\lambda$ extends to $C$ while  $P(\beta)$ fixes $\lambda$ and has coprime order to that of
$C$, we conclude that $\lambda$ can be
extended to  $C P(\beta)$, (see Theorem 6.26 in \cite{is}).
 So we can apply the 
 Main Theorem in \cite{lew} to the groups $P(\beta) Q(\bn)$,
 $ P(\beta) Q_1$ and $Q_1$.
We conclude that $\bn$ extends  to $P(\beta) Q(\bn)$ as 
its $P(\beta)$-Glauberman correspondent
$\lambda$ can be extended to $ P(\beta) C =P(\beta) N(P(\beta) \tin Q(\bn))$.
We write $\beta^{\nu, e}$ for an  extension of $\bn$ to $Q(\bn)$.

To complete the proof of the theorem it remains to show that 
$Q(\beta) \leq Q(\bn)$.
The group $(Q_1 P) (\beta) = Q_1 P(\beta) $ is a normal subgroup of 
$Q(\beta) P(\beta) $, as $Q$ normalizes $Q_1P $.
Hence Frattini's argument implies that 
$$Q(\beta) = Q_1 N(P(\beta) \tin Q(\beta)).
$$
Therefore  $Q(\beta) \leq Q_1 N(P(\beta) \tin Q)= Q_1 C$.
But  we have already seen that $ Q_1 P(\beta)$ is a normal subgroup of
 $Q(\bn) P(\beta) $.
 Hence the  Frattini argument implies that 
$$
Q(\bn) = Q_1 N(P(\beta) \tin Q(\bn))
= Q_1 C.
$$ 
Thus $Q(\beta) \leq Q(\bn)$ and the theorem follows.
\end{proof}

\section{ Proof of Theorem E}
We first need some lemmas.
\begin{lemma}\mylabel{l.e1}
Assume that $G$ is a finite group and that $S \leq H $ are subgroups of
 $G$ with  $S$ normal in $G$. Assume further that
$\theta \in \Irr(H)$ lies above $\lambda \in \Irr(S)$. Let
 $\theta_{\lambda} \in \Irr(H(\lambda))$
denote the unique $\lambda$-Clifford correspondent
of $\theta$.
If $\theta$ extends to its stabilizer $G(\theta)$ in $G$, then
$\theta_{\lambda}$ also  extends to
$ G(\theta, \lambda)$.
\end{lemma}

\begin{proof}
A straight forward application of Clifford Theory implies that
$G(\theta, \lambda) \leq  G(\theta_{\lambda})$.
 Furthermore, as $G(\theta)$ fixes $\theta$ it permutes among themselves the
members of the  $H$-conjugacy
class of characters in $\Irr(S)$ lying under $\theta$. Since $\lambda \in \Irr(S)$
lies under $\theta$ we get
\begin{subequations}\mylabel{prel.t5.1}
\begin{equation}
G(\theta) = H \cdot G(\theta, \lambda) \leq  H \cdot G(\theta_{\lambda}).
\end{equation}
In addition,
\begin{equation}
G(\theta, \lambda) \cap H = H( \lambda).
\end{equation}
\end{subequations}

Let $\theta^e \in \Irr(G(\theta))$ be an extension of $\theta$ to $G(\theta)$.
Then  $\theta^e$ lies above $\lambda $. Let $\Psi \in \Irr(G(\theta, \lambda))$
denote the unique $\lambda$-Clifford correspondent of $\theta^e \in \Irr(G(\theta))$.
So $\Psi $ lies above $\lambda$ and induces $\theta^e$.
 Therefore,
$$(\Psi^{G(\theta)})|_H = \theta^e |_H = \theta.$$
Mackey's Theorem, along with \eqref{prel.t5.1},  implies that
$$(\Psi^{G(\theta)})|_H = (\Psi|_{H(\lambda)})^H.$$
Hence $ (\Psi|_{H(\lambda)})^H= \theta$ is an irreducible character of
$H$. So the restriction $\Psi|_{H(\lambda)}$
is an irreducible character of $H(\lambda)$ that induces $\theta$
and lies above $\lambda$ (as $\Psi $ lies above $\lambda$).
We conclude that $\Psi|_{H(\lambda)} $ is the $\lambda$-Clifford correspondent of
$\theta$.  Hence $\Psi|_{H(\lambda)} = \theta_{\lambda}$.
Thus $\Psi$ is an extension of $\theta_{\lambda}$  to $G(\theta, \lambda)$,
and the lemma  follows.
\end{proof}

Applying this lemma to direct linear reductions we can prove 
\begin{lemma} \mylabel{l.e1n}
Let $T = (G, N , \psi) $ be a triple, and $T' = (G', N', \psi')$ be a linear reduction of $T$. Let $H$ be any
 subgroup of $G$ with $N \leq H$. If $\zeta'$ is the central character of $T'$, then $H'= H \cap G'  = H (\zeta')$. Also 
 if $\psi$ extends to $H(\psi)$ then $\psi'$ extends to $H'(\psi)= H'(\psi')$.
\end{lemma}

\begin{proof}
Assume first that $T'$ is a direct linear reduction of $T$. So $ T' = T(\lambda)$ where  $\lambda$ is a linear character of  a 
normal subgroup $L$ of $G$ with $Z(T) \leq L \leq N$. Furthermore, $\lambda$ lies above the central character 
$\zeta \in \Irr(Z(T))$ of $T$.  
In addition, if $Z(T')$ is the center of $T'$, and $\zeta' \in \Lin(Z(T')$ its central character, then 
$Z(T) \leq L \leq Z(T') \leq N$ while $\zeta'$ lies above $\lambda$. 

According to Proposition \ref{l7}, we have  $G(\lambda) = G' = G(\zeta')$. 
Hence if $H'= H \cap G'$ then  $ H(\lambda)  = H ' = H(\zeta')$. 
This, along with Lemma \ref{l.e1},  implies that $\psi'$ extends to $H(\psi, \lambda) = H'(\psi)$ whenever 
$\psi$ extends to $H(\psi)$. But $H(\psi, \lambda) = H (\psi', \lambda) = H' (\psi')$, 
 since $\psi'$ is the $\lambda$-Clifford correspondent of $\psi$. 
This implies  the lemma in the case that $T'$ is a direct linear reduction. 
If $T'$ is a linear reduction, then repeated applications of the above argument completes 
the proof of the lemma.
\end{proof}

\begin{lemma}\mylabel{l.e2}
Let $Q$ be a normal  $q$-subgroup of a finite group $G$, for some prime $q$.
Assume further  that $E \unlhd G$ with $E \leq Q$, while
$\lambda \in \Irr(E)$
and   $\chi \in \Irr(Q| \lambda)$.
 If $A$ is a $q'$-subgroup of $G$ then there exists a
  $Q$-conjugate $\lambda_1 \in \Irr(E)$ of $\lambda$
that is $A(\chi)$-invariant and lies  under $\chi$.
\end{lemma}

 \begin{proof}
 Clifford's Theorem implies that $\chi$ lies above the $Q$-conjugacy
class of $\lambda$.
The $\pi$-group $A(\chi)$ fixes
 $\chi$, and normalizes $E$,  as the latter is normal in $G$. Hence
 $A(\chi)$ permutes among
themselves the $Q$-conjugates of $\lambda$. As  $(|A(\chi)|, |Q|) = 1$,
Glauberman's Lemma (Lemma 13.8 in \cite{is}) implies that $A(\chi)$ fixes at least
one character  $\lambda_1$ of the $Q$-conjugates of $\lambda$.
\end{proof}

The above lemma implies
\begin{proposition}\mylabel{r.e1}
Let $Q$ be a normal $q$-subgroup of $G$, while $A$ is any  $q'$-subgroup   of $G$. If $\chi \in \Irr(Q)$ 
we write $T$ for the triple $T = (G, Q, \chi)$.  Then there exists  a chain of linear subtriples $T_0, T_1, \dots, T_n$
of $T$ starting with $T_0 = T$ and ending with a linear limit $T_n = T'$ of $T$, so that for all $i=0, 1, \dots, n$,
the central character $\zeta^{(T_i)}$ of $T_i$  is $A(\chi)$-invariant. 
 We call such a $T'$, an $A(\chi)$-invariant linear limit of $T$. 
\end{proposition}

\begin{proof}
 Let  $T_1= (G_1,  Q_1, \chi_1)$ be a direct linear reduction   of $T= (G, Q, \chi)$. 
Then $T_1 = T(\lambda_1)$, where $\lambda_1 $ is a linear character of a normal subgroup $E_1$ with 
$ E_1 \leq Q$. According to Lemma \ref{l.e2}  there is  a $Q$-conjugate of $\lambda_1$ that is 
$A(\chi)$-invariant. Thus, without loss, 
  we may assume that $\lambda_1$ is $A(\chi)$-invariant.
So  $A(\chi) \leq G_1 = G(\lambda_1)$. Because $A(\chi)$ fixes $\lambda_1$,  while
$\chi_1 \in \Irr(Q_1)$ is the $\lambda_1$-Clifford correspondent of $\chi$ in $Q_1= Q(\lambda_1)$,
we get $A(\chi)(\chi_1)= A(\chi)$. Let $Z_1 $ be the central subgroup of $T_1$, so $E_1 \leq Z_1 \leq Q_1$. 
If  $\zeta_1  \in \Lin(Z_1)$ is the central character of $T_1$, then $\zeta_1$ lies above $\lambda_1$ 
while $\chi_1 |_{Z(T_1)}= \chi_1(1) \cdot \zeta_1$. Hence $A(\chi) = A(\chi, \chi_1)$ also fixes $\zeta_1$.

Let  $T_2= (G_2,  Q_2, \chi_2)$ be  a direct linear reduction of $T_1$.
So $T_2 = T_1 (\lambda_2)$, for some linear character $\lambda_2 \in \Irr(E_2)$, where $E_2 \unlhd G_1$. 
Since $A(\chi) = A(\chi, \chi_1)$, Lemma \ref{l.e2} implies that we can pick $\lambda_2$ to be $A(\chi)$-invariant.
Hence the $\lambda_2$-Clifford correspondent $\chi_2$ of $\chi_1$ is also $A(\chi)$-invariant.
Therefore, the central character $\zeta_2$ of $T_2$ is also $A(\chi)$-invariant, while $A(\chi)= A(\chi, \chi_2)$.
We proceed similarly at every direct linear reduction until we reach $T'$.
\end{proof}

\begin{proposition}\mylabel{l.e4}
Let $G$ be  finite group of odd order such that $G =N K $,  where
 $N$ is a normal subgroup of $G$  and $(|G/ N|, |N|) = 1$. Let $H  = N \cap K$ and let $\theta$ be any irreducible
$K$-invariant  character of $H$  that induces an  irreducible character $\theta^N$ of $N$. 
Then $\theta$ has a unique canonical extension, $\theta^e$, to $K$ such that 
$(|K/ H|, o(\theta^e)) = 1$ (where $o(\theta^e) $ is the determinantal order of $\theta^e $, see 
p. 88 in \cite{is}).  Also $\theta^N$  has a unique canonical extension, $(\theta^N)^e$, to $G$  such that 
$(|G/ N|, o((\theta^N)^e)) = 1$.  Furthermore, $\theta^e$   induces 
$$
(\theta^e) ^G = (\theta ^N) ^e.
 $$
\end{proposition}

\begin{proof}
Let $\pi$ be the set of primes that divide $|N|$.  Then
$|K/H|= |G/N|$ is a $\pi'$-number,   and thus is coprime to $|H|$.
Because $\theta \in \Irr(H)$ is $K$-invariant,  there exists a unique extension $\theta^e$ to $K$
such that 
\begin{equation}\mylabel{pr.e0}
o(\theta)= o(\theta^e),
\end{equation}
by Corollary 6.28 in \cite{is}.

According to Corollary 4.3 in \cite{isa3},  induction defines a bijection \,$\Irr(K |\theta) \to \Irr(G |\theta^N)$.
Therefore,
\begin{equation}\mylabel{pr.e1}
\chi:=( \theta^e)^G  \in \Irr(G | \theta^N).
\end{equation}
But $\theta^N$ is $G$-invariant since  $\theta$ is $K$-invariant  and $G= NK$. This, and the fact that $(|N|, |G/N|)=1$,
 implies that $\theta^N$ 
extends to $G$.  
Let $\Psi = (\theta^N)^e \in \Irr(G)$ be the unique extension of $\theta ^N$
such that  $o(\Psi)= o(\theta^N)$
is a $\pi$-number. Since $\chi$ lies above $\theta^N$, Gallagher's theorem 
 (see Corollary 6.17 in \cite{is}) implies  that
$$\chi= \mu \cdot \Psi,$$
for some $\mu \in \Irr(G/N)$. We compute the degree $\deg(\chi)$ in two ways.
First
$$
\deg(\chi)= \deg(\mu) \cdot \deg(\Psi)= \deg(\mu) \cdot \deg(\theta^N)= \deg(\mu)  \cdot |N:H| \cdot \deg(\theta).
$$
As $\chi= (\theta^e)^G$ we also have that
$$
\deg(\chi) = |G:K| \cdot \deg(\theta^e)= |G:K|\cdot  \deg(\theta) = |N:H| \cdot \deg(\theta).
$$
We conclude that $\deg(\mu)=1$. Thus $\mu \in \Lin(G/N)$.
Therefore
\begin{equation}\mylabel{pr.e2}
\det(\chi) = \mu^{\Psi(1)} \det(\Psi).
\end{equation}

We can now compute $o(\chi) $ in two ways.  First,
$o(\Psi) = o(\theta^N) $  and  $\Psi(1) = \theta^N(1)$ are  $\pi$-numbers.
Since $\mu \in \Irr(G/N)$,  we get that $o(\mu)$ is a $\pi'$-number.
Therefore,  \eqref{pr.e2} implies that  the $\pi'$-number $o(\mu)$ divides $o(\chi)$.

On the other hand, \eqref{pr.e1} and  Lemma 2.2 in \cite{isa4} imply that
$$
o(\chi) = o((\theta^e)^G)  \text{ divides }  2\cdot  o(\theta^e).
$$
As $G$ has odd order, we get that
$o(\chi)$ divides $o(\theta^e) $. In view of \eqref{pr.e0},  we have $o(\theta^e) =
o(\theta)$, while $ o(\theta)  \bigm |   |H|$. We conclude that
$ o(\chi) $ is a $\pi$-number.

Hence  the only way  the $\pi'$-number $o(\mu)$ can divide $o(\chi)$,
 is if $o(\mu) = 1$.  So $\mu=1$,   and
$$
(\theta^e)^G = \chi = \Psi  = (\theta^N)^e,
$$
as desired.
\end{proof}

\begin{lemma}\mylabel{l.e5}
Let $ E \unlhd Q$ be subgroups of a finite group $G$, and 
$\zeta , \beta$ be irreducible characters of $E$ and $ Q$, respectively, with $\beta$ lying above $\zeta$.
Let $T$ be  any subgroup of $G$ that normalizes both  $E$ and $Q$. 
If  $Q \leq T \leq G(\beta)$ then  $T = T(\zeta) \cdot Q$.  If in addition $T$ is a Sylow subgroup of $G(\beta)$
then $T(\zeta)$ is a Sylow subgroup of $G(\beta, \zeta)$.  
\end{lemma}

\begin{proof}
Since $T$ fixes $\beta$, it permutes among themselves the $Q$-conjugacy
 class of characters  in $\Irr(E)$ lying under $\beta$. So $T \leq T(\zeta) \cdot Q$.
The other inclusion is trivial.

Now assume that $T$ is a $q$-Sylow subgroup of $G(\beta)$, for some prime $q$. 
Since $G(\beta) =G(\beta, \zeta)Q$ with $G(\beta,\zeta) \cap Q = Q(\zeta)$,
 the index $[G(\beta) : G(\beta,\zeta)]$ equals the index $[Q : Q(\zeta)]$.
 Similarly, $[Q : Q(\zeta)] = [T :T(\zeta)]$. Hence the index of $T(\zeta) =  G(\beta,\zeta) \cap T$ in $G(\beta, \zeta)$ 
is a $q'$-number, which implies that  $T(\zeta)$ is a Sylow $q$-subgroup of $G(\beta,\zeta)$.
\end{proof}

The following is a well known result that we will use
 repeatedly  below.
\begin{lemma}\mylabel{ext}
Let $M = A \ltimes B$  be a finite solvable  group where $(|A|, |B|) = 1$. 
Let $\beta \in \Irr(B)$ be  $M$-invariant.
Then there is a unique canonical extension $\beta^e$ of $\beta$ to $M$.
 Furthermore,  for any irreducible character $\chi$ 
of $M$ lying above $\beta$  there exists a unique 
irreducible character $\alpha$ of $A$ such that $\chi = \alpha \cdot  \beta^e$, 
where 
$$
\chi(x \cdot  y) = \alpha(x) \cdot \beta^e(xy),
$$ 
 for all $x \in A$  and $ y \in B$. We will write $\chi = \alpha \ltimes \beta$, to denote the
 above product. 
\end{lemma}

\begin{proof}
The existence of a unique canonical extension  $\beta^e$ of $\beta$  follows 
from  Corollary 6.28 in \cite{is}. The rest of the lemma follows from 
  Corollary 6.17 in \cite{is}.
\end{proof}

\begin{defn}\mylabel{prd}
If $M =  A \ltimes B$  where $(|A|, |B|)= 1$,  and 
$\alpha \in \Irr(A)$ while $1_B$  is the trivial character of $B$, 
 we write $\alpha \ltimes  1_B$   (or simply $\alpha \ltimes 1$ if $B$ is clear), 
for  the unique irreducible 
character of $M$ defined as $\alpha \ltimes  1 (x \cdot y ) = \alpha (x)$ , for all $x \in A $ and
 $y \in B$. Notice that if $M \leq G$ then $G(\alpha ) \leq G(\alpha \ltimes  1)$, while 
$G(\alpha) = G(\alpha \ltimes  1)$   if $G$ normalizes $A$.Furthermore, if $M \unlhd G$ then Frattini's Argument implies that  $G(\alpha \ltimes  1) = G(\alpha) \cdot B$.
\end{defn}

In view of the above definition we slightly generalize Proposition \ref{r.e1} to 
\begin{proposition}\mylabel{finv}
Let $M= P \ltimes B$ be a normal subgroup of $G$, where $P$ is a $p$-group and $B$ is a  normal 
 $p'$-subgroup of $G$.
Let $\chi  \in \Irr(M)$ where  $\chi = \alpha \ltimes 1$ and $\alpha \in \Irr(P)$. 
Let $T$ be the triple  $T= (G, M, \chi)$.
If $A$ is any $p'$-subgroup of $G$ then 
there exists an $A(\chi)$-invariant linear limit $T'$ of $T$, i.e., there exists a chain of linear
 subtriples  $T_i= (G_i, M_i , \chi_i) $, for $i=0, \dots, n$
 starting with $T =T_0$ and ending with the linear limit $T_n = T'$
 of $T$ so that for all $i=0, 1, \dots, n$, the central character $\zeta_i$ of $T_i$ is 
$A(\chi)$-invariant. Furthermore, $T_i = T_{i-1}(\lambda_i)$ for all $i=1, \dots, n$ where 
$\lambda_i$ is an $A(\chi)$-invariant linear character of $L_i$ lying under $\chi_{i-1}$ and 
above $\zeta_{i-1}$, and $L_i$ is a normal subgroup of $G_{i-1}$  with
 $B \leq Z(T_{i-1}) \leq L_i \leq M_{i-1}$.
\end{proposition}

\begin{proof}
Let $\bar{T} $ be the triple $\bar{T}= (\bar{G}, \bar{M}, \bar{\chi})$, where $\bar{G}= G/ B$, 
$\bar{M}= M/ B \cong P$ and $\bar{\chi} \in \Irr(\bar{M})$ is the unique irreducible character of $\bar{M}$ 
that inflates to $\chi \in \Irr(M)$. We also write $\bar{A}$ for the quotient group 
$\bar{A} = (A(\chi) B) / B$. Clearly $\bar{A}$ fixes $\bar{\chi}$, and its order is coprime to that 
of $\bar{M}$. Hence Proposition \ref{r.e1} implies that there exists a chain of
$\bar{A}$-invariant  linear subtriples  $\{ \bar{T}_i\} _{i=0}^n$, so that $\bar{T}_0= \bar{T}$, 
$\bar{T_i}$ is a direct linear reduction of $\bar{T}_{i-1}$, while  $\bar{T}_n$ is a linear limit of $\bar{T}$ 
and the central character $\bar{\zeta}_i$ of each $\bar{T}_i$ is $\bar{A}$-invariant.

Let  $\bar{T}_i= (\bar{G}_i, \bar{M}_i, \bar{\chi}_i)$, for $i=0, 1, \dots, n$.  
Then  $\bar{G}_i = G_i/ B$ , $\bar{M}_i = M_i/ B$, where $G_i$ and $M_i$ are subgroups 
of $G$ with $M_i \leq M$. In addition,  $\bar{\chi}_i $ inflates to a unique irreducible 
character $\chi_i$ of $M_i$.
 So $T_i= (G_i, M_i, \chi_i)$ is a subtriple of $T$ for all such $i$.
Furthermore, it is easy to see (with the use of Proposition \ref{char.z})
that if $Z(\bar{T}_i)$  is the center of $\bar{T}_i$ then $Z(\bar{T}_i)= Z_i / B$ where $Z_i$ is  
the center  of $T_i$. In addition, the central character $\zeta^{(\bar{T}_i)} \in \Irr(Z(\bar{T}_i)$
 of $\bar{T}_i$ 
is inflated from the central character $\zeta_i \in \Irr(Z_i)$ of $T_i$.
Hence the fact that $\bar{\zeta}_i$ is $\bar{A}$-invariant implies that $\zeta_i$
 is $A(\chi)$-invariant, for all $i=1, \dots, n$.

 Even more, $T_i$ is a direct linear reduction of $T_{i-1}$, for all $i=1, \dots, n$.
Indeed, the fact that $\bar{T}_i$ is a direct linear reduction of $\bar{T}_{i-1}$ implies that 
$\bar{T}_i = \bar{T}_{i-1}(\bar{\lambda}_i)$, where $\bar{\lambda}_i$ is a linear character of 
the normal subgroup $\bar{L}_{i}$ of $\bar{G}_{i-1}$ that  satisfies  
$Z(\bar{T}_{i-1}) \leq \bar{L}_i \leq \bar{M}_{i-1}$, while $\bar{\lambda}_i$ lies under 
$\bar{\chi}_{i-1}$. Hence $\bar{L}_i= L_i/ B$, where $L_i$ is a normal subgroup of $G_{i-1}$ with 
$Z_{i-1} \leq L_i \leq M_{i-1}$ and $\bar{\lambda}_i$ inflates to a unique linear 
 character $\lambda_i$ of $L_i$ that lies under $\chi_{i-1}$. Now it is easily to see that  
$\bar{G}_{i} = \bar{G}_{i-1} (\bar{\lambda}_i)= G_{i-1}(\lambda) / B$.
So $G_i = G_{i-1}(\lambda)$ and similarly, $M_i= M_{i-1}(\lambda)$. Because $\bar{\chi}_{i}
\in \Irr(\bar{M}_i)$ induces $\bar{\chi}_{i-1}$ to $\bar{M}_{i-1}$, we have that 
$\chi_{i} \in \Irr(M_i)$ induces $\chi_{i-1}$ to $M_{i-1}$. Also $\chi_{i}$ lies above $\lambda_i$ 
and thus $\chi_i$ is the $\lambda_i$-Clifford correspondent of $\chi_{i-1}$, for all $i=1, \dots, n$. 
We conclude that $T_i= T_{i-1}(\lambda_i)$, for all $i=1, \dots, n$. 
In addition,   $\lambda_i$ is $A(\chi)$-invariant since $\bar{\lambda}_i$ is $\bar{A}$ invariant for all such $i$.

To see that $T_n$ is a linear limit of $T$ observe that  the procedure described above works  both ways.
So any direct linear reduction   of $T_n$ determines  a direct linear reduction  of $\bar{T}_n$. As the latter triple
 is irreducible, we have that $T_n$ is a linear limit of $T$. 
\end{proof}

We first prove:
\begin{theorem}\mylabel{e0}
Assume that  $G$ is a  finite group. 
Let $L\leq M \leq N$ be  normal subgroups of $G$,  so that  $N/ M, M/ L$ and $L$ 
are nilpotent group while $N$ has order $p^a q^b$ for distinct odd primes  $p, q$.
So $L = L_p \times Q$ where $Q$ is the  $q$-Sylow and $L_p$ the $p$-Sylow group subgroup of $L$.
Let  $M_p$ be a $p$-Sylow subgroup of $M$ and let  $H = M_p L$. Assume that   
$\phi \in \Irr(L)$ and $\theta  \in \Irr(H)$ lies above $\phi$. 
Let $\theta_{\phi}  \in \Irr(H(\beta))$ 
be the $\phi$-Clifford correspondent of $\theta$. Write $\phi = \phi_p \times \beta$, for $\phi_p \in \Irr(L_p)$ 
and $\beta \in \Irr(Q)$. 
 Assume further that 
 all irreducible characters of $G$ 
lying above $\beta$ are monomial, while $\beta$ extends to $G(\theta_{\phi})$. Then  
there exists a linear limit $(G', H', \theta')$ of $(G, H, \theta)$ so that 
$$
[N_q', H_p'] \leq K',
$$
where $K'$ is the kernel of the triple $(G', H', \phi')$, while $H_p'$ is a $p$-Sylow 
subgroup
of $H'$ and $N_q'$ is a $q$-Sylow subgroup of $ N' = G' \cap N$.
\end{theorem}

\begin{proof}
First note that $H= M_p L $ is  a normal subgroup of $G$, because $H/L$ is the $p$-Sylow subgroup of the 
nilpotent group $M/L$.  
We fix  $\theta \in \Irr(H)$ and $\phi = \phi_p \times \beta \in \Irr(L)$ satisfying
 the hypothesis of the theorem.
 We also pick a $p$-Sylow subgroup $A$  of $G$ so that $A(\beta)$  is a $p$-Sylow 
subgroup of $G(\beta)$.

Because $H/Q$ is a $p$-group, $H$ is the semidirect product
$H = (A\cap H) \ltimes Q$, where $A \cap H$ is a $p$-Sylow subgroup of $H$.
Let $P=(A\cap H)(\phi)$. Then $H(\phi) = P \ltimes Q$.
Furthermore, $\beta $ extends to $H(\phi) \leq H(\beta)$, since $(|H/Q|, |Q|)= 1$.
Hence, if $\theta_{\phi} \in \Irr(H(\phi))$ is the $\phi$-Clifford correspondent of
$\theta$, then (see Lemma \ref{ext}) there exists a unique irreducible character $\alpha \in \Irr(P)$ such that
$$
\theta_{\phi} = \alpha \ltimes \beta,
$$
Because $H \unlhd G$, Frattini's argument implies that $G(\theta_{\phi}) = G(\alpha, \beta) Q$.
Note also that $\alpha$ lies above $\phi_p \in \Irr(L_p)$ as $\theta_{\phi}$  lies above $\phi$. 
Hence $G(\alpha) \leq G(\phi_p) \leq G$.

At this point we pick a $q$-Sylow   subgroup $B$ of $G$ so that 
$B(\alpha)$ is a $q$-Sylow   subgroup of $G(\alpha)$ while  $B(\alpha, \beta)$ is a 
$q$-Sylow subgroup of $G(\alpha, \beta)$.

Let $S= (G, Q, \beta)$.
Because $A(\beta)$ has coprime order to that of $Q$, Proposition  \ref{r.e1}
implies that we can get an $A(\beta)$-invariant  linear limit $S_1= (G_1, Q_1, \beta_1)$
 of  $S$. Hence the central character $\zeta_1 $ of $S_1$ is $A(\beta)$-invariant.
Furthermore, $G_1= G(\zeta_1)$ while $\beta_1$ is the unique character of $  Q_1 = Q(\zeta_1)$ 
lying above $\zeta_1$,  by Proposition \ref{l7}.
We write $L_1 = L \cap G_1$,  $H_1 = G_1\cap H = H(\zeta_1), M_1 = M \cap G_1= M(\zeta_1)$ and 
$N_1 =  G_1 \cap N = N(\zeta_1)$.
Observe that $L_1 = L_p \times Q_1$, while $N_1/ M_1$ is a nilpotent group and $M_1/ H_1$ is a $q$-group. 
Let $\theta_1\in \Irr(H_1)$ be  the $S_1$-reduction  of $\theta$.
So $\theta_1$  lies above $\phi_1$  the $S_1$-reduction of $\phi$. 
Note that $\phi_1 = \phi_p \times \beta_1$. So $\theta_1$ lies above 
 $\beta_1 $ and thus above $\zeta_1 $, and induces $\theta $.
Clearly all the irreducible characters of $G_1$ that lie above $\beta_1$ ( and equivalently above $\zeta_1$) 
are still monomial  (see Proposition \ref{l4}).

Since $G_1 = G(\zeta_1)$ we have $A(\beta) \leq G_1$.  Also $G_1 = G(\zeta_1 ) = G(\zeta_1, \beta_1)$
 by Proposition \ref{l1.5}.
Thus $G_1 \leq G(\beta)$, since $\beta_1$ induces $\beta$.  Also $H_1 \leq H(\beta)$. 
The group $A(\beta)$ is a $p$-Sylow subgroup of $G(\beta)$ contained in $G_1$. 
Hence $A(\beta)$ is a $p$-Sylow subgroup of $G_1$.  Furthermore, 
  the $p$-Sylow subgroup $P$ of $H(\beta)$ is contained in $H_1$,
 and thus $P$  is also a $p$-Sylow subgroup of $H_1$.
So $H_1 = P \ltimes Q_1$.
Because $\beta$ extends to $G(\theta_{\phi})= Q G(\alpha, \beta)$, 
 Lemma \ref{l.e1} implies that $\beta_1$ extends to $(Q G(\alpha, \beta)) (\zeta_1, \beta_1) $.
Hence $\beta_1$ extends to $Q_1 G_1(\alpha)$ (note that $G_1(\beta) = G_1(\beta_1) = G_1$).
 Let $\beta_1^e$ be the  unique  canonical extension of $\beta_1 $ to $H_1$.
Then there exists a unique irreducible character $\alpha' $ of $P$ so that $\theta_1 =
\alpha \ltimes \beta_1= 
\alpha' \cdot \beta_1^e$. It is not hard to see that $\theta_1^{H(\beta)} = \alpha'
\cdot (\beta_1^e)^{H(\beta)}$ (see Exercise 5.3 in \cite{is}). Since $\beta_1$
induces $\beta \in \Irr(Q)$,  Lemma  \ref{l.e4}  implies that 
$(\beta_1^e)^{H(\beta)} = \beta^e$.
So $\theta_1^{H(\beta)}= \alpha' \cdot \beta^e$. But  $\theta_1^{H(\beta)}$
lies above $\beta$ and induces $\theta  \in \Irr(H)$ (because $\theta_1$ does). Hence $\theta_1^{H(\beta)}= 
\theta_{\beta}$, and thus $\alpha' = \alpha$. Therefore, 
\begin{equation}\mylabel{et}
\theta_1 = \alpha \ltimes \beta_1, 
\end{equation}
which implies that $G_1(\theta_1) = G_1(\alpha, \beta_1) Q_1 = G_1(\alpha) Q_1$.

Let $K_1 = \Ker(\zeta_1)$ and   $\br{G}{1}= G_1/ K_1$, $\br{H}{1} = H_1/ K_1$ and  
$\br{Q}{1}= Q_1/ K_1$.
Then $\br{Q}{1} \leq \br{H}{1} \leq \br{N}{1}$ are normal subgroups of $\br{G}{1}$. 
Furthermore,  the irreducible characters 
$\beta_1, \theta_1$  of $Q_1$ and $  H_1$  are inflated  from unique irreducible characters 
$\br{\beta}{1}$ and $ \br{ \theta}{1}$   of $\br{Q}{1}$ and $ \br{H}{1}$, all respectively.  
If $Z_1$ is the center of the triple $S_1= (G_1, Q_1, \beta_1)$
then we write $\br{Z}{1}$ for the quotient group $Z_1/ K_1$. We also write  $\br{\zeta}{1}$ for 
the unique irreducible 
character of $\br{Z}{1}$ that inflates to the central character $\zeta_1 \in \Irr(Z_1)$ of $S_1$. 
According to Proposition \ref{l1}, the cyclic  group $\br{Z}{1}$ is the center of $\br{Q}{1}$, affords the faithful 
$\br{G}{1}$-invariant linear character $\br{\zeta}{1}$, and it is maximal among the abelian normal subgroups of
$\br{G}{1}$ contained in $\br{Q}{1}$. So $\bar{Q}_1$ satisfies Condition X.
Since  every irreducible character of $G_1$ lying above $\zeta_1$ is monomial, 
we get that every irreducible character of $\br{G}{1}$ lying above $\br{\zeta}{1}$, and equivalently above 
$\br{\beta}{1}$,   is also  monomial.
Because the $q$-special character $\br{\zeta}{1}$ is $\br{G}{1}$-invariant,  there exists a $q$-special
 character $\br{\chi}{1}$ of 
$\br{G}{1}$ lying above $\br{\zeta}{1}$.  So $\br{\chi}{1}$  is  monomial.
Let $\bar{P}= (P K_1)/ K_1$. Then 
$\bar{P}\cong P$ while $\bar{P} \cdot \br{ Q}{1} = \br{H}{1}  \unlhd \br{G}{1}$.
Hence we can apply Theorem B to the groups $\br{Q}{1}, \bar{P}$ and $\br{G}{1}$ in the place of $P, S$ 
and $G$ respectively.
 We conclude that $\bar{P}$ centralizes $\br{Q}{1}$, and thus $\bar{H}_1$ is a nilpotent group. Furthermore, if 
$\bar{\alpha}$ is the unique irreducible character of $\bar{P}$ that inflates
 to $\alpha \ltimes 1\in \Irr(P K_1)$, then 
$$
\br{H}{1} = \bar{P} \times \br{Q}{1} \text{ and } \br{\theta}{1} = \bar{\alpha} \times \br{\beta}{1}. 
$$
So  both $\bar{P}$ and $\br{Q}{1}$ are normal subgroups of $\br{G}{1}$. Also,    
 it is easy to see that $G_1(\alpha) / K_1 = \br{G}{1}(\bar{\alpha})$, because $\br{G}{1}$ normalizes 
$\bar{P}$ and thus $\br{G}{1}(\bar{\alpha}) \leq G_1(\alpha \ltimes 1)/ K_1 \leq G_1(\alpha)/ K_1$ (the
 other inclusion holds trivially).
Because $\beta_1$ extends to $Q_1G_1(\alpha)$, we get that $\br{\beta}{1}$ extends to 
$(Q_1G_1(\alpha))/ K_1$. So $\br{\beta}{1}$  extends to $\br{G}{1} (\bar{\alpha})$ 
(note that $\br{Q}{1}$ centralizes $\bar{P}$ and
 thus it is a subgroup of $\br{G}{1} (\bar{\alpha})$).

Take  $\br{E}{1}$  to be the  triple $\br{E}{1} = (\br{G}{1}, \bar{P}, \bar{\alpha})$.
If $B_1= B \cap G_1 = B(\zeta_1)$, then $\br{B}{1} = B_1/ K_1$ is a $q$-subgroup of $\br{G}{1}$. 
Thus we can form  a  $\br{B}{1}(\bar{\alpha})$-invariant   linear limit $\br{E}{2}$ of $\br{E}{1}$.
So $\br{E}{2}= (\br{G}{2}, \br{P}{2}, \br{\alpha}{2})$, where $\br{P}{2} \leq \bar{P}$ 
and $\br{\alpha}{2}$ induces $\bar{\alpha} $ to $\bar{P}$.
Hence $\br{P}{2} = P_2 / K_1$ with $P_2 \leq P$, while $\br{\alpha}{2} $ inflates 
to a unique irreducible character $\alpha_2 $ of $P_2$ that induces $\alpha$ to $P$.
In addition, $\br{G}{2} = G_2/ K_1$ where $G_2$ is a subgroup of $G_1$.
  The central character $\br{\zeta}{2}$ of $\br{E}{2}$  is $\br{B}{1}(\bar{\alpha})$-invariant. 
Also  $\br{G}{2} = \br{G}{1}(\br{\zeta}{2})$  and $\br{P}{2} = \bar{P}(\br{\zeta}{2})$. 
The center $\br{Z}{2}$ of $\br{E}{2}$ clearly contains $K_1$, even more the kernel $\br{K}{2}=  \Ker(\br{\zeta}{2})$  
of $\br{\zeta}{2}$
contains $K_1$.  According to  Propositions \ref{l1.5} and \ref{l7} the character 
$\br{\alpha}{2}$ is $\br{G}{2}$-invariant and it is  the 
only character of $\br{P}{2}$ that lies above $\br{\zeta}{2} \in \Irr(\br{Z}{2})$.
So 
\begin{equation}\mylabel{bga}
\br{G}{2}= \br{G}{2}(\br{\alpha}{2}) = \br{G}{2}(\bar{\alpha}) \leq \br{G}{1}(\bar{\alpha}), 
\end{equation}
where the last equality follows from  the fact that $\br{\alpha}{2}$ 
is the only  irreducible character
of $\br{P}{2}$ that lies above $\br{\zeta}{2}$ and induces $\bar{\alpha}$ to $\bar{P}$. 

We next observe that the $\br{E}{2}$-reductions leave $\br{Q}{1}$ unchanged, i.e., $\br{Q}{1}(\br{\zeta}{2}) = \br{Q}{1}$, 
because $\br{Q}{1}$ centralizes $\bar{P}$. Hence $\br{Q}{1} \leq \br{B}{1}( \bar{\alpha}) \leq  \br{G}{2}$.
Furthermore, the $\br{E}{2}$-reduction of $\br{H}{1}= \bar{P} \times  \br{Q}{1}$  equals  $\br{H}{2} = \br{P}{2} 
\times \br{Q}{1}$.
 In addition, the $\br{E}{2}$ reduction of $\br{\theta}{1} \in \Irr(\br{H}{1})$ 
equals $\br{\theta}{2}  = \br{\alpha}{2} \times \br{\beta}{1}$.
Note that since $\br{\beta}{1}$ extends to $\br{G}{1}(\bar{\alpha})$, equation \ref{bga} implies that 
$\br{\beta}{1}$ extends to $\br{G}{2}$.
Because every irreducible character of $\br{G}{1}$ lying above $\br{\beta}{1}$  is monomial, Proposition \ref{l4}
implies that  every irreducible character of 
$\br{G}{2}$ lying above $\br{\theta}{2}$ is still monomial.

We look at the quotient group $\br{G}{2}/ \br{K}{2}$.
 Its subgroup  $\br{H}{2}/ \br{K}{2}$ is a nilpotent normal subgroup, and
splits as  $\br{H}{2}/ \br{K}{2} = \br{P}{2}/ \br{K}{2} \times  (\br{Q}{1} \br{K}{2}) / \br{K}{2}$, 
where $(\br{Q}{1} \br{K}{2}) / \br{K}{2}$, is naturally 
isomorphic to $\br{Q}{1}$.   We identify these two isomorphic groups and we  consider $\bar{\beta}_1$ as 
an irreducible character 
of $(\br{Q}{1} \br{K}{2}) / \br{K}{2}$. As we have seen $\bar{\beta}_1$ extends to $\br{G}{2}$, 
hence $\bar{\beta}_1$ (considered as a character of $(\br{Q}{1} \br{K}{2}) / \br{K}{2}$
extends to $\br{G}{2}/ \br{K}{2}$.  
Let $\hat{\beta}_1 \in \Irr(\br{G}{2}/ \br{K}{2}) $  be such an extension.  Then $\hat{\beta}_1$ is 
a $q$-special character of $\br{G}{2}/ \br{K}{2}$.
The irreducible character $\br{\alpha}{2}$ of $\br{P}{2}$ 
is inflated from a unique irreducible character $\bar{\alpha'}_2$ of 
the $p$-group $\br{P}{2}/ \br{K}{2}$.  Because $\br{\alpha}{2}$ is 
$\br{G}{2}$-invariant,  $\bar{\alpha'}_2$ is a $\br{G}{2}/ \br{K}{2}$-invariant 
$p$-special character of $\br{P}{2}/ \br{K}{2}$. Therefore   there 
exists a $p$-special irreducible character  $\hat{\alpha}_2$ of
 $\br{G}{2}/ \br{K}{2}$ lying above $\bar{\alpha'}_2$. Then the product 
$\hat{\chi} = \hat{\alpha}_2 \cdot \hat{\beta}_1$ is an irreducible character of 
$\br{G}{2}/ \br{K}{2}$ that lies above $\bar{\alpha'}_2 \times \bar{\beta}_1  \in \Irr(\br{H}{2}/ \br{K}{2})$. 
 In addition $\hat{\chi}_q = \hat{\beta}_1(1)$.
 Observe also that $\hat{\chi}$ is monomial, because every irreducible character of $\br{G}{2}$ lying above  
$\br{\theta}{2}$ is monomial, while $\bar{\alpha'}_2 \times \bar{\beta}_1$ inflates to $\br{\theta}{2} \in \Irr(\br{H}{2})$.
Furthermore, Proposition \ref{l1} implies that 
$\br{Z}{2} / \br{K}{2}$ is the center of $\br{P}{2}/ \br{K}{2}$, it is maximal among the abelian normal subgroups of 
$\br{G}{2}/ \br{K}{2}$ 
contained in $\br{P}{2}/ \br{K}{2}$  and it affords a faithful $\br{G}{2}/ \br{K}{2}$-invariant  irreducible 
 character that inflates 
to $\br{\zeta}{2} \in \Irr(\br{Z}{2})$.
We conclude that all the hypothesis of Theorem B are satisfied for the groups $\br{G}{2}/ \br{K}{2}, \br{P}{2}/ \br{K}{2}$
  and $(\br{Q}{1} \br{K}{2})/ \br{K}{2}$ in the place of $G, P$ and $S$, respectively. 
Hence  any $p'$-subgroup $\mathcal{Q}$ of $\br{G}{2}/ \br{K}{2}$ 
 centralizes $\br{P}{2}/ \br{K}{2}$ 
provided that $\mathcal{Q} \cdot \br{P}{2}/ \br{K}{2}$ is a normal subgroup of $\br{G}{2}/ \br{K}{2}$.

Let $\br{N}{2}$ be the $\br{E}{2}$-reduction of $\br{N}{1}$  (that is $\br{N}{2} = \br{N}{1} (\br{\zeta}{2})$ ), 
and $\br{M}{2}$ the 
$\br{E}{2}$ reduction of $\br{M}{1}$. Then $\br{N}{2}/ \br{M}{2}$ is still a nilpotent group, 
while $\br{M}{2}/ \br{H}{2}$ is a $q$-group.
  Let  $U$ be a  $q$-Sylow subgroup of $\br{N}{2}$,  then 
$U \br{M}{2}  = U \br{H}{2}$ is a normal subgroup of $\br{G}{2}$, because $(U \br{M}{2} ) / \br{M}{2}$ 
is the $q$-Sylow subgroup of $\br{N}{2}/ \br{M}{2}$.
 So $(U \br{H}{2}) / \br{K}{2} \unlhd \br{G}{2}$.
Furthermore $(U \br{H}{2})/ \br{K}{2} = (U \br{K}{2}) / \br{K}{2} \ltimes \br{P}{2}/ \br{K}{2}$.
Hence the $q$-group $(U \br{K}{2}) / \br{K}{2}$ centralizes  $\br{P}{2}/ \br{K}{2}$, i.e, 
$[U, \br{P}{2}] \leq \br{K}{2}$. So 
$$
[U, P_2] \leq K_1 K_2.
$$

It is easy to see that the triple $(G_2, P_2 \ltimes Q_1, \alpha_2 \ltimes \beta_1)$ is a linear reduction of
 $(G, P\ltimes Q, \theta_{\phi})$. 
This completes the proof of the  theorem.
\end{proof}

For later use we observe that what we actually proved in  Theorem \ref{e0} is
\begin{remark}\mylabel{inv.limit}
Assume the hypothesis and the notation of Theorem \ref{e0}. Then any 
 $A(\beta)$-invariant limit $S_1=(G_1, B_1, \beta_1)$
 of $S= (G, B, \beta)$ makes $(G_1 \cap H)/ K_1$  nilpotent,
 where $A$ is a $p$-Sylow subgroup
 of $G$ and $A(\beta)$ is a $p$-Sylow subgroup of $G(\beta)$ and $K_1 $ is the kernel of $S_1$.
Furthermore, $H_1 = P \ltimes Q_1$, where $P$ is a $p$-Sylow subgroup of $H(\phi)$. 
Let $\br{E}{1}$ be the factor triple $\br{E}{1}= (\br{G}{1}, \bar{P}, \bar{\alpha})$, where $\br{G}{1}= G_1/ K_1$ 
and $\bar{P}= (PK_1) / K_1 \cong P$. If  
  $B$ is  a $p'$-Hall subgroup 
of $G$  so that $B(\alpha)$  is a $p'$-Hall subgroup of $G(\alpha)$  and $B(\alpha, \beta)$
 is a $p'$-Hall subgroup of $G(\alpha, \beta)$, let $\br{B}{1} = (B \cap G_1) / K_1$. Then 
 any $\br{B}{1}(\bar{\alpha})$-invariant 
linear limit $\br{E}{2}= (\br{G}{2}, \br{P}{2}, \br{\alpha}{2})$
 of $\br{E}{1}$ satisfies
 $[(\br{N}{1}  \cap \br{G}{2})_q, \br{P}{2}] \leq \br{K}{2}$, 
where $\br{K}{2}$  is the kernel of $\br{E}{2}$ while 
$\br{N}{1}= (N\cap G_1) / K_1$ and $\br{H}{1}= H_1/ K_1$.
\end{remark}

\begin{corollary}\mylabel{ec0}
Assume the hypothesis of Theorem \ref{e0}. Let $B$ be a $q$-Sylow subgroup of $G$ so that 
$B(\alpha)$ is a $q$-Sylow subgroup of $G(\alpha)$ while $B(\alpha, \beta)$ is a $q$-Sylow subgroup of 
$G(\alpha, \beta)$. Assume further that $A$ is a $p$-Sylow subgroup of $G$ with $A(\beta)$ being 
a $p$-Sylow subgroup of $G(\beta)$, while $P = (A\cap H)(\beta)$. Let 
 $\zeta_1$ be  the central character of 
the $A(\beta)$-invariant direct linear reduction
$S_1= (G_1, Q_1, \beta_1)$ of $S= (G, Q, \beta)$. Then $G(\beta) = G(\beta, \zeta_1) \cdot Q$.
If in addition $C = C_Q(P)$ then $G(\alpha, \beta ) = G(\alpha, \beta, \zeta_1) \cdot C$. 
Thus 
 $B(\alpha, \beta, \zeta_1)\cdot K_1= B_1(\alpha ) K_1$ is a $p'$-Hall  subgroup 
of both  $G_1(\alpha) \cdot K_1$. 
\end{corollary}

\begin{proof}

Since $S_1$ is an $A(\beta)$-invariant linear limit of $S$, there exists a chain of linear subtriples $D_i = (\hat{G}_i, \hat {Q}_i,
 \hat{\beta}_i)$ of $S$, for  $i=0, 1, \dots, n$, such that  
$S=D_0  \geq D_1 \geq \cdots \geq D_n= S_1$, with $D_i$ being a direct linear reduction of
 $D_{i-1}$, for all $i=1, \dots, n$. Hence $D_i = D_{i-1}(\lambda_i)$ where $\lambda_i $ is a linear character of
a normal subgroup $L_i$ of $\hat{G}_{i-1}$, that lies under $\hat{\beta}_{i-1}$.  In addition,
$\hat{G}_i = \hat{G}_{i-1}(\lambda_{i})$ and $\hat{Q}_i = \hat{Q}_{i-1}(\lambda_i)$ 
while $\hat{\beta}_i \in \Irr(\hat{Q}_i)$ is the  $\lambda_i$-Clifford correspondent of
$\hat{\beta}_{i-1} \in \Irr(\hat{Q}_{i-1})$.    Thus 
$$
Z(D_{i-1}) \leq L_i \leq Z(D_i) \leq \hat{Q}_i \leq \hat{Q}_{i-1},
$$
 while $\zeta^{(D_i)}$ is an extension of  both $\lambda_i$ and $\zeta^{(D_{i-1})}$, and lies under $\hat{\beta}_i$.
According to Proposition \ref{l7}, we have $\hat{G}_i = G(\zeta^{(D_i)})$ and $\hat{Q}_i = Q(\zeta^{(D_i)})$.
 The fact that the linear limit $S_1= D_n$ is $A(\beta)$-invariant, implies
that $\zeta^{(D_i)}$ is $A(\beta)$-invariant for all $i=0, 1, \dots, n$.

According to Lemma \ref{l.e5}, applied to the groups $L_{i+1} \unlhd \hat{Q}_i \unlhd \hat{G}_i$, we have 
$G(\hat{\beta}_i) = G (\hat{\beta}_i, \lambda_{i+1}) \cdot \hat{Q}_i$, for all $i =0, 1, \dots, n-1$.
Since $G(\zeta^{(D_i)}) = \hat{G}_i \leq G(\lambda_i)$, while $\zeta^{(D_i)}$ extends $\lambda_i$, 
 we also have $G(\hat{\beta}_i, \lambda_{i+1}) = G(\hat{\beta}_i, \zeta^{(D_{i+1})})$. 
In addition,  the latter group normalizes $\hat{Q}_{i+1}= Q(\zeta^{(D_{i+1})})$ 
and  fixes its irreducible character  $\hat{\beta}_{i+1}$. Hence 
$G(\hat{\beta}_i, \zeta^{(D_{i+1})}) \leq G(\hat{\beta}_{i+1}, \zeta^{(D_{i+1})})$. The other inclusion also holds, because 
$\hat{\beta}_{i+1}$ induces $\hat{\beta}_i$ and $\hat{Q}_i = Q(\zeta^{(D_i)})$.
So for all $i = 0, 1, \dots, n-1$ we get 
$$
G(\hat{\beta}_i) =  G(\hat{\beta}_{i+1}, \zeta^{(D_{i+1})}) \cdot \hat{Q}_i
$$
Therefore, 
$$
G(\beta) = G(\hat{\beta}_1)(\zeta^{(D_1)}) \cdot Q = G(\hat{\beta}_2)(\zeta^{(D_2)},
 \zeta^{(D_1)}) \cdot \hat{Q}_1 \cdot Q= \cdots
= G(\hat{\beta}_n) (\zeta^{(D_n)}, \dots, \zeta^{(D_1)}) \cdot \hat{Q}_{n-1}\cdot  \cdots \cdot  Q.
$$
But $\zeta^{(D_n)}$ extends all the previous central characters, while $\hat{Q }_i$ is a subgroup of $Q$,
 for all $i = 1, \dots, n-1$. 
Thus $G(\beta) = G(\hat{\beta}_n, \zeta^{(D_n)}) \cdot Q$.  As  $\hat{\beta}_n$ induces $\beta$ we get 
$G(\hat{\beta}_n, \zeta^{(D_n)}) \cdot Q\leq G(\beta, \zeta^{(D_n)}) \cdot Q \leq G(\beta)$. So 
$G(\beta) = G(\beta, \zeta_1) \cdot Q$, since $\zeta_1 = \zeta^{(D_n)}$.
 This completes the proof of the first part of the corollary.        

 Because $P  \leq A(\beta) $ the characters $\zeta^{(D_i)}$ are   $P$-invariant.  As   $P$ also fixes $\beta$ we 
conclude that $P$ fixes $\hat{\beta}_i$, for all $i=1, \dots, n$.    
Let $\beta_i ^* \in \Irr(C_{\hat{Q}_i}(P))$ be the $P$-Glauberman correspondent of $\hat{\beta}_i \in \Irr(\hat{Q}_i)$,
 for all $i=1, \dots, n$.
Similarly we define $\zeta_i^* \in \Irr(C_{Z(D_i)}(P))$    to be the $P$-Glauberman correspondent of $\zeta^{(D_i)} $,
 for all such $i$. Then $\zeta_{i+1}^*$ lies under $\beta_i^*$, as $\zeta^{(D_{i+1})}$ lies under $\hat{\beta}_i$.
The group $G(\alpha, \beta, \zeta^{(D_i)})$  fixes $\zeta^{(D_0)}, \zeta^{(D_1)}, \dots,  \zeta^{(D_{i-1})}$  
since $\zeta^{(D_i)}$ 
is an  extension of all these characters.
 Because $\hat{\beta}_i$ is the only character of $\hat{Q}_i$ lying above $\zeta^{(D_i)}$ (by Proposition 
\ref{l7}), we conclude that $G(\alpha, \beta, \zeta^{(D_i)})$  fixes the Glauberman correspondents 
$\beta^*_i$, for all $i=0, 1, \dots, n$. In addition,  the group $G(\alpha, \beta, \zeta^{(D_i)})$ normalizes 
both $C_{Z(D_i)}(P)$ and $C_{\hat{Q_i}}(P)$. 
 Therefore   Lemma \ref{l.e5}  implies
\begin{equation}\mylabel{e.l6}
G(\alpha, \beta, \zeta^{(D_i)}) = G(\alpha, \beta, \zeta^{(D_i)}, \zeta_{i+1}^*) \cdot C_{Q(\zeta^{(D_i)})}(P), 
\end{equation}
for all $i=0, 1, \dots, n-1$. (For $i=0$ the above equation becomes $G(\alpha, \beta) = G(\alpha, \beta, \zeta_1^*) \cdot C$.)
The group $L_{i+1}$ was picked to be a normal subgroup of  $G_i= G(\zeta^{(D_i)}$.
 Hence $G(\alpha, \beta, \zeta^{(D_i)} ) $ normalizes $L_{i+1}$, 
as well as $P$. Thus $G(\alpha, \beta, \zeta^{(D_i)}, \zeta_{i+1}^*) = G(\alpha, \beta, \zeta^{(D_i)}, \zeta^{(D_{i+1})} )$. 
But $\zeta^{(D{i+1})}$ is an extension of $\zeta^{(D_i)}$ and therefore $G(\alpha, \beta, \zeta^{(D_i)}, \zeta^{(D_{i+1})}) = 
G(\alpha, \beta, \zeta^{(D_{i+1})} )$ for all $i=1, \dots, n-1$. So \eqref{e.l6} implies
$$
G(\alpha, \beta, \zeta^{(D_i)})   = G(\alpha, \beta, \zeta^{(D_{i+1})}) \cdot C_{Q(\zeta^{(D_i)})}(P), 
$$
 for all $i=1, \dots, n-1$, while for $i=0$ we have 
$G(\alpha, \beta) = G(\alpha, \beta, \zeta^{(D_1)}) \cdot C$
Therefore, 
$$
G(\alpha, \beta) = G(\alpha, \beta, \zeta^{(D_n)} ) \cdot C_{Q(\zeta^{(D_1)})}(P) \cdot  C_{Q(\zeta^{(D_2)})}(P) \cdot 
\cdots  \cdot C_{Q(\zeta^{(D_{n-1})})}(P) \cdot C = G(\alpha, \beta, \zeta^{(D_n)}) \cdot C.
$$
Since $\zeta^{(D_n)}= \zeta_1$,  the second part of the   corollary follows.

The group $B $ of Theorem \ref{e0} was picked so that $B(\alpha, \beta)$ is a $p'$-Hall subgroup of 
$G(\alpha, \beta)$. Since the latter group equals $G(\alpha, \beta, \zeta_1) \cdot C$ we conclude
(by looking at the indexes)  that 
$B(\alpha, \beta, \zeta_1)$ is a $p'$-Hall subgroup of $G(\alpha, \beta, \zeta_1)$. 
Of course $B(\alpha, \beta, \zeta_1) = B_1(\alpha)$ and $G(\alpha, \beta, \zeta_1) = G_1(\alpha)$, as
$G_1= G(\zeta_1)$ and it is a subgroup of $G(\beta)$.
Because $B_1(\alpha)$ is a $p'$-Hall subgroup of $G_1(\alpha)$ and $K_1  $ is a $p'$-subgroup of $B_1$ we
get that $B_1(\alpha) \cdot K_1$ is a $p'$-Hall subgroup of $G_1(\alpha) \cdot K_1$.
This completes the proof of the  corollary.
\end{proof}

We can remove the additional hypothesis on part b) of Theorem \ref{e0}
 that wants $\beta$ extendible to $Q G(\alpha, \beta)$,
 without loosing any of its 
conclusions, provided firstly that the group $G$ has order  $p^a q^b$, for
 some odd primes $p$ and $q$,   and secondly that 
we have plenty of monomial characters in $\Irr(G)$. 
 Thus we can prove the following result.
\begin{theorem}\mylabel{e1}
Assume that  $G$ is a monomial  group of order $p^a q^b$, where $p, q$ are two odd primes.
Let $L\leq M \leq N$ be  normal subgroups of $G$,  so that  $N/ M, M/ L$ and $L$ 
are nilpotent groups. So $L = L_p \times Q$ where $Q$ is the $q$-Sylow subgroup of $L$ and 
$L_p$ its $p$-Sylow subgroup.
Let  $M_p$ be a $p$-Sylow subgroup of $M$ and let  $H = M_p L$. Assume that   
$\phi \in \Irr(L)$ and $\theta  \in \Irr(H)$ lies above $\phi$ 
  and let $\theta_{\phi}  \in \Irr(H(\beta))$ 
be the $\phi$-Clifford correspondent of $\theta$. Write $\phi = \phi_p \times \beta$, for $\phi_p \in \Irr(L_p)$ 
and $\beta \in \Irr(Q)$. 
 Let $A$ be a $p$-Sylow subgroup of $G$ so 
that $A(\beta)$ is a $p$-Sylow subgroup of $G(\beta)$.  If 
 $S_1= (G_1, Q_1, \beta_1)$ is  an $A(\beta)$-invariant  
linear limit of $S= (G, Q, \beta)$, then  $(G_1 \cap H)/ K_1 = \br{H}{1} = \bar{P} \times \br{Q}{1}$,
where  $K_1$ is the kernel of $S_1$ and $\bar{P}= (PK_1)  / K_1$ while $\br{Q}{1}= Q_1/ K_1$.
Furthermore,  there exists a linear reduction
 $\br{E^*}{1}=(\br{G^*}{1}, \bar{P}^*, \bar{\alpha}^* )$ 
of $\br{E}{1}=(\br{G}{1}, \bar{P},\bar{\alpha})$
 so that $[(\br{G^*}{1} \cap \br{N}{1})_q, (\br{G^*}{1} \cap \br{H}{1})_p]  \leq  
\Ker(\br{E^*}{1})$.
\end{theorem}

\begin{proof}
The group  $H= M_p L $ is  a normal subgroup of $G$.
We fix  $\theta \in \Irr(H)$ and $\phi = \phi_p \times \beta \in \Irr(L)$ 
satisfying the hypothesis of the theorem.
 We also pick a $p$-Sylow subgroup $A$  of $G$ so that $A(\beta)$  is a $p$-Sylow 
subgroup of $G(\beta)$.

Because $H/Q$ is a $p$-group, $H$ is the semidirect product
$H = (A\cap H) \ltimes Q$, where $A \cap H$ is a $p$-Sylow subgroup of $H$.
Let $P=(A\cap H)(\phi)$. Then $H(\phi) = P \ltimes Q$.
Furthermore, $\beta $ extends to $H(\phi) \leq H(\beta)$, since $(|H/Q|, |Q|)= 1$.
Hence, if $\theta_{\phi} \in \Irr(H(\phi))$ is the $\phi$-Clifford correspondent of
$\theta$, then (see Lemma \ref{ext}) there exists a unique irreducible character $\alpha \in \Irr(P)$ such that
$$
\theta_{\phi} = \alpha \ltimes \beta,
$$
Because $H \unlhd G$, Frattini's argument implies that $G(\theta_{\phi}) = G(\alpha, \beta) Q$.
Note also that $\alpha$ lies above $\phi_p \in \Irr(L_p)$ as $\theta_{\phi}$  lies above $\phi$. 
Hence $G(\alpha) \leq G(\phi_p) \leq G$.

We   pick a $q$-Sylow subgroup $B$ of $G$ so that 
$B(\alpha)$ is a $q$-Sylow subgroup of $G(\alpha )$ while  $B(\alpha, \beta)$ is a 
$q$-Sylow subgroup of $G(\alpha,  \beta)$.

According to Theorem  \ref{e0}  and Remark \ref{inv.limit} there exists an $A(\beta)$-invariant linear
 limit $S_1= (G_1, Q_1, \beta_1)$ 
of $S= (G, Q, \beta)$ so that the  quotient group $H_1 / K_1$ is nilpotent, where 
$H_1 = G_1 \cap H$ and $K_1$ is the kernel of $S_1$.
Furthermore,  $H_1 = P \ltimes Q_1$ while 
the $S_1$-reduction $\theta_1$ of $\theta$ equals $\theta_1= \alpha \ltimes \beta_1$.
 Let $B_1= B \cap G_1$, and 
$  N_1 = N \cap G_1$. We also write $\br{G}{1}, \br{N}{1}, \br{Q}{1}, \br{B}{1}$ and 
  $\bar{P}$ for the quotient groups
 $G_1/ K_1, N_1/ K_1, Q_1/ K_1, B_1/ K_1$ and $(PK_1)/ K_1$,  respectively. 
If $\iota$ is the natural epimorphism   of $G_1$ onto $\br{G}{1}$, then 
$\iota$ sends $P$ isomorphically onto $\bar{P}$ and $\alpha \in \Irr(P)$ to some character $\bar{\alpha}
\in \Irr(\bar{P})$.
 Let $F = N_G(P)$ and $F_1 = N_{G_1} (P) = F \cap G_1$, 
then Frattini's argument implies that $G_1 = F_1 K_1$. Furthermore, $C_(Q_1)(P) = F_1 \cap Q_1$ 
covers $\br{Q}{1}$. Hence $\iota$ sends $F_1$ onto $\br{G}{1}$ with kernel 
$C_{K_1}(P)= K_1 \cap F_1$.

Let $R$ be the triple $R = (F, P, \alpha)$, where $F = N_G(P)$. We take a $B(\alpha)$-invariant linear limit 
$R^*= (F^*, P^*, \alpha^*)$ of $R$. Then Proposition \ref{l7} implies that $F^*= F(\alpha^*)= 
N_{G(\alpha^*)}(P) \leq G(\alpha)$. Because $B(\alpha)$ is a $q$-Sylow subgroup of $G(\alpha)$  
that fixes $\alpha^*$, it is also a $q$-Sylow subgroup of $F^*$.
We write $\zeta^*$ for the central character of $R^*$. We also write  $Z^*= Z(R^*)$ and $K^* = 
\Ker(R^*)$.

 Let $\br{E}{1}= (\br{G}{1}, \bar{P}, \bar{\alpha})$.
Then  any direct linear reduction $R'= (F', P', \alpha')$ of $R$ determines a direct linear reduction 
$\br{E'}{1}$ of $\br{E}{1}$ in the following way. Assume that $L \unlhd P$ is a normal subgroup $F$
  contained in $P$, 
while $\lambda \in \Lin(L)$ lies under $\alpha$ so that $P' = P(\lambda) $ and $\alpha'$ is the 
$\lambda$-Clifford correspondent of $\alpha$. Then
 $L K_1 \unlhd P K_1$, is a normal subgroup of $G_1 \leq F K_1$, 
 while the character $\lambda \ltimes 1_{K_1}$ lies under $\alpha \ltimes 1_{K_1}$. Clearly 
$\alpha' \ltimes 1_{K_1}$ is the $\lambda \ltimes 1_{K_1}$-Clifford correspondent of $\alpha \ltimes 1_{K_1}$.
Hence $\bar{L} = (L K_1) / K_1$ is a normal subgroup of $\br{G}{1}$ contained in 
$\bar{P}$, while the unique irreducible character $\bar{\lambda} \in \Irr(\bar{L})$ 
that inflates to $\lambda \ltimes 1_{K_1}$ lies under $\bar{\alpha} \in \Irr(\bar{P})$. 
 So the triple 
$\br{E'}{1}= (\br{G'}{1}, \bar{P'}, \bar{\alpha'})$
 is a direct linear reduction of $\br{E}{1}$,
 where $ \br{G'}{1}= \br{G}{1}(\bar{\lambda})= G_1(\lambda) / K_1 \leq  (F'K_1)/  K_1$.
Assume  in addition that $R'$ is $B(\alpha)$-invariant, that is $\lambda$ 
is a $B(\alpha)$-invariant 
character.  
Then the character $\lambda \ltimes 1_{K_1}$ is $B_1(\alpha \ltimes 1_{k_1})$-invariant, since 
$B_1(\alpha \ltimes 1_{k_1}) = B_1(\alpha) K_1 \leq B(\alpha) K_1$.  
Thus the character $\bar{\lambda}$ is also $\bar{B}{1}(\bar{\alpha})$-invariant. 
 Hence $R^*$ determines the linear reduction $\br{E^*}{1} =(\br{G^*}{1}, \bar{P}^*, \bar{\alpha^*})$ 
of $\br{E}{1}$, which is 
$\br{B}{1}(\bar{\alpha})$-invariant. Also 
$\br{G^*}{1} = G_1^*/ K_1$ where 
$G_1^* \leq (F^* K_1) \cap G_1$.
Furthermore, 
$\bar{P}^* $ is isomorphic to $P^*$, while  $\alpha^* \in \Irr(P^*)$ is mapped 
under that isomorphism to 
$\bar{\alpha^*} \in \Irr(\bar{P}^*)$.  
According to  Corollary  \ref{ec0}
$B_1(\alpha) \cdot K_1$ is a $q$-Sylow subgroup of $G_1(\alpha) \cdot K_1$. 
Hence $\br{B}{1}(\bar{\alpha})$ is a $q$-Sylow subgroup of $\br{G}{1}(\bar{\alpha})$. 
Because $\br{E^*}{1}$ is $\br{B}{1}(\bar{\alpha})$-invariant the latter group is also a subgroup of 
$\br{G^*}{1}$. But $\br{G^*}{1} \leq \br{G}{1}(\bar{\alpha})$, since $G_1^* \leq (F^* K_1) \cap G_1 
\leq G_1(\alpha) K_1$.
We conclude that 
\begin{equation}\mylabel{qsyl}
\br{B}{1}(\bar{\alpha}) \in \Syl_q(\br{G^*}{1}).
\end{equation}
Note also that 
\begin{equation}\mylabel{zc} 
Z(\br{E^*}{1}) \geq (Z^* \cdot K_1)/ K_1 \text{ and } \Ker(\br{E^*}{1}) \geq ( K^* \cdot K_1)/ K_1. 
\end{equation}
%Furthermore, the central character $\zeta^{E^*}$  of $E^*$ is an extension of 
%$\zeta^* \ltimes 1_{K_1}$ and thus it is an extension of $\zeta^*$ to $Z(E^*)$.

We  apply   Theorem D to the $q$-groups  $Q \unlhd B$, the $p$-group $A \cap H$,
 and the character $\beta \in \Irr(Q)$. (Clearly $B$ normalizes $(A\cap H ) \cdot Q = H$).
This way we get an irreducible character $\bn $ of $Q$  that extends to $B(\bn)$ and satisfies
$P = (A \cap H)(\beta) = (A\cap H)(\bn)$, and $B(\beta) \leq B(\bn)$. In addition, we get that
$N_B(P) \leq B(\bn)$. Hence $B(\alpha) \leq B(\bn)$.
Note also that $B(\alpha, \bn)Q = B(\alpha \ltimes 1_Q, \bn)$.
Because $B(\alpha )$ is a $q$-Sylow subgroup of $G(\alpha)$, we conclude that 
$B(\alpha, \bn) = B(\alpha)$ is a $q$-Sylow subgroup of $G(\alpha, \bn)$.
So 
\begin{equation}\mylabel{l.e7}
B(\alpha, \beta) \leq B(\alpha) = B(\alpha, \beta^{\nu}) \in \Syl_q(G(\alpha, \beta^{\nu})).
\end{equation}
Note that $\beta^{\nu}$ extends not only  to $B( \bn)$ but to $QG(\alpha, \bn)$. 
Indeed, $\bn$ extends to any $R$, where $R/ Q$ is an  $r$-Sylow subgroup of $(QG(\alpha, \bn))/Q$, for any $r \ne q$
by Corollary 8.16 in \cite{is}.  It  also extends to $QB(\alpha, \bn)$ where $(Q B(\alpha, \bn))/ Q$
 is a $q$-Sylow subgroup of $(QG(\alpha, \bn))/ Q$. Hence Corollary 11.31 in \cite{is} implies that 
$\bn$ extends to $QG(\alpha, \bn)$.

What is important about this new character is that  $H(\beta) = H(\bn) = P \ltimes Q$, 
(that is, the $p$-group $P$  remains the same for the two characters $\beta$ and $\bn$).
Furthermore, 
the product $\alpha \ltimes \beta^{\nu}$ is an irreducible character of
$H(\bn)$ lying above $\bn$. Hence Clifford's Theorem implies that $\alpha \ltimes \beta^{\nu}$
induces  an irreducible character $\theta^{\nu}$ of $H$. Also, 
$G(\theta^{\nu}_{\beta^{\nu}}) = QG(\alpha, \beta^{\nu})$.
So the groups $Q \leq H \leq G$ and the characters $\bn \in \Irr(Q)$ and
 $ \theta^{\nu} \in \Irr(H)$ 
satisfy all the hypothesis of Theorem \ref{e0}. 
Then there exists a linear limit  $S_1^{\nu}= (\gn_1, \qn_1, \bn_1)$  of $\sn= (G, Q, \bn)$ such that  
$\hn_1 = P \ltimes \qn_1$  while $\hn_1/ \kn_1$ is a nilpotent group, where $\kn_1$ is the kernel of $\sn_1$.
Let $B_1^{\nu} = B \cap G_1^{\nu}$. Then 
\begin{equation}\mylabel{bbb}
B_1^{\nu} = B_1^{\nu}(\bn_1)= B_1^{\nu}(\bn), 
\end{equation}
by Proposition \ref{l1.5}.
We also write $\br{G^{\nu}}{1}, \br{N^{\nu}}{1}$, $\br{\qn}{1}, \br{B^{\nu}}{1}$ and 
  $\bar{P}$ for the quotient groups
 $G_1^{\nu}/ K_1^{\nu}, N^{\nu}_1/ \kn_1$, $ 
\qn_1/ \kn_1, B^{\nu}_1/ \kn_1$ and $(P\kn_1)/ \kn_1$,  respectively.

 Let $\br{\en}{1}= (\br{\gn}{1}, \bar{P}, \bar{\alpha})$
Because  $P \kn_1 \unlhd \gn_1 $, Frattini's argument also implies that 
$\gn_1 =N_{\gn_1}(P) \kn_1 \leq F \kn_1$. So, as with the triple $\br{E}{1}$ and its  reduction  $\br{E^*}{1}$, 
the linear limit $R^{*} = (F^* , P^*, \alpha^*)$ determines a linear reduction 
$\br{E^{\nu^*}}{1} =(\br{G^{\nu^*}}{1}, \bar{P}^*, \bar{\alpha^*})$ 
of $\br{E^{\nu}}{1}$, which is 
$\br{B^{\nu}}{1}(\bar{\alpha})$-invariant. Note that $\bar{P}^* = (P^*K_1)/ K_1 \cong P^*$.   
Also if $\br{G^{\nu^*}}{1} = G^{\nu^*}_1 / \kn_1$ 
for some subgroup  $G^{\nu^*}_1$ of $\gn_1$, then  $\br{G^{\nu^*}}{1} \leq F^* \kn_1$ and 
in addition
\begin{equation}\mylabel{nqsyl}
\br{B^{\nu}}{1}(\bar{\alpha})  \in \Syl_q(\br{G^{\nu^*}}{1}).
\end{equation}
 Furthermore,  similarly to \eqref{zc} we get 
\begin{equation}\mylabel{zcn} 
Z(\br{E^{\nu^*}}{1}) \geq (Z^* \cdot \kn_1)/ \kn_1 \text{ and } \Ker(\br{E^{\nu^*}}{1}) \geq 
( K^{\nu^*} \cdot \kn_1)/ \kn_1. 
\end{equation}

Now let $\br{\en}{2} = (\br{\gn}{2}, \br{\pn}{2}, \br{\an}{2})$ be a 
 $\br{B^{\nu}}{1}(\bar{\alpha})$-invariant linear limit of $\br{E^{\nu^*}}{1}$. 
Then $\br{\en}{2}$ is also a  $\br{B^{\nu}}{1}(\bar{\alpha})$-invariant linear limit of $\br{\en}{1}$. 
Hence Remark \ref{inv.limit} implies that $[(\br{N^{\nu}}{1} \cap \br{\gn}{2})_q, \br{\pn}{2}] \leq
 \br{\kn}{2}$ is nilpotent, where $\br{\kn}{2}$ is the kernel of $\br{\en}{2}$.
(Note that we have used the fact that 
 the group $B$ is a $q$-Sylow subgroup of 
$G$ so that $B(\alpha)$ is a $q$-Sylow subgroup of $G(\alpha)$ while $B(\alpha, \beta^{\nu})$ is a
$q$-Sylow subgroup of $G(\alpha, \beta^{\nu})$.) 
By \eqref{nqsyl} the group $\br{B^{\nu}}{1}(\bar{\alpha})$ is a $q$-Sylow subgroup of $\br{G^{\nu^*}}{1}$. 
Because  $\br{\en}{2}$ is  a 
 $\br{B^{\nu}}{1}(\bar{\alpha})$-invariant linear limit of $\br{E^{\nu^*}}{1}$, we get that  
$\br{B^{\nu}}{1}(\bar{\alpha})$ is contained in $\br{\gn}{2} \leq  \br{G^{\nu^*}}{1} $, and thus it is a
 $q$-Sylow subgroup of $\br{\gn}{2}$.
Hence 
\begin{equation}\mylabel{e.l12}
[\br{B^{\nu}}{1}(\bar{\alpha}) \cap \br{N}{1}, \br{\pn}{2}] \leq \br{\kn}{2}, 
\end{equation}
where $\br{\kn}{2}$ is the kernel of $\br{\en}{2}$.
In addition, \eqref{zcn} along with Remark \ref{l00} implies 
\begin{subequations}\mylabel{zp} 
\begin{equation}
Z^* \cong (Z^* K_1^{\nu})/ \kn_1 \leq Z(\br{E^{\nu^*}}{1}) \leq Z(\br{\en}{2})
\leq \br{\pn}{2} \leq \bar{P}^* \cong P^* 
\end{equation} 
and 
\begin{equation}\mylabel{zpb}
(K^* \kn_1)/ \kn_1 \leq \Ker(\br{E^{\nu^*}}{1}) \leq \br{\kn}{2}.
\end{equation}
\end{subequations}
We identify $Z^*$  with its isomorphic image $\bar{Z^*}= (Z^* \kn_1)/ \kn_1$ inside 
$Z(\br{E^{\nu^*}}{1})$. Under this isomorphism the central character $\zeta^* \in \Lin(Z^*)$
 of $R^*$ is mapped to a linear character of $\bar{Z^*}$ that lies under the cental 
character $\zeta^{\nu^*}_1 \in \Lin(Z(\br{E^{\nu^*}}{1}))$ of $\br{E^{\nu^*}}{1}$.

Let $V := P^*/ Z^*= P^* / Z(R^*)$. 
 Then $V$ is an anisotropic $F^*/P^*$-group, by  Proposition \ref{aniso}.
(The $F^*/ P^*$-invariant bilinear form 
$c : V \times V \to C ^{\times } $  is defined (see \eqref{form} ) 
as $c(\bar{x}, \bar{y}) = \zeta^*([x, y])$, for all $\bar{x}, \bar{ y}  \in V$).
Thus $V$ written additively, is the direct sum
$$
V = V_1 \dotplus V_2,
$$
of the perpendicular $F^*/P^*$-groups $V_1= C_V(N^*)$ and $V_2 = [V, N^*]$, 
where $N^* = N \cap F^*$. 
Because  $B(\alpha)$ is a $q$-Sylow subgroup of $F^*$, the  
 group $Q^*= B(\alpha) \cap N^*= B(\alpha) \cap N$  is a $q$-Sylow subgroup of $N^*$. So
$$
N^* = Q^*\ltimes P^*.
 $$
  Therefore the direct summands $V_1, V_2$ of $V$ are
\begin{equation}\mylabel{e.l13}
V_1 = C_V(Q^*) \text{ and }  V_2 = [V, Q^*].
\end{equation}
Both $V_1$ and $V_2$ are anisotropic $F^*/ P^*$-groups.

Now let $U = \br{\pn}{2} /Z(\br{E^{\nu}}{2})$.
 Then $U$ is isomorphic to a section   of $V$, by \eqref{zp}.
Furthermore, $U$ is isomorphic
to the orthogonal direct sum $U =U_1 \dotplus U_2$, where
in view of \eqref{e.l13}   we get
$$
U_1= C_U(Q^*) \text{ and } U_2 = [U,Q^* ], 
$$
According to 
Corollary \ref{ec0} (applied to the $\nu$-groups), 
we have
$G_1^{\nu}(\alpha, \bn) C_Q(P) = G(\alpha, \bn)$. This, along with 
\eqref{bbb} and  \eqref{l.e7} implies 
\begin{equation}\mylabel{b1cb}
B_1^{\nu}(\alpha) \cdot C_Q(P)=  B_1^{\nu}(\alpha, \beta^{\nu}) \cdot C_Q(P) =
B(\alpha, \bn)= B(\alpha).
\end{equation}
Hence the image of $B(\alpha ) \cap N$ in the automorphism group of $P$ equals that of 
$ B_1^{\nu}(\alpha) \cap N$. So 
$U_2= [U, B_1^{\nu}(\alpha) \cap N]$. In view of  
\eqref{e.l12} and \eqref{zpb},   this latter group is trivial.
Hence  $U_2 = 0$ and $U_1= C_U(Q^*) =C_U(B_1^{\nu}(\alpha, \bn) \cap N)$.

Because  $V=P^*/ Z^*$ 
is an abelian anisotropic $F^*/ P^*$-group, the 
group  $\bar{P}^*/ Z(\br{E^{\nu^*}}{1})$  is also an abelian group.
It also affords a bilinear  $\br{G^{\nu^*}}{1}/ \bar{P}^*$-invariant 
form (see \eqref{form})) defined as 
 $\hat{c}  (\bar{u}, \bar{v}) = \zeta^{\nu^*}_1([u, v])$, for all
 $\bar{u}, \bar{ v}  \in \bar{P}^* / Z(\br{E^{\nu^*}}{1})$.
Identifying $P^*$ with $\bar{P}^*$ and $Z^*$ with $ \bar{Z^*}$ we see that 
 $[\bar{P}^*, \bar{P}^* ] \leq \bar{Z^*}  \leq Z(\br{E^{\nu^*}}{1})$. 
Hence
$\hat{c} (\bar{u}, \bar{v}) = \zeta^* ([u, v])$, 
for all $\bar{u}, \bar{ v}  \in \bar{ P^*}/ Z(\br{E^{\nu^*}}{1})$. 
In addition, since  
$G^{\nu^*}_{1} \leq F^* K_1^{\nu}$ the factor group $\br{G^{\nu^*}}{1} / \bar{P}^*  $ 
is isomorphic to a subgroup of $F^* / P^*$. Hence  $\bar{P}^*/ Z(\br{E^{\nu^*}}{1})$ 
is an abelian symplectic $\br{G^{\nu^*}}{1} / \bar{P^* } $. Thus  
  we can apply Proposition \ref{l6} to the linear limit $\br{E^{\nu}}{2}$ of
 $\br{E^{\nu^*}}{1}$.
So $U$  is isomorphic to $L^{\perp}/ L$, where $L$ is 
 maximal among the  $\br{G^{\nu^*}}{1} / \bar{P}^*$-invariant  isotropic  subgroups  of  
$\bar{P}^* / Z(\br{E^{\nu^*}}{1})$.
Furthermore, $L = Z(\br{E^{\nu}}{2})/ Z(\br{E^{\nu^*}}{1})$ and $L^{\perp} 
= \br{P^{\nu}}{2}/ Z(\br{E^{\nu^*}}{1})$.
If  $L = L_1 \dotplus L_2$, with $L_1 = C_L(Q^*)$ and $L_2= [L, Q^*]$ then the facts that 
$ L^{\perp} / L= U$ while $U_2 = 0$ implies that $L_2^{\perp} = L$.
  Hence $L_2$ is a self perpendicular  
$\br{G^{\nu^*}}{1}/ \bar{P}^*$-invariant subgroup of $\bar{P}^*/ Z(\br{E^{\nu^*}}{1}) $.

Since $\bar{P}^*/ Z(\br{E^{\nu^*}}{1}) \cong \frac{P^*/ Z^*} {Z(\br{E^{\nu^*}}{1})/ Z^*}$, 
 we get that   $L_2$ is isomorphic to a self perpendicular subgroup $W_2$ of $V$. Hence
$W_2 \leq V_2 = [V, Q^*]$, as $L_2 = [L, Q^*]$. Because $L_2$ is 
$\br{G^{\nu^*}}{1}/ \bar{P}^*$-invariant, it is $\br{B^{\nu}}{1}(\bar{\alpha})$-invariant, by 
\eqref{nqsyl}).
According to \eqref{b1cb} the  image of $B^{\nu}_1(\alpha) $ in $\Aut(P)$ 
equals that of $B(\alpha)$.
Therefore $W_2$ is a self perpendicular $B(\alpha )$-invariant subgroup of $V_2$. 
Hence  $V_2$ is 
a hyperbolic $B(\alpha)$-group. 
But $B(\alpha)$ is a $q$-Sylow subgroup of $F^*$. Hence Theorem 3.2 in \cite{da} implies that 
$V_2$ is a hyperbolic $F^*$-group. 
Because  it is also an anisotropic $F^*/ P^*$-group, we conclude that $V_2 $ is $0$. 
Hence $Q^*$ centralizes $V = P^*/ Z^*$. 
In addition, $Z^*/ K^*$ is a cyclic central $p$-section of $F^*$.
So the $q$-subgroup $Q^*$ of $F^*$ centralizes both $p$-groups  $V = P^*/ Z^*$ and $Z^*/ K^*$. 
We conclude  that $Q^*$ centralizes $P^*/ K^*$. Since $Q^* = B(\alpha) \cap N$ we get  
$[P^* , B(\alpha) \cap N] \leq K^*$. Thus 
$[(P^*K_1)/ K_1 , (B(\alpha)K_1 \cap N)/ K_1] \leq (K^*K_1)/ K_1$. This and \eqref{zc} implies 
that 
$$
[\bar{P}^* , \br{B}{1}(\bar{\alpha}) \cap \br{N}{1}] \leq  \Ker(\br{E^*}{1}).
$$  
But  $\br{B}{1}(\bar{\alpha})  \in \Syl_q(\br{G^*}{1})$, by \eqref{qsyl}.
We conclude that a $q$-Sylow subgroup of $\br{N^*}{1}$ centralizes $\bar{P}^*$ module
$\Ker(\br{E^*}{1})$. Hence the theorem follows.
\end{proof}

Note that  we pick the linear limit $\br{E^*}{1}$ in the statement of 
Theorem \ref{e1} in the following  way
\begin{corollary}\mylabel{fco}
Assume the hypothesis and notation of Theorem \ref{e1}. 
Let $B$ be a $q$-Sylow subgroup
of $G$ so that $B(\alpha) $ and $B(\alpha, \beta)$ are $q$-Sylow subgroups of $G(\alpha)$
 and $G(\alpha, \beta)$, respectively. Let $F= N_G(P)$. Then any 
$B(\alpha)$-invariant linear limit $(F^*, P^*, \alpha^*)$ 
of $(F, P, \alpha)$ determines naturally a linear reduction 
$\br{E^*}{1}=(\br{G^*}{1}, \bar{P}^*, \bar{\alpha^*} )$ 
of $\br{E}{1}=(\br{G}{1}, \bar{P},\bar{\alpha})$
 so that $[(\br{G^*}{1} \cap \br{N}{1})_q, \bar{P}^*]  \leq  
\Ker(\br{E^*}{1})$, while $\bar{P}^*\bar{P}^*$ is a $p$-Sylow subgroup of $\br{H^*}{1}$.
\end{corollary}

We can now prove Theorem E that we restate here. 
\begin{E}
Assume that $G$ is a finite odd monomial $p, q$-group, while $N$ is a normal subgroup of $G$
of nilpotent length $3$. So there exists  $L \leq M $  normal subgroups 
of $G$ contained in $N$  so that 
$L, M/L$ as well as $N/ M$ are nilpotent groups.
Let $\chi$ be any  irreducible character of $M$. 
 Then there exists a linear limit $T'= (G', M', \chi')$ 
of the triple $T= (G, M, \chi)$ with  
$N'/ \Ker(T')$ being nilpotent, where $N' = N \cap G'$, and 
$\Ker(T')$   is the kernel   of the triple $T'$.
Therefore $N$ is a monomial group.
\end{E}

\begin{proof}
Let $\phi \in \Irr(L)$ lying under $\chi \in \Irr(M)$. 
Because $L = L_p \times L_q$ is nilpotent,  $\phi$ splits us 
$\phi= \eta \times \beta$, where $\eta \in \Irr(L_p)$ and $\beta \in \Irr(L_q)$.  
The groups $H= M_p \ltimes L_q$  and $J= M_q \ltimes L_p$ are  normal subgroups of $G$. 
Furthermore, $H(\phi) = P \ltimes L_q$ where $P = M_p(\phi)$ is a $p$-Sylow subgroup of
 $H(\phi)$, 
and $J(\phi)= O \ltimes L_p$ where $O= M_q(\phi)$ is a $q$-Sylow subgroup of $J(\phi)$.
Clearly $P \geq L_p$ and $O \geq L_q$.

Let $\theta \in \Irr(H)$ be any irreducible character of $H$ lying above $\phi$ and under $\chi$.
Then the $\phi$-Clifford correspondent $\theta_{\phi}$ of $\theta $ equals 
$$
\theta_{\phi}= \alpha \ltimes  \beta,  
$$
where $\alpha \in \Irr(P)$  is uniquely determined by $\theta$ 
and lies above $\eta \in \Irr(L_p)$. Furthermore, $\alpha$ restricted to 
$L_p$ is a multiple of $\eta$.
Similarly we pick an irreducible character $\psi \in \Irr(J)$ lying above $\phi$ and under 
$\chi$. Its $\phi$-Clifford correspondent satisfies
$$ 
\psi_{\phi} =  \gamma \ltimes \eta, 
$$
where $\gamma \in \Irr(O)$ is uniquely determined by $\psi$ and 
lies above $\beta \in \Irr(L_q)$. In addition,  $\gamma$ restricted to $L_q$ is a multiple of 
$\eta$. 
It is easy to see that 
\begin{subequations}\mylabel{abpq}
\begin{equation}\mylabel{abq}
G(\alpha, \gamma) \leq G(\alpha, \beta) \leq G(\alpha) \leq G(\eta)
\end{equation}
and 
\begin{equation}\mylabel{abp}
G(\alpha, \gamma) \leq G(\eta, \gamma) \leq G(\gamma) \leq G(\beta).
\end{equation}
\end{subequations}

Now we pick a $p$-sylow subgroup $A$ and a $q$-Sylow subgroup $B$ of $G$ so that 
$A$ intersected with every group in \eqref{abp} is a $p$-Sylow subgroup of that group, 
and $B$ intersected with any group in \eqref{abq} is a $q$-Sylow subgroup of that group.

Let $(G', L_q', \beta')$ be an $A(\beta)$-invariant linear limit of 
$(G, L_q, \beta)$. Let also $(G'', L_p'', \eta'')$ be a $B(\eta)$-invariant 
linear limit of $(G, L_p, \eta)$. So $H' = H\cap G' = P \ltimes L'$, and $J' = 
J \cap G'' = O \ltimes L_p''$. Also $\theta' = \alpha \ltimes \beta' \in \Irr(H')$ 
is the $G'$-reduction of $\theta$, while $\psi'' = \beta \ltimes \eta'' \in \Irr(J'')$   
is the $G''$-reduction of $\psi$. 
If $G_1 = G' \cap G''$, $L_1 = L_p'' \times L_q'$ 
and $\phi_1 = \eta'' \times \beta'$  then 
 $(G_1, L_1, \phi_1)$  is a linear limit of $(G, L, \phi)$, by Remark \ref{l2}. 
In addition, (see  Lemma \ref{C.l2})  the factor group  $M_1/ K_1$ is a nilpotent
 group where $K_1$ is the kernel of $(G_1, L_1, \phi_{1})$, and $M_1 = M \cap G_1$. 
Note that $K_1 = K'' \times K'$ where $K''$ is the kernel of 
$(G'', L_p'', \eta'')$
and $K'$ is the kernel of $(G', L_q', \beta')$.
In addition,  $M_1/ K_1 = (P_1K_1) /K_1 \times (O_1K_1)/ K_1$, where 
$P_1 = P \cap G_1= P \cap G''$, (where the last equality follows from the fact that 
the $A(\beta)$-invariant reductions of $(G, L_q, \beta)$ do not change  $P \leq A(\beta)$)  and 
similarly 
$O_1 = O \cap G_1 = O \cap G'$.
Let $\alpha_1 \in \Irr(P_1)$ and $\gamma_1 \in \Irr(O_1)$ be the
 $G_1$-reductions of $\alpha 
\in \Irr(P)$ and $\gamma \in \Irr(O)$, respectively.
So $\alpha_1$ is actually the $(G'', L_p'', \eta'')$-reduction of $\alpha$ while
$\gamma_1$ is the $(G', L_q', \beta')$-reduction of $\gamma$.
Note also that $(P_1 K_1)/ K_1 = (P_1 K')/ K_1$ while $(O_1 K_1)/ K_1 = (O_1 K'')/ K_1$.
Let $\bar{\alpha}_1 \in \Irr((P_1 K_1)/K_1)$ be the unique 
irreducible character of $(P_1 K_1)/K_1$ that inflates $\alpha \ltimes 1 \in \Irr(P_1 K')$. 
Similarly we define $\bar{\gamma}_1 \in \Irr((O_1K_1)/ K_1)$.

Now let $R= (F, P, \alpha)$ and $S = (I, O, \gamma)$, where $F = N_G(P)$ and $I = N_G(O)$.
Observe that the $B(\eta)$-invariant linear limit $(G'', L_p'', \eta'')$ 
of $(G, L_p, \eta)$ determines naturally a  $B(\eta)$-invariant  linear reduction 
$R'' =(F \cap G'', P \cap G'', \alpha'')$ of $R$, where 
$\alpha''$ is the $(G'', L_p'', \eta'')$-reduction of $\alpha$.  
But $P \cap G'' = P_1$ and thus $R''= (F'', P_1, \alpha_1)$.
Note also that in view of \eqref{abq} the group $B(\alpha)$ is a subgroup of  $B(\eta)$
 Hence the linear 
reduction $R''$ of $R$ is $B(\alpha)$-invariant. 

Similarly, 
 the $A(\beta)$-invariant linear limit $(G', L_q', \beta')$ 
of $(G, L_q, \beta)$ determines naturally an  $A(\beta)$-invariant  linear reduction 
$S' =(I', O_1, \gamma_1)$ of $S$, where $ I' =G' \cap I$.
 In addition,  \eqref{abp} implies that $A(\gamma) \leq A(\beta)$. Hence the linear 
reduction $S'$ of $S$ is $A(\gamma)$-invariant. 

Let $R^*= (F^*, P^*, \alpha^*)$ be a $B(\alpha)$-linear limit of $R''= (F'', P_1, \alpha_1)$ and thus of $R$. 
Similarly we pick $S^*=(I^*, O^*, \gamma^*)$ to be an  $A(\gamma)$-invariant linear limit of 
$S'=(I', O_1, \gamma_1)$ and thus of $S$.

Now we can apply Corollary \ref{fco}. Let $ \bar{G}''= G''/ K''$ then $(PK'')/ K''= \bar{P}$ 
is isomorphic to   $P$ and under this isomorphism $\alpha \in \Irr(P)$ gets mapped to 
$\bar{\alpha} \in \Irr(\bar{P})$. If 
$\bar{E} = (\bar{G}'',  \bar{P}, \bar{\alpha})$, then Corollary \ref{fco} implies that  
 $R^*$ determines naturally a linear reduction 
$\bar{E}^*= (\bar{G}^*, \bar{P}^*, \bar{\alpha}^*)$ so that 
\begin{equation}\mylabel{fee1}
[(\bar{G}^* \cap \bar{N}^*)_p, \bar{P}^* ] \leq \Ker(\bar{E}^*),
\end{equation}
where $\bar{P}^*$ is a $p$-Sylow subgroup of $\bar{H}^*= \bar{G}^* \cap (H''/ K'')$. 
Observe that because $R^*$ is a linear limit of $R''$, the reductions done in $\bar{E}$ 
are actually reductions done inside $(P_1K'')/ K'' \leq \bar{P}$.
Because $G_1 = G' \cap G''$ and $K_1 = K' \times K''$,  the linear reduction $\bar{E^*}$ of $\bar{E}$ 
lifts to a linear reduction $U_1$
 of  $U=(G_1/ K_1, (P_1 K_1)/ K_1, \bar{\alpha}_1 )$.
 Note that $G_1/K_1$ is reduced to $G_1/K_1 \cap \bar{G}^* K_1/K_1$.

Similarly we write $\hat{G}'= G'/ K'$ then $(OK')/ K'= \hat{O}$ 
is isomorphic to   $O$ and under this isomorphism $\gamma \in \Irr(O)$ get mapped to 
$\hat{\gamma} \in \Irr(\hat{O})$. If 
$\hat{D} = (\hat{G}',  \hat{O}, \hat{\gamma})$, then Corollary \ref{fco} implies that  
 $S^*$ determines naturally a linear reduction 
$\hat{D}^*= (\hat{G}^*, \hat{O}^*, \hat{\gamma}^*)$ so that 
\begin{equation}\mylabel{fee2}
[(\hat{G}^* \cap \hat{N}^*)_q, \hat{O}^* ] \leq \Ker(\hat{D}^*),
\end{equation}
where $\hat{Q}^*$ is a $q$-Sylow subgroup of $\hat{J}^*= \hat{G}^* \cap (J'/ K')$.
Furthermore, similarly to $U_1$ and $U$ we get that the linear reduction 
$\bar{D}^*$ of $\bar{D}$ lifts to a linear reduction 
$V_1$ of $V = (G_1/ K_1, (O_1K_1)/ K_1, \bar{\gamma}_1)$. 
Also the group $G_1/K_1$ is being reduced to 
$G_1/ K_1 \cap (\hat{G}^*K_1)/ K_1$.

The fact that $M_1/ K_1 = (P_1 K_1)/ K_1 \times (O_1 K_1)/ K_1$ is a nilpotent group, 
along with Remark \ref{l2} provides a linear  reduction $\bar{T}_1$  of 
$\bar{T}=(G_1/ K_1, M_1/ K_1, \bar{\alpha}_1 \times \bar{\gamma}_1)$ using the linear reductions 
$U_1$ and $V_1$. The kernel of $\bar{T}_1$ contains  the group $(\Ker(\bar{E}^*)K_1)/K_1 \times (\Ker(\hat{D}^*) K_1)/K_1$.
So \eqref{fee1} and \eqref{fee2} imply that the $\bar{T}_1$-reduction of the group  $N_1/K_1$ is a nilpotent group 
module the kernel of $\bar{T}_1$.

Because $\bar{T}_1$ lifts to a reduction of the triple $(G_1, M_1, \chi_1)$  and thus of $(G, M, \chi)$, the theorem follows.

\end{proof}


\begin{thebibliography}{9}

\bibitem{be}T. R. Berger, Hall--Higman Type Theorems V.Pacific J.Math.  \textbf{73}, 1-62 (1977)

\bibitem{da} E. C. Dade,  Monomial Characters and Normal Subgroups,
Math.Z. \textbf{178}, 401--420 (1981)

\bibitem{da1} E. C. Dade, Normal subgroups of $M$-groups need not be
$M$-groups. Math. Z. \textbf{133}, 313--317 (1973)

\bibitem{dl}E. C. Dade, M.Loukaki,  Linear limits of irreducible characters, preprint 

\bibitem{do} L. Dornhoff, $M$-groups and $2$-groups, Math. Z. \textbf{100}, 226--256 (1967)


%\bibitem{gl} G. Glauberman, Correspondence of characters for relatively prime operator groups,
%Can. J. Math. \textbf{20}, 1465--1488 (1968)

\bibitem{go} D. Gorenstein, Finite Groups. New York: Harper and Row 1968

\bibitem{is} I. M. Isaacs, Character Theory of Finite Groups. New York-San
Francisco--London: Academic Press 1976


\bibitem{isa3} I. M. Isaacs,  Primitive characters, normal subgroups, and $M$-groups, 
Math. Z. \textbf{177}, no. 2, 267--284, (1981)

\bibitem{isa4} I.M.Isaacs,  Induction and restriction of $\pi$-special 
characters. Can. J. Math. \textbf{38}, 576--604 (1986)

\bibitem{isa5} I. M. Isaacs,  Induction and restriction of $\pi$-partial 
characters and their lifts. Can. J. Math. \textbf{48}, no. 6,  1210--1223 (1996)

\bibitem{isol} I. M. Isaacs,  Characters of Solvable and Symplectic Groups, Amer. J. Math.
 \textbf{95}, 594--635, (1973)


\bibitem{lew} M. L. Lewis, Characters, coprime actions, and operator groups, Arch.Math.
\textbf{69}, 455--460 (1997)

\bibitem{lew2} M. L. Lewis,  Primitive characters of subgroups of $M$-groups, 
Proc. Amer. Math. Soc.  \textbf{125}, no. 1, 27--33,  (1997)

\bibitem{lo} M. Loukaki, Normal subgroups of odd order monomial $p^a q^b$-groups, 
Thesis, University of Illinois at Urbana--Champaign, 2001


\bibitem{wo}O. Manz, T. R. Wolf,  Representations of Solvable Groups, 
Lecture Notes Series \textbf{185} 

\bibitem{na}G. Navarro, Primitive characters of subgroups of $M$-groups, Math.Z.
\textbf{218}, no. 3,  439--445 (1995)

\bibitem{pa}A. Parks,  Nilpotent by supersolvable $M$-groups,  Canad. J. Math. 
\textbf{37}, no. 5,  934--962 (1985)

\bibitem{se} G. Seitz, $M$-groups and the supersolvable residual, Math. Z.
   \textbf{110}, 101--122  (1969)


\bibitem{wa}R.W.van der Waall, On the embedding of minimal non-$M$-groups.
Indag. Math. \textbf{36}, 157--167 (1974)

\end{thebibliography}
\end{document}